\documentclass[onefignum,onetabnum]{siamonline190516}
\pdfoutput=1 

\usepackage{xfrac}  
\usepackage{lipsum}
\usepackage{amsfonts,dsfont,mathrsfs,mathtools,mathabx,extarrows}
\usepackage{esvect}
\usepackage{graphicx}
\usepackage{epstopdf}
\usepackage{algorithmic}
\usepackage{tikz}
\usepackage{graphicx}

\usepackage{enumitem,comment}
\usepackage{subcaption}
\usepackage{graphicx}
\ifpdf
\hypersetup{
  pdftitle={Adjoint},
  pdfauthor={Tan Bui-Thanh}
}
\usepackage{multirow}
\usepackage{colortbl}
\fi
\usetikzlibrary{positioning}

\ifpdf
  \DeclareGraphicsExtensions{.eps,.pdf,.png,.jpg}
\else
  \DeclareGraphicsExtensions{.eps}
\fi
\usepackage{comment}

\usepackage{tikz}
\usetikzlibrary{arrows,shapes}
\usetikzlibrary{decorations.pathreplacing,calligraphy}
\usetikzlibrary{patterns}
\usetikzlibrary{arrows,shapes}
\usetikzlibrary{decorations.pathreplacing,calligraphy}

\newsiamremark{remark}{Remark}
\newsiamremark{hypothesis}{Hypothesis}
\crefname{hypothesis}{Hypothesis}{Hypotheses}
\newsiamthm{claim}{Claim}
\newsiamremark{example}{Example}
\newsiamremark{convention}{Convention}

\headers{Adjoint operators in computational Sciences, Engineering, and Mathematics}{Tan Bui-Thanh}

\title{Adjoint  and Its roles in Sciences, Engineering, and Mathematics: A Tutorial
}

\author{Tan Bui-Thanh\thanks{Department of Aerospace Engineering and
  Engineering Mechanics, The Oden Institute for Computational
  Engineering and Sciences,  the University of Texas at Austin, Texas
  (\email{tanbui@oden.utexas.edu},
  \url{https://users.oden.utexas.edu/\string~tanbui/}).}
}
\usepackage{amsopn}

\usepackage{xr}
\makeatletter
\newcommand*{\addFileDependency}[1]{
  \typeout{(#1)}
  \@addtofilelist{#1}
  \IfFileExists{#1}{}{\typeout{No file #1.}}
}
\makeatother

\newcommand{\Grad} {\ensuremath{\nabla}}  
\newcommand{\Div} {\ensuremath{\nabla\cdot}} 

\newcommand{\nor}[1]{\left\| #1 \right\|} 
\newcommand{\snor}[1]{\left| #1 \right|} 
\newcommand{\LRp}[1]{\left( #1 \right)} 
\newcommand{\LRs}[1]{\left[ #1 \right]} 
\newcommand{\LRa}[1]{\left< #1 \right>} 
\newcommand{\LRc}[1]{\left\{ #1 \right\}} 
\newcommand{\pp}[2]{\frac{\partial #1}{\partial #2}} 

\newcommand{\mc}[1]{\mathcal{#1}} 
\newcommand{\mb}[1]{\mathbf{#1}} 
\newcommand{\mbb}[1]{\mathbb{#1}} 
\newcommand{\ms}[1]{\mathscr{#1}} 
\newcommand{\half}{\frac{1}{2}}

\newcommand{\bs}[1]{\boldsymbol{#1}}

\newcommand{\X}{X}
\newcommand{\Xb}{{\bs{E}}}
\newcommand{\e}{e}
\newcommand{\Y}{Y}
\newcommand{\Yb}{{\bs{G}}}
\newcommand{\g}{g}
\newcommand{\f}{f}
\newcommand{\fb}{\bs{\f}}
\newcommand{\Z}{Z}
\newcommand{\U}{U}
\newcommand{\V}{V}

\newcommand{\B}{{\ms{B}}}
\newcommand{\Bf}{{\mathfrak{B}}}
\newcommand{\Xs}{{\mbb{\X}}}
\newcommand{\Ys}{{\mbb{\Y}}}
\newcommand{\Us}{{\mbb{U}}}
\newcommand{\Vs}{{\mbb{V}}}

\newcommand{\Zs}{{\mbb{\Z}}}
\newcommand{\Fs}{{\mbb{F}}}
\newcommand{\Ls}{{\mbb{L}}}
\newcommand{\Lms}{{\ms{L}}}
\newcommand{\Hs}{{\mbb{H}}}
\newcommand{\Ps}{{\mbb{P}}}
\newcommand{\Ss}{{\mbb{S}}}
\newcommand{\R}{{\mbb{R}}}
\newcommand{\real}{{\R}}
\newcommand{\Gs}{{\mathsf{G}}}
\newcommand{\C}{{\mbb{C}}}
\newcommand{\Cs}{{\mathfrak{C}}}
\renewcommand{\L}{{\mathfrak{L}}}
\newcommand{\ub}{\bs{u}}
\newcommand{\A}{{\ms{A}}}
\newcommand{\Am}{\mc{A}}
\newcommand{\Bm}{\mc{B}}
\newcommand{\Cm}{\mc{C}}
\newcommand{\Mm}{\mc{M}}
\newcommand{\Dm}{\mc{D}}
\newcommand{\Fm}{\mc{F}}
\newcommand{\Hm}{\mc{H}}
\newcommand{\Zm}{\mc{Z}}
\newcommand{\Um}{{\mc{U}}}
\newcommand{\Uset}{{\mathsf{U}}}
\newcommand{\Vm}{{\mc{V}}}
\newcommand{\Wm}{{\mc{W}}}
\renewcommand{\P}{{\ms{P}}}
\renewcommand{\Im}{{I}}
\newcommand{\I}{\ms{I}}
\newcommand{\D}{{\mathscr{D}}}
\newcommand{\Do}{{\mathsf{D}}}
\newcommand{\x}{{x}}
\newcommand{\y}{{y}}
\newcommand{\z}{{z}}
\newcommand{\zb}{\bs{\z}}
\newcommand{\xb}{\bs{\x}}
\newcommand{\yb}{\bs{\y}}
\newcommand{\cb}{\bs{c}}
\renewcommand{\u}{{u}}
\renewcommand{\v}{{v}}
\newcommand{\w}{{w}}
\newcommand{\vb}{\bs{\v}}
\newcommand{\wb}{\bs{\w}}
\newcommand{\defeq}{:=}
\newcommand{\Range}{{\mathsf{R}}}
\newcommand{\N}{\mathsf{N}} 
\newcommand{\zev}{\theta}
\newcommand{\n}{n}
\newcommand{\nb}{\bs{n}}
\newcommand{\m}{m}
\renewcommand{\d}{{d}}
\renewcommand{\vec}{\bs}
\newcommand{\vareps}{{\varepsilon}}
\newcommand{\h}{{h}}
\newcommand{\hb}{{\bs{\h}}}
\newcommand{\bb}{{\bs{b}}}
\newcommand{\rr}{\z}
\newcommand{\rrb}{{\mb{\rr}}}
\renewcommand{\a}{{a}}
\newcommand{\ab}{{\bs{\a}}}
\newcommand{\pOmega}{\partial\Omega}

\newcommand{\verteq}[1]{\rotatebox{90}{$\xlongequal[\rotatebox{-90}{#1}]{}$}}

\newcommand*{\vertbar}{\rule[-1ex]{0.5pt}{2.5ex}}

\newcommand{\dd}[2]{\frac{d #1}{d #2}} 

\newcommand{\Ball}[2]{B_{#1}\LRp{#2}}
\newcommand{\equaldef}{{:=}}

\newcommand{\SubjectTo}{\hspace{0.2cm} \mbox{subject} \hspace{0.15cm} \mbox{to} \hspace{0.2cm}}

\usepackage{hyperref}
\hypersetup{
    colorlinks=true,
    linkcolor=blue,
    filecolor=magenta,      
    urlcolor=red,
    pdftitle={adjoint},
    pdfpagemode=FullScreen,
    }

\newcommand{\eqnlab}[1]{\label{eq:#1}}
\newcommand{\theolab}[1]{\label{theo:#1}}

\newcommand{\propolab}[1]{\label{propo:#1}}

\newcommand{\exalab}[1]{\label{exa:#1}}

\newcommand{\lemref}[1]{\ref{lem:#1}}
\newcommand{\eqnref}[1]{\eqref{eq:#1}}

\newcommand{\eval}[2][\right]{\relax \ifx#1\right\relax \left.\fi#2#1\rvert}
\renewcommand{\epsilon}{\varepsilon}

\def\height{2.5cm}
\def\width{1.8cm}
\def\angle{45}

\begin{document}

\maketitle


\begin{abstract}
    This paper synergizes the roles of adjoint in  various disciplines of mathematics, sciences, and engineering. Though the materials developed and presented are  not new\textemdash as each or some could be found in (or inferred from) publications in different fields\textemdash  we believe this is the first effort to systematically unify these materials on the same mathematical footing starting from the basic definitions. We aim to provide a unified perspective and understanding of adjoint applications.
    As a result, 
this work could give broader views and better insights into the application of adjoint beyond a single community. By rigorously establishing general results and then developing materials specific to each application, we bring forth  details on  how abstract concepts/definitions can be translated into particular applications and the connections among them. This paper is written as an interdisciplinary tutorial on adjoint with discussions and with many examples from different fields including linear algebra (e.g. eigendecomposition and the singular value decomposition), ordinary differential equations (asymptotic stability of an epidemic model), partial differential equations (well-posedness of elliptic, hyperbolic, and Friedrichs' types), neural networks (backpropagation of feed-forward deep neural networks), least squares and inverse problems (with Tikhonov regularization), and PDE-constrained optimization. The exposition covers results and applications in both finite-dimensional and infinite-dimensional Hilbert spaces.
\end{abstract}


\section{Introduction}
Adjoint is ubiquitous in mathematics. According to \cite{HistoryOfDiff83}, the history of adjoint can be traced back to as far as Lagrange in 1766 \cite{lagrangesolution}, in a  memoir extending the letter that he wrote to D'Alembert in January 1765 discussing his new method of solving $n$th-order differential equations.
Also according to \cite{HistoryOfDiff83}, the term ``adjoint equation", first used to call the corresponding equation developed from Lagrange memoir \cite{lagrangesolution}, is due to Fuchs \cite{Fuchs1866}. The adjoint operator, as reported in \cite{Lindstrm2008OnTO}, was introduced by Riesz in his seminal paper\footnote{As also argued in \cite{Lindstrm2008OnTO}, in his seminal paper \cite{riesz1918lineare} Riesz  established {\em functional analysis} as a new mathematical discipline.} \cite{riesz1918lineare} to study the inverse of linear operators. 

Since then adjoint has been pervasive in vast literature across mathematics, engineering, and sciences disciplines. This is not surprising as the adjoint has  many appealing features including i) the adjoint operator typically possesses nicer properties  than the original operator (e.g. the adjoint of a densely defined linear operator is a closed operator though the original operator may not), and ii)
the adjoint equation is always linear even when the original equation is not. Though a comprehensive survey on adjoint accounting for its application in many disciplines/fields (and their sub-disciplines) could be desirable to appreciate the crucial role that adjoint plays, it is perhaps an impossible task.\footnote{Or more precisely, it is rather the task for a book than for a paper.} 

Our objective is to provide a window into the adjoint and its crucial role in certain subsets of computational science, engineering, and mathematics. The exposition is necessarily personal and biased based on topics that we are familiar with. The materials developed and presented are  not new, as each or some could be found in (or inferred from) publications in different fields. Our objective is to systematically unify these materials on the same mathematical footing starting from the basic definitions. This expectantly provides a more unified perspective on the usefulness of adjoint in variety of applications. As a result, 
this work could give broader views and better insights into applications of adjoint beyond one field. By establishing general results and then developing materials specific to each application, we bring forth the details on  how abstract concepts/definitions can be translated into particular applications and the connections among them. This paper is written as a tutorial on adjoint with many examples presented with discussions. Though we strike for a self-contained expository, it is necessary for us to state a few results without proof to keep the length of the paper manageable and to focus on the adjoint and its roles.

The paper is structured as follows. In \cref{sect:notations} we introduce various notations, definitions, and some examples. The paper is then developed into two parts: Part I in finite dimensions (\cref{sect:PartI}) and Part II in infinite dimensions (\cref{sect:PartII}). In order to keep the exposition succinct, definitions and results that are valid for both cases are presented/proved once and when that happens we will explicitly state so. Most of our developments start from abstract operator settings and then reduce to standard finite-dimensional settings in $\real^n$ as a special case. In some cases, such as optimization, the order is reversed as we believe it is more natural that way.
Each section of the paper is equipped with examples on which we show how to apply the preceding abstract theoretical results. We make an effort to include practical examples from different fields including linear algebra (e.g. eigendecomposition and the singular value decomposition), ordinary differential equations (an epidemic model), partial differential equations (of elliptic, hyperbolic, and Friedrichs' types), neural networks (feed-forward deep neural networks), least squares and inverse problems (with Tikhonov regularization), PDE-constrained optimization, etc. Due to the interdisciplinary nature of the paper, we do not attempt to provide an exhaustive list of references but a few for each section to keep the references at a manageable length. 


For \cref{sect:PartI}, we begin with the celebrated Riesz representation theorem that is then used to prove the existence of the adjoints of continuous linear operators. It is followed by the closed range theorem that will be useful in many places later. Built upon these basic materials, we shall develop  several applications of adjoint.  The first application  is in \cref{sect:solvabilityApplication} in which we highlight the role of adjoint in assessing the solvability of linear operator equations before solving them. Perhaps one of the most important applications of adjoint is in the study of eigenvalue problems and this is the main focus of \cref{sect:eigenvalueApplication}. The main result for this section is the spectral decomposition of self-adjoint operators in finite dimensions. Another important one is the tight relationship between the orthogonality of a projection and its self-adjointness, which immediately leads to a generalized Pythagorean theorem. In \cref{sect:leastSquaresApplication}, we start with the classical projection theorem and then deploy it together with the closed range theorem to find the necessary and sufficient condition for the optimality of an abstract linear least squares problem. The next important application is the singular value decomposition (SVD) in \cref{sect:svdFinite}, in which we employ the spectral decomposition in \cref{sect:eigenvalueApplication} to establish an SVD decomposition for general linear operator in finite dimensions. This SVD decomposition is then deployed to provide trivial proofs for the closed range theorem, rank-nullity theorem, and the fundamental theorem of linear algebra for abstract linear operators. We then discuss the equivalence of the SVD of an abstract linear operator and the SVD of its matrix representation. 

Our next application of interest for \cref{sect:PartI} is optimization with equality constraints. This is the main topic of \cref{sect:optimization} in which we expose at length the role of adjoint in optimization theory that is valid for both finite and infinite dimensions. This is accomplished by working with Fr\'echet derivative and its Riesz representation counterpart as the gradient. Though can be further developed (e.g. to second-order optimality conditions) our focus for this expository is on the first-order optimality condition. We recall an implicit function theorem and use it to prove an abstract inverse function theorem, which is then deployed to derive the first-order optimality condition for abstract optimization problems with equality constraints. The important role of adjoint comes into the picture when we prove an abstract Lagrangian multiplier theorem using the closed range theorem. The importance of adjoint is further amplified when we study constrained optimization problems with separable structure. Here, adjoint facilitates an efficient gradient-based optimization algorithm in unconstrained reduced subspace while ensuring the feasibility of the contraints at all times. Though discussed in the separate \cref{sect:DNN}, we show that when applying this reduced space approach for the optimization problem arisen from training deep neural networks (DNNs), we recover the backpropagation algorithm that has been the workhorse in training DNNs. Viewing backpropagation from the reduced space approach provides further insights into the algorithm as we shall argue. The last finite dimensional application that we present in \cref{sect:ODEs} is the stability of autonomous ordinary differential equations (ODEs). The main goal of this section is to exhibit the vital role of adjoint in establishing the necessary and sufficient conditions for the asymptotic stability (in the sense of Lyapunov) of ODEs' equilibiria. 

For the infinite dimensional settings in \cref{sect:PartII}, we start with a more general adjoint definition  for densely defined linear operators. This is useful for all subsections except \cref{sect:illposedApplication}. The first application that we study is the illposedness (in Hadamard's sense) nature of inverting compact operators in \cref{sect:illposedApplication}. To that end, we extend the spectral theorem in \cref{sect:eigenvalueApplication} and SVD decomposition theorem in \cref{sect:svdFinite} to the
Hilbert-Schmidt theorem and a general SVD theorem for
compact (linear) operators in infinite dimensions. The main result is a Picard theorem which states the necessary and sufficient conditions under which the task of inverting a compact operator is solvable. This is then employed to show that inverting a compact operator violates at least the stability condition of well-posedness. We then deploy the Riesz-Fredholm theory to show how Tikhonov regularization can restore the well-posedness at the expense of getting a nearby solution. In \cref{sect:PDEapplication}, we discuss the role of adjoint in establishing the well-posedness of abstract linear operator equations with application to partial differential equations (PDEs). Two main results developed in this section are the Banach-Ne$\check{\text{c}}$as-Babu$\check{\text{s}}$ka theorem and the Lax-Milgram lemma. We then switch to \cref{sect:SLapplication} to study Sturm-Liouville eigenvalue problem and generalized Fourier series in $\Ls^2$. Our exposition is for closed linear operators. Using the Hilbert-Schmidt theorem in \cref{sect:illposedApplication} and the Lax-Milgram lemma in \cref{sect:PDEapplication} we establish a spectral decomposition theorems for quite general linear operators, and then apply them to Sturm-Liouville eigenvalue problems to derive Fourier series and its generalization. The last application of interest is PDE-constrained optimization and this is the main topic of \cref{sect:PDEconstrained}. Here we show how to rigorously translate the abstract Lagrangian multiplier theorem in \cref{sect:optimization} to derive the adjoint equation and the reduced gradient for  prototype elliptic and hyperbolic PDEs. We show that the differential operators of adjoint equations are indeed the adjoint operators that we derive at the beginning of  \cref{sect:PartII}.

\section{Notations}
\label{sect:notations}
In this paper, boldface lower case letters such as $\ub$ are reserved for vector-valued functions in $\R^n$, for some integer $n$. 
Calligraphic uppercase letters such as $\Am$ denote matrices, while script uppercase letters such as $\A$ denote operators and superscript $^T$ denotes the transpose of a matrix or a vector.
Bold blackboard upper cases, i.e. $\Xs$ and $\Ys$, are used for spaces and sets. For example $\Hs^n\LRp{\Omega}:= \LRc{\u \in \Ls^2\LRp{\Omega}: \text{ weak derivative up to order } n \text{ residing in } \Ls^2\LRp{\Omega}}$ are standard $\Ls^2$-based Sobolev's spaces \cite{AdamsFournier03, McLean00}.
Lowercase letters are for scalar-valued functions.  \emph{We also use lowercase letters for results that are valid for both finite and infinite dimensional settings} and  
 boldface uppercase letters are for bases of vector spaces. 
Unless otherwise explicitly specified, all spaces are Hilbert spaces endowed with appropriate inner products and the corresponding induced norms. For example, Hilbert space $\Xs$ is endowed with the inner product $\LRp{\u,\v}_\Xs$ for any $\u,\v \in \Xs$, and the induced norm $\nor{\u}_\Xs = \sqrt{\LRp{\u,\u}_\Xs}$. All spaces are either complex or real. For the former, the inner product is conjugate symmetric, i.e., $\LRp{\u,\v}_\Xs = \overline{\LRp{\v,\u}_\Xs}$, where the overline denotes the complex conjugate.  We use $\Fs$ to denote either the set of real ($\R$) or complex ($\C$) numbers, and $\Re$ to denote the operation of taking the real part of a complex number. We shall frequently identify the dual of any Hilbert space with itself. We define $\L\LRp{\Xs,\Ys}$ as the space of all linear operators from $\Xs$ to $\Ys$, $\Bf\LRp{\Xs,\Ys}$ as the space of all bounded linear  operators from $\Xs$ to $\Ys$, and $\Cs\LRp{\Xs,\Ys}$ as the space of all continuous mapping from $\Xs$ into $\Ys$. By $\Cs^n\LRp{\Xs}$ and $\Cs^\infty_0\LRp{\Xs}$ we mean the space of $n$-times continuously differentiable function on $\Xs$ and the space of test functions (infinitely differentiable functions with compact support).
Superscript $^*$ denotes either the topological dual spaces or adjoint operator or the conjugate transpose of a matrix (or a vector). 
Superscript $^\perp$ stands for the orthogonal complement, and by ``$:=$", we mean ``is defined as". 

\begin{definition}[Linear transformation/map/operator]
Consider two inner product vector spaces, $\Xs$, $(\cdot,\cdot)_\Xs$ and $\Ys$, $(\cdot,\cdot)_\Ys$. Suppose $\A:\Xs \rightarrow \Ys$ satisfies the following
\[
\A\LRp{\alpha \u + \beta \w} = \alpha \A \LRp{\u} + \beta \A \LRp{\w}, 
\]
where $\alpha,\beta \in \C$, and  $\u, \w \in \Xs$.
Then, we call\footnote{Though this may not be universal, transformation, map, and operator are used interchangeably in this paper for simplicity.} $\A$  a linear transformation or map, or an  operator from $\Xs$ to $\Ys$.
\end{definition}
\begin{convention}
\end{convention}
\begin{enumerate}
\item For linear operator $\A$, we write: $\A\u := \A\LRp{\u}$.
    \item The domain of $\A$ is defined as 
    \[\Do(\A)  \defeq \LRc{\u \in \Xs: \A(\u) \, \text{is well-defined}} \subset \Xs. \]
    
    \item The range of $\A$ is defined as
    \[\Range(\A) \defeq \LRc{\A(\u) : \u \in \Do\LRp{\A} } = \LRc{\y \in \Ys : \exists \u \in \Do\LRp{\A} \, \text{and} \, \y = \A(\u) } \subset \Ys.
    \]
    
    \item The kernel of $\A$ (or the null space of $\A$) is defined as
    \[ 
    \N(\A) \defeq \LRc{\u \in \Xs : \A(\u) = \zev},
    \]
    where, throughout the paper, $\zev$ denotes either ``zero" function or ``zero" vector in the appropriate space.
\end{enumerate}

\begin{example}
    Let us consider $\Xs = \real^n, \Ys = \real^m$,  and let $\Am:\Xs \mapsto \Ys$ be an $\real^{m \times n}$ matrix. For $\ub, \vb \in \Xs, \, \alpha , \beta \in \Fs$, we have
    \[ 
     \Am \LRp{\alpha \ub + \beta \vb} = \alpha \Am \ub + \beta \Am \vb.  
    \]
    Thus, any matrix is a linear operator.
\end{example}

\begin{example}
    Let $ \A : \Xs =  \Ls^2 \LRp{0,1} \to \Ys = \real  $ be defined such that  for any $f(t) \in \Ls^2\LRp{0,1}$ and a given function $\omega(t) \in \Ls^2 \LRp{0,1}$, we have
    \[
    \A{f} = \int_0^1 \omega(t) f(t) \, dt. 
    \]
Clearly, for $\alpha , \beta \in \Fs$:
    \begin{align*}
        \A \LRp{\alpha f(t) + \beta g(t)} & = \int_0^1 \omega(t) \LRs{\alpha f(t) + \beta g(t)} \, dx \\
        & =  \alpha\int_0^1  \omega(t)f(t) \, dt +  \beta\int_0^1  \omega(t)  g(t) \, dt \\
        & = \alpha \A\LRp{f(t)} + \beta \A \LRp{g(t)},
    \end{align*}
\end{example}
and thus integrals are linear operators.

        \begin{example}
        
    
Consider $ \Xs=\Ys =  \Ls^2 \LRp{0,1}  $ and $\A: \Do\LRp{\A}\subset \Xs \to \Ys$ such that
\[ \A \u = \frac{d^2}{dt^2} \u(t). \]
For $\alpha , \beta \in \Fs$, we have
\begin{align*}
    \A(\alpha \u + \beta \v) & = \frac{d^2}{dt^2} \LRp{\alpha \u(t) + \beta \v(t)} \\
    & = \alpha \frac{d^2}{dt^2} \u(t) + \beta \frac{d^2}{dt^2} \v(t) \\
    & = \alpha 
    \A\u + \beta \A\v,
\end{align*}
thus differentiation is a linear operator.
\end{example}

\section{Part I: Adjoint operators in finite dimensional Hilbert spaces}
\label{sect:PartI}
In this part, unless otherwise stated, we assume that $\Xs$ and $\Ys$ are finite dimensional  vector spaces, i.e.  $\dim\Xs = n < \infty$ and $\dim\Ys = m < \infty$, where $\dim$ denotes the dimension. 
Recall that if $\A\in \L\LRp{\Xs,\Ys}$ and $\dim\Xs < \infty$, then $\A\in \B\LRp{\Xs,\Ys}$. 
Let $\Xb = \LRc{\e_1,\hdots,\e_n}$ and  $\Yb = \LRc{\g_1,\hdots,\g_m}$ be  orthonormal\footnote{Orthonormality is simply for convenience, but not essential.} bases for $\Xs$ and $\Ys$, respectively. 
For any $\u \in \Xs$, we denote by $\ub^\Xb$ the unique vector of coordinates of $\u$ in $\Xb$, and it is easy to see that $\LRp{\ub^\Xb,\vb^\Xb}_{\Fs^n} = \LRp{\u,\v}_{\Xs}$. The matrix representation of $\A$ with respect to the bases $\Xb$ and $\Yb$ is denoted as $\Am^{\Xb\Yb}$. When there is no ambiguity on the bases that we refer to, we simply ignore the superscripts for both coordinate vector and matrix representation. We shall denote the $i$th element of a vector $\ub$ as $\ub(i)$ and the element at the $i$th row and $j$th column of a matrix $\Am$ as $\Am(i,j)$. We also use $\ub_i$ to denote $\ub\LRp{i}$ and this will be clear from the context. We will use square brackets to express matrices and vectors with a finite number of components. Unless otherwise stated, vectors with finite number of components are column vectors.

This section is organized as follows. 
We begin with the celebrated Riesz representation theorem and  the closed range theorem, upon which we shall develop  several applications of adjoint.  The first application  is in \cref{sect:solvabilityApplication} on the solvability of linear operator equations before solving them. The role of adjoint in the study of eigenvalue problems is given in \cref{sect:eigenvalueApplication}.  In \cref{sect:leastSquaresApplication}, we employ the classical projection theorem together with the closed range theorem to find the necessary and sufficient condition for the optimality of an abstract linear least squares problem. The  singular value decomposition (SVD) is the main subject of  \cref{sect:svdFinite}. The SVD decomposition is then deployed to provide trivial proofs for the closed range theorem, rank-nullity theorem, and the fundamental theorem of linear algebra for abstract linear operators. Optimization with equality constraints is the main topic of \cref{sect:optimization} in which we expose at length the role of adjoint in optimization theory that is valid for both finite and infinite dimensions.  This follows by showing that a reduced spaced approach using adjoint reduces to  backpropagation of deep neural networks in \cref{sect:DNN}. The last finite dimensional application that we present in \cref{sect:ODEs} is role of adjoint in establishing the stability of autonomous ordinary differential equations (ODEs).

\begin{theorem}[Riesz representation theorem]
\label{theo:Riesz}
Let $\Lms$ be a bounded linear functional on a Hilbert space $\Xs$.
There exists a unique $u \in \Xs$ such that
\[
\Lms\LRp{v} = \LRp{u,v}_\Xs, \quad \forall v \in \Xs.
\]
Furthermore, the operator norm of $\Lms$ is given as $\nor{\Lms} \equaldef \sup_{v\in \Xs}\frac{\snor{\Lms\LRp{v}}}{\nor{v}_\Xs}  = \nor{u}_\Xs$. 
\end{theorem}
\begin{proof}
    A general proof that works for both finite and infinite dimensional settings is quite standard and can be found in any functional analysis book (see, e.g., \cite{OdenDemkowicz10, ArbogastBona08, brezis2010functional,Luenberger69,Showalter77}). We provide a short and intuitive proof for finite dimensions.
    We prove the result for $\Fs = \Cs$ as the case $\Fs = \R$ is analogous. Let $\Xb = \LRc{\e_1,\hdots,\e_n}$ be an orthonormal basis for $\Xs$.  Let $\ub$ be the representation of $\u$ in $\Xb$.
we have
\[
\u = \sum_{i=1}^n\ub\LRp{i}\e_i \implies \Lms\u =  \sum_{i=1}^n\ub\LRp{i}\underbrace{\Lms\e_i}_{=: \overline{\bs{\ell}\LRp{i}}} = \LRp{\ub,\bs{\ell}}_{\Fs^n} = \LRp{\u,\ell}_\Xs, \quad \forall \u \in \Xs,
\]
where we have defined $\ell$ through its coordinate vector $\bs{\ell}$ in the basis $\Xb$, and thus it is unique.
\end{proof}

\begin{definition}[Adjoint operator]
    Let $\A \in \Bf\LRp{\Xs,\Ys}$. We say that $\A^*:  \Ys \rightarrow \Xs$ is the adjoint of $A$ iff
\[
\LRp{\A\u,\v}_\Ys = \LRp{\u,\A^*\v}_\Xs, \quad \forall \u \in \Xs, \v \in \Ys.
\]
\label{defi:adjoint}
\end{definition}

\begin{proposition}
Let $\A \in \Bf\LRp{\Xs,\Ys}$.
Then $\A^*$ exists and is unique. Furthermore, it is linear with $\nor{\A^*} = \nor{\A}$,
 where the operator norm is defined as usual, e.g.,
    \[
    \nor{\A} := \sup_{\u \in \Xs} \frac{\nor{\A\u}_\Ys}{\nor{\u}_\Xs} = \sup_{\nor{\u}_\Xs = 1} \nor{\A\u}_\Ys.
    \]
 \label{propo:adjointExistence}
\end{proposition}
\begin{proof}
To see the existence and linearity, we note that owing to the continuity of $\A$ and inner product, $\LRp{\A\u,\v}_\Ys$ is continuous in $\u$. By the Riesz representation \cref{theo:Riesz}, there exists a unique $\ell \in \Xs$ depending on $\A$ and $\v$ such that 
\[
\LRp{\u,\ell\LRp{\A,v}}_\Xs = \LRp{\A\u,\v}_\Ys, \quad \forall \u \in \Xs,
\]
which implies that $\ell\LRp{\A,v}$ is linear in $\v$. Defining $\A^*\v := \ell\LRp{\A,v}$ shows the $\A^*$ exists and linear. We next show that $\A^*$ is continuous (bounded) and $\nor{\A^*} = \nor{\A}$. We have
\[
\nor{\A^*\v}^2_\Xs = \LRp{\A^*\v,\A^*\v}_\Xs = \LRp{\A\A^*\v,\v}_\Ys \le \nor{\A}\nor{\A^*\v}_\Xs\nor{\v}_\Ys,
\]
which shows that $\A^*$ is bounded and $\nor{\A^*} \le \nor{\A}$. Since $\A$ is the adjoint of $\A^*$, following a similar arguement we have $\nor{\A} \le \nor{\A^*}$, and this concludes the proof.
\end{proof}

\begin{remark}
    Note that \cref{theo:Riesz}, \cref{defi:adjoint}, and \cref{propo:adjointExistence} are also valid for infinite-dimensional settings. 
\end{remark}

\begin{example}
 Let $\Us = \R^n$ and $\Vs = \R^n$ be respectively endowed with the inner products $(.,.)_{\real^n}$ and an $\Mm$-weighted inner product $(.,.)_{\real^n,\Mm}$ where $ \LRp{\vb,\wb}_{\real^n,\Mm} \defeq \sum_{i,j} \vb(i) \Mm\LRp{i,j} \wb(j) \defeq \vb^T \Mm \wb,\forall \vb,\wb \in \Vs$, and $\Mm$ is a symmetric and positive definite matrix.
        We need to find the adjoint operator $\Am^*$ of a matrix $\Am: \Us, (.,.)_{\real^n} \rightarrow \Vs, (.,.)_{\real^n,\Mm}$. We have
        \begin{align*}
            (\Am\ub,\mb{v})_\mbb{V} & = \LRp{\Am\ub,\vb}_{\real^n,\Mm} = \LRp{\Am\vb}^T \Mm\bs{v}\\ 
            &= \bs{u}^T \Am^T \Mm\bs{v} = \bs{u}^T (\Am^T \Mm \bs{v}) =(\bs{u}, \Am^T \Mm \bs{v})_{\real^n} \\
            & = (\bs{u}, \Am^T \Mm \bs{v})_{\mbb{U}}.
        \end{align*}
        By the definition of adjoint operator, we have $\Am^* = \Am^T \Mm$.
        \label{exa:matrixAdjoint}
\end{example}

\begin{example}
      Now, let us consider $ \A : \Us =  \text{Span} \{1,\ x,\ x^2\} \subset \Xs = \Ls^2\LRp{-1,1} \to \real^2   $ such that  the map A is defined as
    \[
    \u(x) \in \Us \mapsto \A \u=\begin{bmatrix}
\int_{-1}^1 \u(x)\ dx \\ 
\int_{-1}^1 (2x+1) u(x)\ dx
\end{bmatrix} 
,
\]
and the inner product on $\Ls^2\LRp{-1,1}$, and hence on $\Us$, is defined as
\[
\LRp{\u,\v}_{\Ls^2\LRp{-1,1}} := \int_{-1}^1\u(x){\v(x)}\,dx.
\]

We have
\begin{multline*}
\LRp{\A\u,\ \bs{\v}}_{\real^2}=\vb_1\int_{-1}^1 \u(x)\ dx+\vb_2 \int_{-1}^1 (2x+1) \ \u(x)\ dx    \\ 
=\int_{-1}^1 \u(x) \LRs{\vb_1+(2x+1)\vb_2}dx=  \int_{-1}^1 u(x) \LRs{1, \ (2x+1)}\bs{v}\ dx.
\end{multline*}
Therefore, by definition $\A^*=\LRs{1, \ (2x+1)}$.
\exalab{quadraticFit}
\end{example}

\begin{example}
      Let us consider $ \A : \Us =  \text{Span} \{1,\ x,\ x^2\} \subset \Xs = \Ls^2\LRp{-1,1} \to \real^2   $ such that  the map A is defined as
    \[
    \u(x) \in \Us \mapsto \A \u=\begin{bmatrix}
\u\LRp{x_1} \\ 
\u\LRp{x_2}
\end{bmatrix},
\]
where $x_1, x_2 \in \LRp{-1,1}$ and $x_1 \ne \x_2$. 
The inner product on $\Ls^2\LRp{-1,1}$, and hence on $\Us$, is defined as
\[
\LRp{\u,\v}_{\Ls^2\LRp{-1,1}} := \int_{-1}^1\u(x){\v(x)}\,dx.
\]

We have
\[
\LRp{\A\u,\ \bs{\v}}_{\real^2}=\bs{\v}_1 \u(x_1) +  \vb_2 \u(x_2) = \int_{-1}^1 \u\LRp{x}\LRs{\delta\LRp{x-x_1}, \delta\LRp{x-x_2}}\vb dx.
\]
Therefore, by definition $\A^*=\LRs{\delta\LRp{x-x_1}, \delta\LRp{x-x_2}}$. Here, we have defined $\delta\LRp{x-y} \in \Us$ via
\begin{equation}
\int_{-1}^1\delta\LRp{x-y}\u(\x)\,dx = \u(y), \quad \forall \u \in \Us,
\eqnlab{DiracDelta2}
\end{equation}
and thus
\[
\delta\LRp{x-y} = 15\frac{\LRp{3y^2-1}}{8}x^2 + 3\frac{y}{2}x + 3\frac{\LRp{3-5y^2}}{8},
\]
by testing \cref{eq:DiracDelta2} with $\u \in\LRc{1,x,x^2}$.
\end{example}

\begin{example}
    Consider $\Ps^n\LRs{0,1}$ the set of complex-valued polynomial of order at most $n$ on $\LRs{0,1}$. We define $\A:\Us := \Ps^n\LRs{0,1} \subset \Ls^2\LRp{0,1} \to \Us$ as 
\[
\A\u := x \u' := x\dd{\u}{x},
\]
    and the inner product on $\Ls^2\LRp{0,1}$, and hence on $\Us$, is defined as
\[
\LRp{\u,\v}_{\Ls^2\LRp{0,1}} := \int_{0}^1\u(x)\overline{\v(x)}\,dx.
\]
By integration by parts we have
\[
\LRp{\A\u,\v}_{\Ls^2\LRp{0,1}} = \int_{0}^1x\u'\overline{\v}\,dx = u(1)\overline{\v(1)} + \LRp{u,-\LRp{x\v}'}_{\Ls^2\LRp{0,1}} = 
\LRp{u,\delta\LRp{x-1}\v(1)-\LRp{x\v}'}_{\Ls^2\LRp{0,1}},
\]
which by definition gives
\[
\A^*\v = \delta\LRp{x-1}\v(1)-\LRp{x\v}',
\]
where we have defined $\delta\LRp{x-1} \in \Us$ via
\begin{equation}
\int_{0}^1\u(\x)\overline{\delta\LRp{x-1}}\,dx = \u(1), \quad \forall \u \in \Us.
\eqnlab{DiracDeltan}
\end{equation}
Clearly, we can find $\delta\LRp{x-1}$ as the unique polynomial of degree at most $n$
 by testing \cref{eq:DiracDeltan} with $\u \in \LRc{1,x,\hdots,x^n}$.
\end{example}


\begin{proposition}
    Let $\Xb$ and $\Yb$ be orthonormal bases of $\Xs$ and $\Ys$, respectively, and $\dim\Xs = n$ and $\dim\Ys = m$. Let $\Am$ and $\Bm$ be the matrix representations of $\A$ and $\A^*$ with respect to the bases $\Xb$ and $\Yb$. Then
    \[
    \Bm = \Am^*,
    \]
    where $\Am^*$ be the conjugate transpose of $\Am$.
    \label{propo:adjointT}
\end{proposition}
\begin{proof}
    By the definition of adjoint and matrix representation,  for any $\u \in \Xs$ and $\v \in \Ys$ we have
    \[
\LRp{\ub,\Bm\vb}_{\Fs^n} = \LRp{\u,\A^*v}_\Xs  = \LRp{\A\u,v}_\Ys = \LRp{\Am\ub,\vb}_{\Fs^m}
    = \LRp{\ub,\Am^*\vb}_{\Fs^n},
    \]
    which concludes the proof.
\end{proof}


The following \cref{defi:orthoComplement}, \cref{defi:closure}, \cref{coro:perpClosure}, \cref{theo:CRT}, and \cref{coro:directSum} are valid for both finite and infinite dimensions.
\begin{definition}[Orthogonal complement]
    Let $\Ss \subset \Xs$, the orthogonal complement $\Ss^\perp$  of $\Ss$ is defined as
    \[
    \Ss^\perp := \LRc{\u \in \Xs: \LRp{\u,\w}_\Xs = 0, \forall \w \in \Ss }.
    \]
    \label{defi:orthoComplement}
\end{definition}

A direct consequence of the definition is that $\Ss^\perp$ is a closed subspace of $\Xs$ and that $\Ss \cap \Ss^\perp = \LRc{\theta}$.

\begin{definition}[Closure]
    Let $\Ss \in \Xs$. The closure $\overline{\Ss}$ of $\Ss$ is the smallest closed set containing $\Ss$.
    \label{defi:closure}
\end{definition}
Note that we use the overline ``$\overline{\cdot}$" to denote both complex conjugate and the closure, but it should be clear from the context.
\begin{corollary}
    There holds: $\LRp{\Ss^\perp}^\perp = \overline{\Ss}$.
    \label{coro:perpClosure}
\end{corollary}
\begin{proof}
    See \cite[proposition 1 of chapter 3]{Luenberger69}.
\end{proof}

Next is an important theorem (see, e.g., \cite{Yosida1995,brezis2010functional, Banach1978-mx}).
\begin{theorem}[The closed range theorem]
Let $\A: \Xs \to \Ys$. 
The following hold:
\begin{itemize}
\item $\LRs{\Range\LRp{\A}}^\perp = \N\LRp{\A^*}$.
\item $\overline{\Range\LRp{\A}} = \LRs{\N\LRp{\A^*}}^\perp$.
\item $\LRs{\Range\LRp{\A^*}}^\perp = \N\LRp{\A}$.
\item $\overline{\Range\LRp{\A^*}} = \LRs{\N\LRp{\A}}^\perp$.\footnote{Since we consider only Hilbert spaces, which are reflexive, $\overline{\Range\LRp{\A^*}} = \LRs{\N\LRp{\A}}^\perp$ holds. In general, $\overline{\Range\LRp{\A^*}} \subset \LRs{\N\LRp{\A}}^\perp$: see \cite[Corollary 2.18]{brezis2010functional}.}
\end{itemize}
If $\Range\LRp{\A}$ is closed, so is $\Range\LRp{\A^*}$, and we can replace
$\overline{\Range\LRp{\A}}$ and $\overline{\Range\LRp{\A}^*}$ by $\Range\LRp{\A}$ and $\Range\LRp{\A^*}$, respectively, in the above results.
\label{theo:CRT}
\end{theorem}
\begin{proof}
The second assertion is the direct consequence of the first assertion and \cref{coro:perpClosure}. The third and fourth assertions follow the first and the second, and the fact that $\LRp{\A^*}^* = \A$. So, we only need to prove the first assertion. Let $\z \in \N\LRp{\A^*}$ and $y \in \Range\LRp{A}$. Then $y = \A x$ for some $x \in \Xs$. We have
\[
\LRp{\z,y}_\Ys = \LRp{\z,\A x}_\Ys = \LRp{\A^*\z, x}_\Xs = 0, \forall y \in \Range\LRp{A},
\]
which says that $\N\LRp{\A^*} \subset \LRs{\Range\LRp{\A}}^\perp$. Now take $\z \in \LRs{\Range\LRp{\A}}^\perp$, we have
\[
\LRp{\A^*\z, x}_\Xs = \LRp{\z,\A x}_\Ys  = 0, \quad \forall x \in \Xs,
\]
which implies that $\A^*\z = 0$, which in turn shows $\LRs{\Range\LRp{\A}}^\perp \subset \N\LRp{\A^*}$.
\end{proof}
\begin{corollary}
There hold:
\begin{itemize}
    \item $\Xs = \N\LRp{\A}\oplus\overline{\Range\LRp{\A^*}}$
        \item $\Ys = \N\LRp{\A^*}\oplus\overline{\Range\LRp{\A}}$   
\end{itemize}
\label{coro:directSum}
\end{corollary}
We note that for finite dimensional vector spaces $\Xs$ and $\Ys$, $\Range\LRp{\A}$, and hence $\Range\LRp{\A^*}$, is obviously closed, and thus \cref{theo:CRT} and \cref{coro:directSum} clearly hold.
We will see that the proof of the closed range \cref{theo:CRT}, hence the \cref{coro:directSum}, for finite dimensional cases is trivial using the SVD decomposition in \cref{theo:SVDfinite}.

\subsection{Application of adjoint to the solvability of linear operator equations}
\label{sect:solvabilityApplication}
In this section, we study the existence  of a solution of the following linear operator equation
\begin{equation}
\A\u = \f,   
\eqnlab{linearEqn}
\end{equation}
where $\A:\Xs \to \Ys$.
\begin{lemma}
$ $
\begin{itemize}
    \item {\bf Existence.} The linear equation \cref{eq:linearEqn} has a solution iff $\y \in \N\LRp{\A^*}^\perp$.
    \item {\bf Uniqueness.} The solution of  \cref{eq:linearEqn} is unique iff $\N\LRp{\A} = \LRc{\theta}$.
    \item If $\dim\Xs = \dim\Ys$, the uniqueness is equivalent to existence.
\end{itemize}
\end{lemma}
\begin{proof}
The existence is the direct consequence of \cref{theo:CRT}. The uniqueness is clear. The proof of the third assertion is as follows. We have
    \begin{align*}
        \N(\A) & = \LRc{\theta} & \text{ Uniqueness}\\
         &\Updownarrow &        \\
         \dim \N(\A) &= 0 \\
        & \Updownarrow  &  \dim \mbb{X} = \dim \mbb{Y} \text{ and } \dim \mbb{X} = \dim \N(\A) + \dim \Range(\A) \\
        \dim \mbb{Y} = \dim \mbb{X} & = \dim \Range(\A) &\\
        &\Updownarrow & \Range(\A) \subseteq \V\\
        \Range(\A) &= \mbb{Y} & \text{ Existence for any } \y \in \Ys
    \end{align*}  
\end{proof}

\begin{remark}
The existence condition can be simply $\y \in \R\LRp{\A}$. However, it is easier to work with $\N\LRp{\A^*}^\perp$ as it gives us equations (see \cref{defi:orthoComplement}) to determine/characterize  $\N\LRp{\A^*}^\perp$.
\end{remark}

\begin{example}
     Consider the operator $\A$ defined in \cref{exa:quadraticFit}, and we are interested in studying the existence of a solution\footnote{Note that this is an operator setting for the problem of fitting a quadratic polynomial $\u(x)$ with two pieces of information about $\u(x)$.} for $\A\u=\fb$. We start by recalling from \cref{exa:quadraticFit} that
\[
\A^*=\LRs{1, \ (2x+1)}.
\]
 To compute $\N(\A^*)$, we pick any $v \in \N(\A^*)$, i.e. $\A^*\vec{v} = 0$. We have
 
 \[v_1+(2x+1)v_2=0 \quad \forall x,
 \]
which implies 
 \[v_1=v_2=0. 
 \]
 Thus,
 \[\N(\A^*)=\{\theta\}. 
 \]
 
  Since $\theta$ is orthogonal to any $\fb$, we conclude that $\A\u=\fb$ always has a solution.

    
    





 
 
\end{example}

\begin{example}
    Let $\Am = \begin{bmatrix}     1 & 2 \\ 1 & 2     \end{bmatrix}$, and $ \fb = \begin{bmatrix} 1 \\ 0 \end{bmatrix}$. The question is if there is a solution to the equation $\Am\ub = \fb$. Consider ${\mbb{U}} = {\mbb{V}} = \real^2
    $ with the standard Euclidean inner product $(\cdot,\cdot)_{\R^2}$. We know from \cref{exa:matrixAdjoint} that the adjoint $\Am^*$ is given by
    \[ \Am^* = \Am^T = \begin{bmatrix} 1 & 1 \\ 2 & 2 \end{bmatrix}.\]
     Let's determine the null space of $\Am^*$. If $\mb{z} \in \N(\Am^*)$ then 
    \[\Am^* \mb{z} = \theta \implies \zb_1 + \zb_2 = 0,\]
    which yields $\N(\Am^*) = \LRc{\begin{bmatrix} \alpha \\ -\alpha \end{bmatrix} : \quad \forall \alpha \in \real}$. However, for any $\zb \in \N\LRp{\Am^*}$ we have
    \[\LRp{\fb,\zb}_{\R^2} = 
    \LRp{\begin{bmatrix} 1 \\ 0 \end{bmatrix},\begin{bmatrix} \alpha\\-\alpha \end{bmatrix}}_{\R^2} = \alpha \quad \forall \alpha \in \real, \]
    that is, $\fb \not\perp \N(\Am^*)$. Thus, $\Am\ub = \fb$ does not have a solution.
    \label{exa:quadraticFitNonUniqueness}
\end{example}


\subsection{Application of adjoint to eigenvalue problems}
\label{sect:eigenvalueApplication}

\begin{definition}[Eigenvalue problem]
Let $A: \Xs \rightarrow \Xs$ be a linear operator.
\begin{equation}
\A \u = \lambda \u \quad \quad \forall\lambda \in \C,\ \u \in \Xs
\label{eq:eigProblem}
\end{equation}
is called an eigenvalue problem 
if there exists a nontrivial  pair $\LRp{\lambda,x}$ ({\em x is not a zero vector/function but $\lambda$ could be zero}). In particular:
\begin{itemize}
    \item $\lambda$ is called an eigenvalue.
    \item $x$ is called an eigenfunction, associated with the eigenvalue $\lambda$, of $\A$. If $\Xs$ is a finite-dimensional space, i.e., $\Xs = \real^n$, $x$ is typically called eigenvector.
\end{itemize}
\label{defi:eigen}
\end{definition}
\begin{definition}[Self-adjoint operator]
    If $\A^* = \A$, then $\A$ is called self-adjoint.
    \label{defi:selfadjoint}
\end{definition}
\begin{lemma}
Let $\A: \Xs \rightarrow \Xs$ be a linear operator and $\A$ is self-adjoint. Then:
\begin{enumerate}
    \item Eigenvalues of $\A$ are real.
    \item Eigenfunctions corresponding to distinct eigenvalues are orthogonal to each other. That is, if $\LRp{\lambda,\u}$ and $\LRp{\alpha,\v}$ are two eigen-pairs and $\lambda \neq \alpha$ then $\LRp{\u,\v}_\Xs= 0$.
\end{enumerate}
\label{lem:realEigen}
\end{lemma}
\begin{proof}
For the first assertion, we have $\lambda\LRp{\u,\v}_\Xs = \LRp{\lambda\u,\u}_\Xs = \LRp{\A\u,\u}_\Xs = \LRp{\u,\A\u}_\Xs = \LRp{\u,\lambda\u}_\Xs = \overline{\lambda}\LRp{\u,\lambda\u}_\Xs$. Thus, $\LRp{\lambda - \overline{\lambda}}\LRp{\u,\lambda\u}_\Xs = 0$, and this implies $\lambda = \overline{\lambda}$, or $\lambda$ is real. For the second assertion, we observe that $\lambda\LRp{\u,\v}_\Xs = \LRp{\A\u,\v}_\Xs = \LRp{\u,\A\v}_\Xs = \alpha\LRp{\u,\v}_\Xs$. Therefore, $\LRp{\lambda - \alpha}\LRp{\u,\v}_\Xs = 0$, and this implies $\LRp{\u,\v}_\Xs = 0$.
\end{proof}

\begin{remark}
    { Note that \cref{defi:eigen}, \cref{defi:selfadjoint}, and \cref{lem:realEigen} hold for both finite and infinite dimensional cases.}
\end{remark}

\begin{proposition}
    Let $\A:\Xs \to \Xs$ be a linear operator. Then, $\A$ has at least one eigenvalue.
    \propolab{eigExistence}
\end{proposition}
\begin{proof}
    Let $n$ be the dimension of $\Xs$ and $\Xb$ be a basis. Let $\Am$ and $\ub$ be the matrix and vector presentation of $\A$ and $\u$ in the basis $\Xb$. The matrix representation of the eigenvalue problem\cref{eq:eigProblem} reads
    \[
    \Am\ub = \lambda \ub,
    \]
    that is, $\lambda$ is an eigenvalue of $\A$ iff it is also an eigenvalue of $\Am$. Since $\det\LRp{\Am - \lambda\Im} = 0$ has $n$ roots for $\lambda$ (including repeated ones), there is at least one eigenvalue.
\end{proof}

\begin{theorem}
    Let $\A:\Xs \to \Xs$ be a self-adjoint linear  operator. Then, an orthonormal basis of $\Xs$ can be constructed from eigenfunctions of $\A$. 
    \theolab{eigComplete}
\end{theorem}
\begin{proof}
    \cref{propo:eigExistence} implies that $\A$ has at least one eigenfunctions. Let $\Ss$ be the span of all eigenfunctions of $\A$. Owing to \cref{lem:realEigen}, it is sufficient to show that $\Ss = \Xs$. If $\Ss^\perp = \LRc{\theta}$, then clearly $\Ss = \Xs$. Now suppose $\Ss^\perp \ne \LRc{\theta}$. If $\u \in \Ss^\perp$ and $\LRp{\lambda,\v}$ be an eigen-pair of $\A$, then $\LRp{\A\u,\v}_\Xs = \LRp{\u,\A\v}_\Xs = \lambda\LRp{\u,\v}_\Xs = 0$ since $\v \in \Ss$. Thus $\A: \Ss^\perp\to \Ss^\perp$. By \cref{propo:eigExistence}, $\A$ has an eigenfunction $\w$ in $\Ss^\perp$, which means that $\w \in \Ss\cap\Ss^\perp$, which in turn implies $\w = \theta$, a contradiction. We conclude that $\Ss^\perp = \LRc{\theta}$ and this concludes the proof.
\end{proof}
\begin{corollary}[Spectral decomposition of self-adjoint operators in finite dimensions]
Let $\dim\Xs = n$ and $\A: \Xs \rightarrow \Xs$ be a linear and self-adjoint operator. There exists $n$ real values, $\lambda_1\geq \lambda_2\geq\lambda_3 \dots \geq\lambda_n$ and orthonormal vectors $\u_1, \u_2,\ \dots \u_n$ such that:
\begin{enumerate}
    \item $ \A\u_i=\lambda_i\u_i$. 
    \item For any $\x \in \Xs$ we have 
    \[ 
    \A\x=\sum_{i=1}^n \lambda_i \LRp{\x,\u_i}_{\Xs}  \u_i, \implies \A= \sum_{i=1}^n \lambda_i \LRp{\cdot,\u_i}_\X\ u_i,
    \]
    that is, $\A$ is completely determined by its eigenpairs.
\end{enumerate}   
\label{coro:spectralFinite}
\end{corollary}
\begin{proof}
The first assertion is clear due to
    \cref{lem:realEigen}, \cref{propo:eigExistence}, and \cref{theo:eigComplete}. The second assertion is also obvious since
    \begin{equation}
    \x = \sum_{i = 1}^n\LRp{\x,\u}_{\Xs}\u_i,
    \label{eq:eigComplete}
    \end{equation}
    due to \cref{theo:eigComplete}.
\end{proof}

\begin{example}[Eigen-decomposition of self-adjoint (Hermitian) matrices]
    Let $\Am: \Fs^n \to \Fs^n$ be a self-adjoint  matrix. Applying \cref{coro:spectralFinite} we conclude that $\Am$ has $n$ real eigenvalues $\lambda_1\geq \lambda_2\geq\lambda_3 \dots \geq\lambda_n$ and orthonormal eigenvectors $\ub_1, \ub_2,\ \dots \ub_n$ such that $\Am\ub_i = \lambda_i\ub_i$ and
    \[
    \Am= \sum_{i=1}^n \lambda_i \LRp{\cdot,\ub_i}_{\Fs^n} \ub_i = \sum_{i=1}^n \lambda_i \ub_i \ub_i^* = \U \Lambda \U^*,
    \]
    where $\U \in \Fs^{n\times n}$ is the eigenmatrix whose columns are eigenvectors of $\Am$ and $\Lambda$ is a diagonal matrix whose diagonals are the corresponding eigenvalues of $\Am$. Thus, the standard eigendecomposition for self-adjoint matrices is a special case of \cref{coro:spectralFinite}. 
\end{example}

We next discuss the relationship between self-adjointness and orthogonal projection. 
\begin{definition}[Orthogonal projection]
    A linear operator $\P: \Xs \to \Xs$ is a projection if $\P^2 := \P\P = \P$. If, in addition, $\Range\LRp{\P} \perp \N\LRp{\P}$, then $\P$ is an orthogonal projection.
\end{definition}
\begin{proposition}
    A projection $\P: \Xs \to \Xs$ is orthogonal iff $\P$ is self-adjoint.
\end{proposition}
\begin{proof}
For any $\x \in \N\LRp{\P}$ we have
\[
    \LRp{\P\y, \x}_{\Xs} = \LRp{\y, \P^*\x}_{\Xs} = \LRp{\y, \P\x}_{\Xs} = 0, \quad \forall \y \in \Xs,
\]  
which ends the proof. Another way to see this is to use the result $\N\LRp{\P} \perp \Range\LRp{\P^*}$ from the closed range \cref{theo:CRT}, but we omit the details for brevity.
\end{proof}
     
Inspired by the spectral decomposition of a self-adjoint operator, we define
\[
\P := \sum_{i = 1}^n\LRp{\cdot,\u_i}_{\Xs}\u_i,
\]
where $\LRc{\u_1,\cdots,\u_n}$ is an orthonormal basis of $\Xs$. It is easy to verify that: i) $\P^2 := \P\P = \P$, and ii) $\P$ is self-adjoint. Thus $\P$ is an orthogonal projection into $\Xs$.  More generally, we can verify that
\[
\P^r := \sum_{i = 1}^r\LRp{\cdot,\u_i}_{\Xs}\u_i,
\]
where $r \le n$, is an orthogonal projection into $\Range\LRp{\P^r}$ spanned by $\LRc{\u_1,\cdots,\u_r}$.
Any orthogonal projection $\P$ orthogonally projects $\Xs$ into $\Range\LRp{\P}$
while 
\[
\I -\P,
\]
where $\I$ is the identity operator,
is the orthogonal projection into $\N\LRp{\P}$. Indeed, by the self-adjointness, we have
\[
\LRp{\P\x,y - \P\y}_{\Xs}
= \LRp{\x, \P\y - \P^2\y}_{\Xs} = 0.
\]
We thus have the following generalized Pythagorean theorem for any $\x \in \Xs$ and an orthogonal projection $\P$:
\begin{equation}
\nor{\x}_\Xs^2 = \nor{\LRp{\I - \P}\x}_\Xs^2 + \nor{\P\x}_\Xs^2.
\label{eq:Pythagorean}
\end{equation}

\subsection{Application of adjoint to linear least squares problems}
\label{sect:leastSquaresApplication}
We start with the classical projection theorem that holds for both finite and infinite dimensional settings.
\begin{theorem}[Projection theorem]
    Let $\Ss$ be a  subspace of a pre-Hilbert\footnote{A pre-Hilbert space is an incomplete metric space with an inner product.} space $\Ys$. Let $\y \in \Ys$. Then, $\u \in \Ss$ is the unique minimizer of $\inf_{\w \in \Ss}\nor{\w-\y}_\Ys$ iff $\LRp{\y- \u} \perp \Ss$. 
    The existence of the minimizer $\u$ is guaranteed if $\Ys$ is Hilbert and $\Ss$ is closed. 
    \label{theo:projection}
\end{theorem}
\begin{proof}
    We follow closely the proof by contradiction in \cite[Theorem 2]{Luenberger69}. Suppose there exists $\v \in \Ss$ that is not orthogonal to $\LRp{\y-\u}$. We can assume that $\LRp{\y-\u,\v}_\Xs = \epsilon \ne 0$ and $\nor{\v}_\Xs = 1$. We have
    \[
    \nor{\y - \u-\epsilon\v}_\Ys^2 = \nor{\y - \u}_\Ys^2 +\snor{\epsilon}^2 - 2\epsilon\overline{\epsilon} =
    \nor{\y - \u}_\Ys^2 - \snor{\epsilon}^2 < \nor{\y - \u}_\Ys^2,
    \]
    contradicting the fact that $\u$ is a minimizer. Conversely, let $\LRp{\y-\u} \perp \Ss$, by the Pythagorean identity \cref{eq:Pythagorean}, for any $\v \in \Ss$ we have
    \[
    \nor{\y-\v}_\Ys^2 = \nor{\y-\u + \u -\v}_\Ys^2 = \nor{\y-\u}_\Ys^2 + \nor{\u-\v}_\Ys^2 \ge \nor{\y-\u}_\Ys^2,
    \]
    which shows that $\u$ is the unique minimizer. The existence proof is lengthy and hence is omitted. 
\end{proof}
The following corollary is also valid for both finite and infinite dimensions.
\begin{corollary}[Linear least squares]
    Let $\Xs, \Ys$ be pre-Hilbert and $\A:\Xs \to \Ys$ be linear. Then, for any $\y \in \Ys$, $\tilde{\x} \in \Xs$ is a minimizer of $\inf_{\x \in \Xs}\nor{\A\x-\y}_\Ys$ iff
    \begin{equation}
        \A^*\A\tilde{\x} = \A^*\y.
        \label{eq:LSsolution}
    \end{equation}
    Furthermore, if $\A$ is injective then the minimizer $\tilde{\x}$ is unique.
    \label{coro:linearLS}
\end{corollary}
\begin{proof}
    Note that $\inf_{\x \in \Xs}\nor{\A\x-\y}_\Ys$ is equivalent to $\inf_{\w \in \Range\LRp{\A}}\nor{\w-\y}_\Ys$. Applying \cref{theo:projection} we know that  $\u \in \Range\LRp{\A}$ is a minimizer of $\inf_{\w \in \Range\LRp{\A}}\nor{\w-\y}_\Ys$ iff $\LRp{\y-\u}\perp \Range\LRp{\A}$, that is, by the Closed Range \cref{theo:CRT},
    \[
    \LRp{\y-\u} \in \N\LRp{\A^*} \Leftrightarrow \A^*\LRp{\y - \u} = \theta.
    \]
    Since $\u \in \Range\LRp{\A}$, there exists $\tilde{\x} \in \Xs$ such that 
    \[
    \A^*\LRp{\y - \A\tilde{\x}} = \theta,
    \]
    which concludes the first assertion.

    For the second assertion, 
    the uniqueness in \cref{theo:projection} leads to the uniqueness of $\tilde{\xb}$ when $\A$ is injective. Another way to see this is to note that $\N\LRp{\A^*\A} = \N\LRp{\A}$ and thus the injectivity of $\A$ is equivalent to the injectivity of $\A^*\A$. The least square solution in \cref{eq:LSsolution} is therefore unique.
\end{proof}

The beauty here is that \cref{eq:LSsolution} is exactly the first order optimality condition that is typically obtained by requiring the derivative of $\nor{\A\x-\y}_\Ys$, with respective to $\x$, to vanish (see the same result via derivative for linear least squares problem in \cref{exa:LSsolutionDerivative}). When $\dim\LRp{\Xs} < \infty$, then $\Range\LRp{\A}$ is finite dimension and hence closed in $\Ys$. If, additionally, $\Ys$ is Hilbert then the existence of $\tilde{\x}$ is guaranteed by \cref{theo:projection}.

\begin{example}[linear least squares in finite dimensions]
Consider the operator $\A$ defined in \cref{exa:quadraticFitNonUniqueness} and we are interested in minimizing $\nor{\A\x-\yb}_{\R^2}$ for some given $\yb \in \R^2$. By \cref{coro:linearLS}, a minimizer $\tilde{\x}$ must satisfy
\[
\A^*\A\tilde{\x} = \A^*\yb.
\]
Since $\A$ is not injective, there are multiple minimizers for this problem. This is consistent with the non-uniqueness in \cref{exa:quadraticFitNonUniqueness}. Note that we can reformulate the operator form $\nor{\A\x-\yb}_{\R^2}$ by an equivalent matrix representation form. Indeed, let $\Am$ be the matrix representation of $\A$ with respect to an orthonormal basis in $\Us =  \text{Span} \{1,\ x,\ x^2\}$ and the canonical basis of $\R^2$. Let $\xb$ be the coordinate vector of $\x$ in the same orthonormal basis of $\Us =  \text{Span} \{1,\ x,\ x^2\}$, we have $\nor{\A\x-\yb}_{\R^2}$ = $\nor{\Am\xb-\yb}_{\R^2}$. Again, $\Am$ is not injective and a solution is not unique.
\end{example}

\begin{example}[linear least squares with matrices]
Consider $\Am:\Fs^n \to \Fs^m$ and $\yb \in \Fs^m$. Applying \cref{coro:linearLS} we have that there exists $\tilde{\xb} \in \Fs^n$ minimizing $\nor{\Am\xb-\yb}_{\Fs^m}$ and
\begin{equation*}
\Am^*\Am\tilde{\xb} = \Am^*\yb.    
\end{equation*}
When $\Am$ is injective or equivalently full column rank, we have the uniqueness of $\tilde{\xb}$ again by \cref{coro:linearLS}. 
Another way to see this is that in this case $\Am^*\Am$ is an invertible matrix which implies the uniqueness of the minimizer $\tilde{\xb}$.
\end{example}

\subsection{Application to the singular value decomposition (SVD)}
\label{sect:svdFinite}
Consider $\A: \Xs \rightarrow \Ys$ with $\dim(\Xs)=n$ and $\dim(\Ys)=m$.  It is clear that the operators $\A^*\A: \Xs \rightarrow \Xs$ and $\A\A^*: \Ys \rightarrow \Ys$ are linear and self-adjoint.  The spectral decomposition 
 of self-adjoint operator in \cref{coro:spectralFinite}
 states that $\exists \lambda_1\geq \lambda_2\geq \hdots \ge \lambda_n$ and an orthonormal basis $\u_1,\ \u_2, \hdots, \u_n$ of $\Xs$ such that 
 \[
 \A^*\A\u_i=\lambda_i\u_i,
 \]
 which implies, after taking the inner product both sides with $\u_i$,
 \[
 \nor{\A\u_i}^2_\Ys = \lambda_i\nor{\u_i}^2_\Ys,
 \]
 and thus $\lambda_i \ge 0$.
We  define $\sigma^2_i = \lambda_i$ and have
\begin{equation}
    \A^*\A\u_i=\sigma^2_i\u_i. 
    \label{eq:eig_AAs}
\end{equation}

\begin{theorem}[Singular Value Decomposition (SVD)]
Let $\A: \Xs \rightarrow \Ys$ with $\dim(\Xs)=n$ and $\dim(\Ys)=m$. Then, there exist $\LRc{\sigma_i,\u_i,\v_i}$ (the singular triplets of $\A$) with $\sigma_1\geq \hdots \ge \sigma_i\ge \dots \ge\sigma_k \geq 0$ and $k=min \{ n,\ m\}$, an orthonormal basis $\{\u_1, \u_2,\hdots, \u_n\}$ of $\Xs$, and an orthonormal basis $\{\v_1, \v_2,\hdots, \v_m\}$ of $\Ys$ such that: 
\begin{enumerate}
    \item $ \A\u_i={\sigma_i}\v_i$ for $i=1,\dots, r$ and $ \A\u_i = \theta$ for $i = r+1,\hdots,n$,
    \item $ \A^*\v_j={\sigma_j}\u_j $ for $ j=1,\hdots, r$ and and $ \A^*\v_j = \theta$ for $j = r+1,\hdots,m$,
    \item $\A$ is completely determined by its singular triplets in the following sense: for any $\x \in \Xs$, we have 
    \[ 
    \A\x=\sum_{i=1}^r \sigma_i \LRp{\x,\u_i}_{\Xs}  \v_i, \implies \A= \sum_{i=1}^r \sigma_i \LRp{\cdot,\u_i}_\Xs \v_i,
    \]
\end{enumerate}
where $r$ is the maximum index for which $\sigma_r > 0$.
\label{theo:SVDfinite}
\end{theorem}
\begin{proof}
    Starting from \cref{eq:eig_AAs}, let $r$ be the maximum index for which $\sigma_r > 0$ and define $\v_i=\frac{1}{{\sigma_i}}\A\u_i$ for $i\le r$, so that
\begin{equation}
    \A\u_i={\sigma_i}\v_i.
    \eqnlab{svd_first}
\end{equation}
Substituting \cref{eq:svd_first} into \cref{eq:eig_AAs} gives
\begin{equation}
     \A^*\v_i={\sigma_i}\u_i.
     \eqnlab{svd_second}
\end{equation}
We claim that $\LRc{\sigma_i^2,\v_i}$ for $i\le r$  are the eigenpairs of the self-adjoint operator $\A\A^*$. To see this, applying $\A$ to both sides of \cref{eq:svd_second} to arrive at
\begin{equation}
\A\A^*\v_i={\sigma_i}\A\u_i=\sigma_i^2\v_i.    
\label{eq:eigenAAs}
\end{equation} 
That is, for every eigenpair of $\A^*\A$ corresponding to a non-zero eigenvalue we have an eigenpair of $\A\A^*$ with the same eigenvalue. By the same token, we can show that for every eigenpair of $\A\A^*$ corresponding to a non-zero eigenvalue we have an eigenpair of $\A^*\A$ with the same eigenvalue. As a result, the rest of eigenvalues  of $\A^*\A$ and $\A\A^*$ with indices larger than $r$ must be $0$. The orthonormality of $\{\u_1, \u_2,\hdots, \u_n\}$  and $\{\v_1, \v_2,\hdots, \v_m\}$ is the direct consequence of the spectral decomposition of self-adjoint operators in \cref{coro:spectralFinite}.
The third assertion is clear owing to the first assertion and \cref{eq:eigComplete}.
\end{proof}

\begin{example}[SVD for the closed range theorem, the rank-nullity theorem, and the fundamental theorem of linear algebra]
The SVD decomposition in \cref{theo:SVDfinite} allows us to provide trivial proofs of various important results in finite dimensions including the closed range \cref{theo:CRT}, the rank-nullity theorem and the fundamental theorem of linear algebra. While these results are typically presented for matrices, it is not more difficult to do so for generic linear operators using our general setting as we shall show. To begin, we note that ${\u}_{r+1}, \hdots, {\u}_{n}$ are orthonormal eigenfunctions corresponding to $0$ eigenvalues of $A^*A$. We conclude that
    \[
    {span\{{\u}_{r+1},\hdots, {\u}_{n}\}=\N\LRp{\A^*\A} = \N\LRp{\A},
    }\]
which, together with the fact that $span\{{\u}_{1},\hdots, {\u}_{n}\} = \Xs$ implies
\[
{span\{\A{\u}_{1}, \hdots, \A{\u}_{r}\}=\Range(\A)},
\]
which in turn yield
\[
{span\{\v_1, \hdots, {\v}_{r}\}=\Range(\A)},
\]
since $\A\u_i = \sigma_i\v_i$, $i=1,\hdots,r$.

Similarly, we have
 \[
    {span\{{\v}_{r+1},\hdots, {\v}_{m}\}=\N\LRp{\A\A^*} = \N\LRp{\A^*},
    }\]
and
\[
{span\{\u_1, \hdots, {\u}_{r}\}=\Range(\A^*)}.
\]
With these conclusions in hand, the assertions in the closed range \cref{theo:CRT} and \cref{coro:directSum} are now trivial.

The SVD also provides an obvious proof for the rank-nullity theorem \cite{LinearAlgebraAxler} as
\begin{subequations}
\begin{align}
    \underbrace{\dim(\Xs)}_{=n} &= \underbrace{\dim(\N(\A))}_{=n-r}   + \underbrace{\dim(\Range(\A))}_{=r},
\intertext{and}
     \underbrace{\dim(\Ys)}_{=m} &= \underbrace{\dim (\N(\A^*))}_{=m-r}   + \underbrace{\dim(\Range(\A^*))}_{=r},
\end{align}
\label{eq:RankNullity}
\end{subequations}
which, together with the closed range \cref{theo:CRT}, is the basis for the fundamental theorem of linear algebra \cite{StrangFTL}. This is demonstrated in \cref{fig:FTA} which shows the important role of the four fundamental subspaces $\Range\LRp{\A}, \N\LRp{\A}, \Range\LRp{\A^*}$, and $ \N\LRp{\A^*}$ on the operation of $\A$ and its adjoint $\A^*$. In particular, $\A$ only acts on the range space of $\A^*$, and the results of its action, the range space of $\A$, coincides with the orthogonal complement of the nullspace of $\A^*$. Conversely, $\A^*$ only acts on the range space of $\A$ and the results of its action, the range space of $\A^*$, coincides with the orthogonal complement of the nullspace of $\A$. In other words, the characterization of $\A^*$ completely determines the action of $\A$ and vice versa. We emphasize that  if we remove the dimensions and assume that $\Range\LRp{\A}$ is closed, then the two diagrams in \cref{fig:FTA} also hold for infinite dimensional Hilbert spaces.
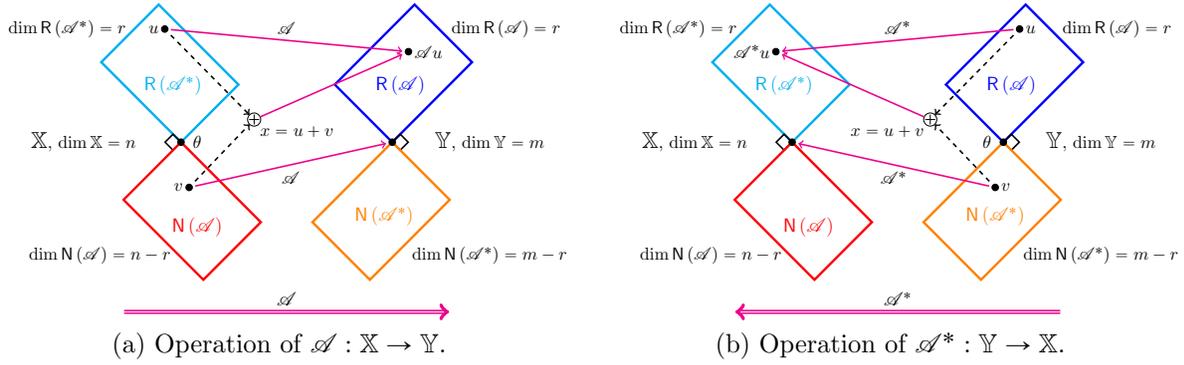
\begin{figure}[h!t!b!]
     \centering
     \begin{subfigure}[b]{0.48\textwidth}
         \centering
         \scalebox{0.6}{
         \begin{tikzpicture}[
    cross/.style={path picture={ \draw[black, shorten <=1pt, shorten >=1pt, line width=.5pt] (path picture bounding box.south) -- (path picture bounding box.north); \draw[black, shorten <=1pt, shorten >=1pt, line width=.5pt] (path picture bounding box.west) -- (path picture bounding box.east);}}
]
    \draw[shift={(1.3*\width,0*\height)}, rotate around={\angle:(0,0)}, color = blue, line width=1.5](0,0) rectangle (\height,\width) {};
    \draw[shift={(1.3*\width,0*\height)}, rotate around={\angle-90:(0,0)}, color = black, line width=1.](0,0) rectangle (0.1*\height,0.1*\height) {};
    \draw[shift={(1.3*\width,0*\height)}, rotate around={\angle:(0,0)}, color = orange, line width=1.5, shift={(-\height,-\width)}](0,0) rectangle (\height,\width) {};

    \draw[shift={(-1.3*\width,0*\height)}, rotate around={-\angle:(0,0)}, color = cyan, line width=1.5, shift={(-\height,0)}](0,0) rectangle (\height,\width) {};
    \draw[shift={(-1.3*\width,0*\height)}, rotate around={\angle+90:(0,0)}, color = black, line width=1.](0,0) rectangle (0.1*\height,0.1*\height) {};
    \draw[shift={(-1.3*\width,0*\height)}, rotate around={-\angle:(0,0)}, color = red, line width=1.5, shift={(0,-\width)}](0,0) rectangle (\height,\width) {};

    \filldraw[black] (-1.3*\width,0*\height) circle (0.04*\width) node[right] () {$\text{ } \theta$};

    \filldraw[black] (1.3*\width,0*\height) circle (0.04*\width) node[] (Right Intersect) {};
    \node[draw=none] at (1.3*\width,0*\height) (Right Intersect) {};
    \filldraw[black] (1.5*\width, 0.8*\height) circle (0.04*\width) node[right] (node b) {$\A u$};
    \node[draw = none] at (1.5*\width, 0.8*\height)  (node b)  {};
    \filldraw[black] (-1.2*\width,-0.4*\height) circle (0.04*\width) node[left] (node vh) {$v$};
    \node[draw = none] at (-1.2*\width,-0.4*\height) (node vh) {};
    \filldraw[black] (-1.5*\width, 1.*\height) circle (0.04*\width) node[left] (node p) {$u$};
    \node[draw = none] at (-1.5*\width, 1.*\height)  (node p)  {};
    \node[draw = none, below right] at (-.4*\width, 0.2*\height) (node x) {$x = u + v$};
    \node[draw, circle, cross, minimum width=.01*\height,line width=.2pt, scale=0.8] at (-.4*\width, 0.2*\height) (node x)  {};

    \node[draw = none] at (-2.5*\width, 0) () {{\Large $\mathbb{X}$}, $\dim\mathbb{X} = n$};
    \node[draw = none] at (-2.7*\width, 1.*\height) () {$\dim{\Range\LRp{\A^*}} = r$};
    \node[draw = none] at (-2.3*\width, -1.*\height) () {$\dim\N\LRp{\A} = n-r$};

    \node[draw = none] at (2.5*\width, 0) () {{\Large $\mathbb{Y}$}, $\dim\mathbb{Y} = m$};
    \node[draw = none] at (2.7*\width, 1.*\height) () {$\dim\Range\LRp{\A} =  r$};
    \node[draw = none] at (2.5*\width, -1.*\height) () {$\dim\N\LRp{\A^*}= m-r$};

    \node[draw = none, color = cyan] at (-1.4*\width, 0.5 * \height) () {\large $\Range\LRp{\A^*}$};
    \node[draw = none, color = red] at (-1.1*\width, -0.75 * \height) () {\large $\N\LRp{\A}$};
    \node[draw = none, color = orange] at (1.2*\width, -0.65 * \height) () {\large $\N\LRp{\A^*}$};
    \node[draw = none, color = blue] at (1.4*\width, 0.5 * \height) () {\large $\Range\LRp{\A}$};

    \draw[->, line width = 1pt, color = magenta](node vh)-- node[below, color = black] {$\A$} (Right Intersect);
    \draw[->, line width = 1pt, color = magenta](node x)--(node b);
    \draw[->, line width = 1pt, color = magenta](node p)-- node[above, color = black] {$\A$} (node b) ;

    \draw[->, dashed, line width = 1pt](node p)--(node x);
    \draw[->, dashed, line width = 1pt](node vh)--(node x);
    
    \draw[->, double, line width = 1pt, color = magenta](-2*\width, -1.5*\height)-- node[above, color = black] {\large $\A$} (2*\width, -1.5*\height);

\end{tikzpicture}
         }
         \caption{Operation of $\A: \Xs \to \Ys$.}
         \label{fig:forward}
     \end{subfigure}
     \hfill
\begin{subfigure}[b]{0.48\textwidth}
         \centering
         \scalebox{0.6}{
         \begin{tikzpicture}[
    cross/.style={path picture={ \draw[black, shorten <=1pt, shorten >=1pt, line width=.5pt] (path picture bounding box.south) -- (path picture bounding box.north); \draw[black, shorten <=1pt, shorten >=1pt, line width=.5pt] (path picture bounding box.west) -- (path picture bounding box.east);}}
]
    \draw[shift={(1.3*\width,0*\height)}, rotate around={\angle:(0,0)}, color = blue, line width=1.5](0,0) rectangle (\height,\width) {};
    \draw[shift={(1.3*\width,0*\height)}, rotate around={\angle-90:(0,0)}, color = black, line width=1.](0,0) rectangle (0.1*\height,0.1*\height) {};
    \draw[shift={(1.3*\width,0*\height)}, rotate around={\angle:(0,0)}, color = orange, line width=1.5, shift={(-\height,-\width)}](0,0) rectangle (\height,\width) {};
    
    \draw[shift={(-1.3*\width,0*\height)}, rotate around={-\angle:(0,0)}, color = cyan, line width=1.5, shift={(-\height,0)}](0,0) rectangle (\height,\width) {};
    \draw[shift={(-1.3*\width,0*\height)}, rotate around={\angle+90:(0,0)}, color = black, line width=1.](0,0) rectangle (0.1*\height,0.1*\height) {};
    \draw[shift={(-1.3*\width,0*\height)}, rotate around={-\angle:(0,0)}, color = red, line width=1.5, shift={(0,-\width)}](0,0) rectangle (\height,\width) {};

    \filldraw[black] (1.3*\width,0*\height) circle (0.04*\width) node[left] () {$\theta \text{ }$};

    \filldraw[black] (-1.3*\width,0*\height) circle (0.04*\width) node[] (Right Intersect) {};
    \node[draw=none] at (-1.3*\width,0*\height) (Right Intersect) {};
    \filldraw[black] (-1.5*\width, 0.8*\height) circle (0.04*\width) node[left] (node b) {$\A^* u$};
    \node[draw = none] at (-1.5*\width, 0.8*\height)  (node b)  {};
    \filldraw[black] (1.2*\width,-0.4*\height) circle (0.04*\width) node[right] (node vh) {$v$};
    \node[draw = none] at (1.2*\width,-0.4*\height) (node vh) {};
    \filldraw[black] (1.5*\width, 1.*\height) circle (0.04*\width) node[right] (node p) {$u$};
    \node[draw = none] at (1.5*\width, 1.*\height)  (node p)  {};
    \node[draw = none, below left] at (.4*\width, 0.2*\height) (node x) {$x = u + v$};
    \node[draw, circle, cross, minimum width=.01*\height,line width=.2pt, scale=0.8] at (.4*\width, 0.2*\height) (node x)  {};

    \node[draw = none] at (-2.5*\width, 0) () {{\Large $\mathbb{X}$}, $\dim\mathbb{X} = n$};
    \node[draw = none] at (-2.7*\width, 1.*\height) () {$\dim{\Range\LRp{\A^*}} = r$};
    \node[draw = none] at (-2.3*\width, -1.*\height) () {$\dim\N\LRp{\A} = n-r$};

    \node[draw = none] at (2.5*\width, 0) () {{\Large $\mathbb{Y}$}, $\dim\mathbb{Y} = m$};
    \node[draw = none] at (2.7*\width, 1.*\height) () {$\dim\Range\LRp{\A} =  r$};
    \node[draw = none] at (2.5*\width, -1.*\height) () {$\dim\N\LRp{\A^*}= m-r$};

    \node[draw = none, color = cyan] at (-1.4*\width, 0.5 * \height) () {\large $\Range\LRp{\A^*}$};
    \node[draw = none, color = red] at (-1.1*\width, -0.75 * \height) () {\large $\N\LRp{\A}$};
    \node[draw = none, color = orange] at (1.2*\width, -0.65 * \height) () {\large $\N\LRp{\A^*}$};
    \node[draw = none, color = blue] at (1.4*\width, 0.5 * \height) () {\large $\Range\LRp{\A}$};

    \draw[->, line width = 1pt, color = magenta](node vh)-- node[below, color = black] {$\A^*$} (Right Intersect);
    \draw[->, line width = 1pt, color = magenta](node x)--(node b);
    \draw[->, line width = 1pt, color = magenta](node p)-- node[above, color = black] {$\A^*$} (node b) ;

    \draw[->, dashed, line width = 1pt](node p)--(node x);
    \draw[->, dashed, line width = 1pt](node vh)--(node x);

    \draw[->, double, line width = 1pt, color = magenta](2*\width, -1.5*\height)-- node[above, color = black] {\large $\A^*$} (-2*\width, -1.5*\height);

\end{tikzpicture}
         }
         \caption{Operation of $\A^*: \Ys \to \Xs$.}
         \label{fig:adjoint}
     \end{subfigure}
     \caption{The fundamental theorem of algebra: four fundamental subspaces $\Range\LRp{\A}, \N\LRp{\A}, \Range\LRp{\A^*}, \N\LRp{\A^*}$, and the operation of $\A$ and $\A^*$ viewed from these subspaces. If we remove the dimensions and assume that $\Range\LRp{\A}$ is closed, then the two diagrams also hold for infinite dimensional Hilbert spaces.}
     \label{fig:FTA}
\end{figure}
\end{example}

Clearly, at the heart of the SVD is the eigenvalue decomposition \cref{eq:eigenAAs}, which could be challenging if it is analytically not tractable on the original operators. In that case, one has to resort to numerical methods. For finite-dimensional settings, an easier path is to explore the matrix representation of the linear operator. 
\begin{corollary}[SVD through matrix representation]
    Consider $\A: \Xs \rightarrow \Ys$ with $\dim(\Xs)=n$ and $\dim(\Ys)=m$, and $\Xb$ and $\Yb$ be orthonormal bases of $\Xs$ and $\Ys$, respectively. Let $\LRc{\sigma_i,\u_i,\v_i}$ be the singular triplets of $\A$ with $1 \le i \le k=min \{ n,\ m\}$ where $\{\u_1, \u_2,\hdots, \u_n\}$ and $\{\v_1, \v_2,\hdots, \v_m\}$ be  orthonormal bases of $\Xs$ and $\Ys$, respectively, given in \cref{theo:SVDfinite}.
    Denote $\ub$ and $\vb$ as the coordinate vectors of $\u$ and $\v$ in the bases $\Xb$ and $\Yb$, respectively, and $\Am$ as the matrix representation of $\A$ with respect to the bases $\Xb$ and $\Yb$. Then $\LRc{\sigma_i,\ub_i,\vb_i}$ be the singular triplets of $\Am$ with
    \begin{enumerate}
    \item $ \Am\ub_i={\sigma_i}\vb_i$ for $i=1,\dots, r$ and $ \Am\ub_i = \theta$ for $i = r+1,\hdots,n$,
    \item $ \Am^*\vb_j={\sigma_j}\ub_j $ for $ j=1,\hdots, r$ and and $ \Am^*\vb_j = \theta$ for $j = r+1,\hdots,m$, where $\Am^*$ is the conjugate tranpose of $\Am$.
    \item $\Am$ is completely determined by its singular triplets in the following sense: for any $\xb \in \Fs^n$, we have 
    \[ 
    \Am\xb=\sum_{i=1}^r \sigma_i \ub_i^*\xb\, \vb_i, \implies \Am= \sum_{i=1}^r \sigma_i \vb_i \ub_i^* ,
    \]
\end{enumerate}
where $r$ is the maximum index for which $\sigma_r > 0$. 

Conversely, if $\LRc{\sigma_i,\ub_i,\vb_i}$, with $1 \le i \le k=min \{ n,\ m\}$,  are the singular triplets of $\Am$, then $\LRp{\sigma_i,\u_i,\v_i}$ are the singular triplets of $\A$.
\label{coro:SVDmatrix}
\end{corollary}
\begin{proof}
   The result is obvious owing to the matrix representation of linear operator and the coordinate vector of a vector in the corresponding bases, \cref{propo:adjointT}, and  the fact that two vectors are orthonormal iff their coordinate vectors in an orthogonal basis are orthonormal (as, e.g., $\LRp{\ub_i,\ub_j}_{\Fs^n} = \LRp{\u_i,\v_j}_{\Xs}$).
\end{proof}

\begin{example}[SVD of matrices]
    For a matrix $\Am: \Xs = \real^n \to \Ys = \real^m$ and if we choose $\Xb$ and $\Yb$ as the canonical bases for $\real^n$ and $\real^m$ with the standard Euclidean inner products, respectively, then the matrix representation of $\Am$ is itself and thus the SVD of $\Am$ is given by \cref{coro:SVDmatrix}. In this case, $\Am^* = \Am^T$. Furthermore: i) $\{\ub_1, \ub_2, \hdots, \ub_n\}$ are orthonormal eigenvectors of $\Am^T\Am$; ii) $\{\vb_1, \vb_2,\hdots, \vb_m\}$ are orthonormal eigenvectors of $\Am\Am^T$; iii) $\sigma_i^2, \ \ i=1,\dots r$, are nonzero eigenvalues of $\Am^T\Am$ or $\Am\Am^T$; and iv) from
    $\Am= \sum_{i=1}^r \sigma_i \vb_i \ub_i^*$
    we can write the full SVD form
    \[ \Am
=\underbrace{\left[
  \begin{array}{cccccccc}
    \vertbar & \vertbar &        & \vertbar&        & \vertbar \\
    \vec{v}_{1}    & \vec{v}_{2}    & \ldots & \vec{v}_{r} & \dots & \vec{v}_{m}\\
    \vertbar & \vertbar &        & \vertbar &        & \vertbar 
  \end{array}
\right]
}_{\Vm}
\underbrace{
\begin{bmatrix}
   \sigma_{1} &  &  & & &\\ 
   & \sigma_{2} &  & & &\\ 
   &  &  \ddots & & &\\ 
   &  &   & \sigma_{r} & &\\
    &  &   & &0& \\ 
     &  &   & &  &\ddots
 \end{bmatrix}
 }_{\Sigma}
 \underbrace{\left[
  \begin{array}{cccccccc}
    \vertbar & \vertbar &        & \vertbar&        & \vertbar \\
    \vec{u}_{1}    & \vec{u}_{2}    & \ldots & \vec{u}_{r} & \dots & \vec{u}_{n}\\
    \vertbar & \vertbar &        & \vertbar &        & \vertbar 
  \end{array}
\right]^T
}_{\Um^T}
 ,
 \]
 that is,
 \[
 \Am = \Vm \Sigma \Um^T,
 \]
 or the reduced SVD form
 \[ \Am
=\underbrace{\left[
  \begin{array}{cccccc}
    \vertbar & \vertbar &        & \vertbar \\
    \vec{v}_{1}    & \vec{v}_{2}    & \ldots & \vec{v}_{r} \\
    \vertbar & \vertbar &        & \vertbar 
  \end{array}
\right]
}_{\Vm_r}
\underbrace{
\begin{bmatrix}
   \sigma_{1} &  &  & \\ 
   & \sigma_{2} &  & \\ 
   &  &  \ddots & \\ 
   &  &   & \sigma_{r}
 \end{bmatrix}
 }_{\Sigma_r}
 \underbrace{\left[
  \begin{array}{cccccc}
    \vertbar & \vertbar &        & \vertbar\\
    \vec{u}_{1}    & \vec{u}_{2}    & \ldots & \vec{u}_{r} \\
    \vertbar & \vertbar &        & \vertbar  
  \end{array}
\right]^T
}_{\Um_r^T},
 \]
 that is,
 \[
 \Am = \Vm_r \Sigma_r \Um_r^T.
 \]
\end{example}

\begin{example}
    Now consider the operator $\A:  \Us =  \text{Span} \{1,\ x,\ x^2\} \subset \Xs = \Ls^2\LRp{-1,1} \to \real^2 $ defined in \cref{exa:quadraticFit}. We are going to find the singular value decomposition of $\A$ indirectly via its matrix representation using \cref{coro:SVDmatrix}. Clearly, two orthonormal bases for $\Us$ and $\real^2$ are $\Xb = \LRc{1,x,\half\LRp{3x^2-1}}$, $\Yb = \LRc{\LRs{1,0}^T,\LRs{0,1}^T}$, respectively. It is a simple exercise to show that the matrix representation $\Am$ of $\A$ in these two bases is given by
    \[
    \Am =\begin{bmatrix}
2& 0 & 1 \\ 
2 & 4/3 & 1/3 
\end{bmatrix}
    \]
    and from \cref{propo:adjointT} we know that the matrix presentation $\Am^*$ of the adjoint $\A^*$ is $\Am^* = \Am^T$. From the proof of \cref{theo:SVDfinite}, by computing the eigendecomposition of $\Am\Am^T$ and $\Am^T\Am$, we can find the full SVD of $\Am$ as
    \[
    \Am = 
    \underbrace{\begin{bmatrix}
-0.6701 & 0.7423  \\
0.7423 & -0.6701  
\end{bmatrix}}_{\Vm} 
\underbrace{\begin{bmatrix}
3.1306 & 0 & 0 \\ 
0 & 1.0433 & 0 
\end{bmatrix}}_{\Sigma}
\underbrace{
\begin{bmatrix}
-0.9023 & 0.1385 & -0.4082 \\ 
-0.3162 & -0.8564 & 0.4082 \\
-0.2931 & 0.4974 & 0.8165
\end{bmatrix}^T
}_{\Um^T}.
    \]
Now, \cref{coro:SVDmatrix} shows that the singular values of $\A$ are $\LRc{3.1306, 1.0433}$ together with the left and right singular functions $\u_i = \ub_i(1) + \ub_i(2) x + \frac{\ub_i(3)}{2}\LRp{3x^2 -1}$, $i = 1,2 ,3$, and $\vb_j$, $j = 1, 2$, respectively.
\end{example}

\subsection{Application of adjoint to optimization with equality constraints}
\label{sect:optimization}
The field of optimization is vast (see, e.g., \cite{Luenberger69,Bertsekas99, NocedalWright06, Bertsekas82b} and the references therein) and we restrict ourselves to unconstrained optimization problems and constrained optimization problems with equality constraints.
We begin our development with unconstrained optimization in finite-dimensional spaces. Let $f: \R^n\ni \xb \mapsto f\LRp{\xb} \in \R$ and we are interested in
studying the optimization problem $\min_{\xb\in\R^N}f\LRp{\xb}$. It is
sufficient to consider the case $n = 1$ as results for optimization problems in  higher dimensions (including infinite dimensions, at least the first order optimality conditions) follow as corollaries. We focus on local optimization problems.
\begin{definition}
$\y$ is a (local) minimizer of $f\LRp{x}$ if there exists a open
  neighborhood, i.e. $\Ball{\delta}{\y}\equaldef\LRc{x: \snor{x - y} <
    \delta }$ (ball with radius $\delta$ in $\R$), for some $\delta > 0$,
  such that
\[
f\LRp{x} \ge f\LRp{y}, \quad \forall x \in \Ball{\delta}{y}.
\]
\label{def:minimum}
\end{definition}
We are interested in finding the necessary and sufficient conditions for $y$ to be a minimizer. To that end, we consider the Taylor remainder theorem \cite{Apostol1967-uv,Kline1972-rf} which states that for twice-differentiable function $f\LRp{x}$ and for any $\epsilon \in \R$, there exists $0 < \theta < 1$ such that 
\begin{equation}
f\LRp{y+\vareps} = f\LRp{y} + \vareps
\underbrace{f'\LRp{y}}_{\text{gradient } g\LRp{y}} +
\half\vareps^2\underbrace{f''\LRp{y+\theta\vareps}}_{\text{Hessian } h\LRp{y+\theta\vareps}}.
\label{eq:TaylorExpansionf}
\end{equation}

If $y$ is a minimizer,  what can we say about the gradient $g\LRp{\cdot}$ and
the Hessian $\h\LRp{\cdot}$ at $y$? We are interested in only necessary conditions. Here is an answer.
\begin{lemma}[First and second order necessary conditions for optimality in $\R$]
\label{lem:necessaryR}
Suppose $f\LRp{x}:\R \to \R$ is twice continuously differentiable in a neighborhood of a minimizer $y$. It
is necessary that 
\begin{itemize}
\item[i)] the gradient vanishes, i.e., $g\LRp{y} = 0$, and
\item[ii)] the Hessian is non-negative, i.e., $\h\LRp{y} \ge 0$.
\end{itemize}
\end{lemma}
\begin{proof}
    We carry out the proof by contradiction. For the first assertion,  we suppose that $\f'(\y) < 0$ and note that  we can pick\footnote{Due to the continuity of $\f''\LRp{\x}$, we can define $M := \max_{\x \in \LRs{\y - L,\y + L}}\snor{\f''\LRp{\x}}$, for some sufficiently large $L > 0$, and then simply pick some $0 < \epsilon < -2\frac{\f'\LRp{\y}}{M}$.} $\epsilon > 0$ such that
    \[
     \epsilon \f'(\y) + \frac{\epsilon^2}{2} \f''(\y + \theta \epsilon) < 0,
    \]
    and together with \cref{eq:TaylorExpansionf} we conclude that $\f\LRp{\y+\epsilon} < \f\LRp{\y}$: a contradiction. A similar contradiction argument can be carry out if $\f'\LRp{\y} > 0$. Thus $\f'\LRp{\y} = 0$.

    For the second assertion, suppose $h\LRp{\y} = \f''\LRp{\y} < 0$. By continuity of $\f''\LRp{\x}$ we can choose sufficiently small $\snor{\epsilon}$ such that $f''\LRp{y+\theta\vareps} < 0$. Then \cref{eq:TaylorExpansionf} reduces to
    \[
    f\LRp{y+\vareps} = f\LRp{y} + 
\half\vareps^2f''\LRp{y+\theta\vareps} < f\LRp{y},
    \]
    which is a contradiction, and this concludes the proof.
\end{proof}
\begin{corollary}
Suppose $f\LRp{\xb}:\R^n \to \R$ is twice continuously differentiable in a neighborhood of a minimizer $\yb$, where $\R^n$ is equipped with the standard Euclidean inner product. It
is necessary that
\begin{itemize}
\item[i)] the gradient vanishes, i.e., $\pp{f}{\xb_i}\LRp{\yb} = 0$, $i=1,\hdots,n$ and
\item[ii)] the Hessian matrix is semi-positive definite, i.e., $\Hm\LRp{\yb} \ge 0$.
\end{itemize}
  \label{coro:optRN}
\end{corollary}
\begin{proof}
We prove the first assertion (the first order optimality condition) as the second one follows similarly. Note that argument is general and will be used again in \cref{lem:firstOrderOpt} to derive the first order optimality condition in general vector spaces.
\[
\begin{array}{cr}
        \yb \text{ is a minimizer of } f\LRp{\xb} & \\
         \Downarrow &     \text{ by definition}   \\
          f\LRp{\xb} \ge f\LRp{\yb}, \quad \forall \xb \in \Ball{\delta}{\yb} := \LRc{\xb \in \R^n: \nor{\xb - \yb}_{\R^n} < \delta} & \\
           \Downarrow &     \text{ pick an arbitrary } \vb   \\
           F\LRp{\vareps} :=f\LRp{\yb + \varepsilon\vb} \ge f\LRp{\yb}, \quad \forall \vareps \in \Ball{\delta/\nor{\vb}_{\R^n}}{0} := \LRc{\vareps: \snor{\vareps} < \delta/\nor{\vb}_{\R^n}} & \\
           \Downarrow &     \text{ by definition}   \\
           0 \text{ is a minimizer of } F\LRp{\varepsilon} & \\
           \Downarrow &     \text{ by \cref{lem:necessaryR}} \\
           \eval{\dd{F}{\vareps}}_{\vareps = 0} = 0 &  \\
           \Downarrow &     \text{ by chain rule} \\
           \sum_{i=1}^n \pp{f}{\xb_i}\LRp{\yb} \vb_i = 0 & \\
           \Downarrow &     \vb \text{ is arbitrary } \\
           \pp{f}{\xb_i}\LRp{\yb} = 0, \quad i=1,\hdots,n. & 
\end{array}  
\]
\end{proof}

Since our goal is to {\em establish the necessary conditions for optimality
that is valid for both finite and infinite dimensional settings},
we present a systematic approach on abstract vector space to accomplish this. {\bf To the end of this section, unless otherwise stated, the results are valid for both finite and infinite dimensional settings}.
We begin
with the notion of the dual space $\Xs^*$ consisting of linear and
bounded functionals on $\Xs$. 
For $\ell\in \Xs^*$ and $\u \in \Xs$, we
use the  standard duality pairing
\[
\LRa{\ell,\u}_{\Xs^*\times \Xs} \equiv  \ell\LRp{\u}
\]
to denote the action of $\ell$ on $u$ (or the evaluation of $\ell$ at
$\u$). \emph{For simplicity in writing, we shall conventionally use $\LRa{\ell,\u}_{\Xs}$ to denote a duality pairing instead of $\LRa{\ell,\u}_{\Xs^*\times \Xs}$.}
The object of interest is nonlinear function 
on a vector space $\Xs$, i.e. \emph{functional}:
\[
f: \Xs \ni \u \mapsto f\LRp{\u} \in \R.
\]
 The classical
derivatives are not
well-defined in this case, and this asks for an extension of
derivatives in vector spaces. Though there are other extensions in the
literature (such as G\^ateaux derivative), let us focus on the Fr\'echet derivative extension (see, e.g., \cite{Luenberger69,ArbogastBona08}).

\begin{definition}
Suppose that there is a linear and bounded map $\D f\LRp{u,\cdot}: \Xs \to \R$
such that
\begin{equation}
{f\LRp{u + \v} =  f\LRp{u} + \D f\LRp{u,\v}} + o\LRp{\nor{\v}_\Xs},
\label{eq:TaylorFrechet}
\end{equation}
where the little-oh notation means 
\[
\lim_{\nor{\v}_\Xs \to 0} \frac{o\LRp{\nor{\v}_\Xs}}{\nor{\v}_\Xs} = 0.
\]
Then $\D f\LRp{u;\cdot}$ is called the Fr\'echet derivative of
the functional $f\LRp{\cdot}$ at $\u$, and \emph{we say $f\LRp{\cdot}$ is Fr\'echet differentiable} at $\u$.
\label{def:Frechet}
\end{definition}

When the Fr\'echet derivative exists, we can compute it conveniently as 
\[
\D f\LRp{\u; \v} =\eval{\dd{f}{t}\LRp{\u+ t\v}}_{t = 0} = \lim_{t \to 0}\frac{f\LRp{\u+ t\v} - f\LRp{\u}}{t}
\] 

For convenience, we use $\D f\LRp{\u}$ to denote the Fr\'echet
derivative $\D f\LRp{\u;\cdot}$ when the argument is irrelevant. It is important
to note that by definition the Fr\'echet derivative $\D f\LRp{\u}$ resides in $\Xs^*$ and thus we interchangeably write it in the duality pairing form
\[
  \D f\LRp{\u;\v} = \LRa{ \D f\LRp{\u}, \v}_{\Xs}.
\]
Due to the linear nature of $\D f\LRp{\u}$, we also write
\[
\D f\LRp{\u}\v :=\LRa{ \D f\LRp{\u}, \v}_{\Xs}.
\]

\begin{definition}[Fr\'echet gradient]
Let $\f: \Xs \to \real$. The gradient of $f\LRp{\cdot}$ at $\u$, denoted as $\Grad f\LRp{\u} \in \Xs$, is defined as a function on $\Xs$ such that
\[
\LRp{\Grad f\LRp{\u}, \v}_\Xs = \D f\LRp{\u;\v} = \LRa{ \D f\LRp{\u}, \v}_{\Xs} = \D \f\LRp{\u}\v, \quad \forall \v \in \Xs.
\]
That is, we define the gradient $\Grad f\LRp{\u} \in \Xs$ as  the Riesz representation of the Fr\'echet derivative $\D f\LRp{\u} \in \Xs^*$.
\label{def:gradient}
\end{definition}
With \cref{def:gradient} at hand, we can identify the gradient of the Fr\'echet derivative and we will explore this fact in many results below.

\begin{example}
Consider $f:\Xs \to \R$ where $\Xs \equiv \R^n$ and
$\R^n$ is endowed with a weighted inner product $\LRp{\xb,\yb}_{{\R^n,\Mm}}
= \xb^T\Mm\yb$ with $\Mm$ being a symmetric positive definite
matrix. Suppose that the (classical) partial derivatives $\pp{f}{\xb_i}$, $i=1,\hdots,n$, of the $f$ are continuous. From \cref{eq:TaylorFrechet}, it is easy to see that the Fr\'echet derivative can be
written as
\[
\D f\LRp{\xb,\hb} = \sum_{i=1}^N\pp{f}{\xb_i}h_i = \LRs{\pp{f}{\xb_1},\hdots,\pp{f}{\xb_n}}^T\hb,
\] 
which, together with \cref{def:gradient}, gives
\[
\Grad f\LRp{\xb} = \Mm^{-1}\LRs{\pp{f}{\xb_1},\hdots,\pp{f}{\xb_n}}^T.
\]
\emph{We observe that the Fr\'echet derivative is a special case of
  the directional derivative, and when $\Mm$ is the identity matrix  usual gradient vector $\LRs{\pp{f}{\xb_1},\hdots,\pp{f}{\xb_n}}^T$ is in fact the Riesz
  representation of the Fr\'echet derivative in the standard Euclidean inner product.}
  \label{ex:gradFinite} 
\end{example}

Of course, the Fr\'echet derivative can be directly generalized to mappings
between two different vector spaces. For example, if $c: \Xs \ni u
\mapsto c\LRp{u} \in \Ys$, then the Fr\'echet derivative $\D c\LRp{u}$, when exists,
can be computed as
\[
\D c\LRp{\u} \v := \D c\LRp{\u;\v} := \lim_{t \to 0}\frac{c\LRp{\u + t\v} - c\LRp{\u}}{t}.
\]
The difference is now that $\D c\LRp{\u}$ is a linear and bounded 
map from $\Xs$ to $\Ys$, that is, $\D c\LRp{\u} \in \Bf\LRp{\Xs, \Ys}$.
\begin{example}
Consider a vector-valued function $\mb{c}\LRp{\xb}:\R^n \to \R^m$ where both $\R^n$ and $\R^m$ are endowed with the standard Euclidean inner products. Applying \cref{ex:gradFinite} for each component of $\mb{c}_i$, $i=1,\hdots,m$ we have
\[
\D \cb\LRp{\xb} \vb = 
\begin{bmatrix}
    \pp{\cb_1}{\xb_1} & \pp{\cb_1}{\xb_2} & \hdots & \pp{\cb_1}{\xb_n} \\
    \pp{\cb_2}{\xb_1} & \pp{\cb_2}{\xb_2} & \hdots & \pp{\cb_2}{\xb_n} \\
    \vdots & \vdots & \hdots & \vdots \\
\pp{\cb_m}{\xb_1} & \pp{\cb_m}{\xb_2} & \hdots & \pp{\cb_m}{\xb_n}
\end{bmatrix}
\vb,
\]
which, together with \cref{def:gradient}, we can define
\[
\Grad\cb\LRp{\xb} := \begin{bmatrix}
    \pp{\cb_1}{\xb_1} & \pp{\cb_1}{\xb_2} & \hdots & \pp{\cb_1}{\xb_n} \\
    \pp{\cb_2}{\xb_1} & \pp{\cb_2}{\xb_2} & \hdots & \pp{\cb_2}{\xb_n} \\
    \vdots & \vdots & \hdots & \vdots \\
\pp{\cb_m}{\xb_1} & \pp{\cb_m}{\xb_2} & \hdots & \pp{\cb_m}{\xb_n}
\end{bmatrix},
\]
which is the Riesz representation of $\D\cb\LRp{\xb}$.
\label{ex:gradFiniteVector}
\end{example}

\begin{lemma}[First order optimality condition for unconstrained optimization]
Suppose that $f:\Xs \to \R$ attains its extremum at $u$. Then it is
necessary
that
\begin{equation}
\D f\LRp{\u}\v = 0, \quad \forall \v \in \Xs,
\label{eq:firstOptInfUnc}
\end{equation}
that is, the (first) variation of $f$ at $\u$ in any ``direction'' $\v$
vanishes. In other words, it is necessary that $\D f\LRp{\u} = 0$ or equivalently
\[
\Grad f\LRp{\u} = 0,
\]
by  the Riesz representation \cref{theo:Riesz}. 
\label{lem:firstOrderOpt}
\end{lemma}
\begin{proof}
It is sufficient to assume that $f$ is minimized at $u$, i.e.,
\[
f\LRp{v} \ge f\LRp{u}, \quad \forall v \in \Ball{\delta}{u} := \LRc{\w \in \Xs: \nor{\w-\u}_\Xs < \delta},
\]
which implies that for any $\v$ such that $\vareps\nor{\v}_\Xs <\delta$, we have
\begin{equation}
f\LRp{u+\vareps \v} \ge f\LRp{u}, \quad \forall \vareps \in \Ball{\delta/\nor{\v}_\Xs}{0}.
\label{eq:infToFiniteOpt}
\end{equation}

If we define $F\LRp{\vareps}\equaldef \f\LRp{\u+\vareps \v}$, then
$F\LRp{\cdot}$ is a function in $\vareps$, namely, $F:\R \ni \vareps \mapsto F\LRp{\vareps} \in
\R$. By \cref{def:minimum}, inequality \cref{eq:infToFiniteOpt} is equivalent to saying that $F\LRp{\cdot}$ attains its minimum at $\vareps =
0$. Thus, from the first result of Lemma \lemref{necessaryR}, we have
\[
\eval{\dd{F}{\vareps}}_{\vareps = 0} = 0,
\]
but this is equivalent to $\D f\LRp{\u,\v} = 0$ by \cref{def:Frechet} of Fr\'echet derivative.
\end{proof}
\begin{example}[First order optimality condition for unconstrained optimization in $\R^n$]
    Back to \cref{ex:gradFinite}. Suppose that $f$ attains its minimum at $\xb$. Combining \cref{eq:firstOptInfUnc}, the gradient found in \cref{ex:gradFinite}, and \cref{def:gradient} yields
    \[
    \Grad f\LRp{\xb} = \bs{0},
    \]
    and thus the first-order necessary condition for optimality is given by
    \[
    \LRs{\pp{f}{\xb_1},\hdots,\pp{f}{\xb_n}}^T = \bs{0}.
    \]
\end{example}

\begin{example}
    We now revisit the least squares problem in \cref{coro:linearLS} in the equivalent form: $\inf_{\x \in \Xs}\half\LRp{\A\x-\y,\A\x-\y}_\Ys$. Using \cref{def:Frechet}, the first order optimality condition \cref{eq:firstOptInfUnc} reads
    \[
    2\LRp{\A\v,\A\x-\y}_\Ys = 0, \quad \forall \v \in \Xs,
    \]
that is,
    \[
    \A^*\A\x = \A^*\y,
    \]
    which is consistent with the least squares solution in \cref{eq:LSsolution}.
\label{exa:LSsolutionDerivative}
\end{example}

Up to this point, we have looked at unconstrained optimization problems
and derived the (first order) necessary condition for optimality. We
next discuss optimality conditions for constrained optimization. Let
us consider the following constrained optimization problem
\[
\min_{u\in \Xs} f\LRp{u}, \quad \SubjectTo c\LRp{u} = \theta, \text{ where } c\LRp{\cdot}: \Xs \to \Ys.
\]
If there were no constraint $c\LRp{u} = \theta$, then from \cref{lem:firstOrderOpt} the optimality condition would be
\[
\D f\LRp{\u}\v = 0, \quad \forall \v \in \Xs,
\]
That is, the variation of $f$ at $\u$ in any ``direction'' $\v$
vanishes. However, $\v$ can be no longer arbitrary since the constraint
must be satisfied at $\u + t\v$ for any small $t$. In other words, $\u + t\v$ needs to be \emph{feasible}, i.e.,
\[
c\LRp{\u + t\v} = 0, \quad \text{ for any feasible } \u + t\v.
\]
Suppose $c$ is Fr\'echet differentiable, it is therefore necessary that
\[
\D c\LRp{\u}\v = 0, \text{ for any feasible } \u + t\v.
\]
To rigorously establish this result, we need the inverse function theorem \cite{Luenberger69}, which in turn is a direct consequence of the implicit function theorem \cite{Lang2001-pt,Edwards1995-bq}.
\begin{theorem}[Implicit function theorem]
    Let $c: \Xs\times \Zs \to \Ys$ be continuously Fr\'echet differentiable and $\D_\u c\LRp{\u,\z}: \Xs \to \Ys$, is invertible at a point $\LRs{\u_0,\z_0}^T$, at which $c\LRp{\u_0,\z_0} = 0$. Then, there exist a neighborhood $\Ball{\delta}{\z_0}$ and a continuously Fr\'echet differential function $g: \Ball{\delta}{\z_0} \to \Xs$ such that $c\LRp{g\LRp{\z},\z} = 0$  for all $\z \in \Ball{\delta}{\z_0}$.
    \label{theo:implicit}
\end{theorem}
\begin{theorem}[Inverse function theorem]
Let $f: \Xs \to \Zs$. Assume that $\D f\LRp{\u_0}$ is continuous and maps
$\Xs$ {\bf onto} $\Zs$. Then, there is a neighborhood $\Ball{\delta}{f\LRp{\u_0}}$ of $f\LRp{\u_0}$
such that $f\LRp{\u} = z$ has a unique continuously differentiable solution $\u(\z)$ for
every $z \in \Ball{\delta}{f\LRp{\u_0}}$.
\label{theo:inverse}
\end{theorem}
\begin{proof}
    The result is clear if we define $c\LRp{\u,\z} = \z - \f\LRp{\u}$ and set $\z_0 = \f\LRp{\u_0}$. Then by the implicit function \cref{theo:implicit}, there exists $g: \Ball{\delta}{\z_0} \to \Xs$ such that $0 = c\LRp{\g\LRp{\z},\z} = \z - \f\LRp{\g\LRp{\z}} = 0$ for all $\z \in \Ball{\delta}{\z_0}$. Setting $\u = \g\LRp{\z}$ concludes the proof.
\end{proof}
\begin{lemma}[First order optimality condition for equality constraints]
Suppose $f: \Xs \to \R$ attains its extremum at $u_0$ subject to the constraint
$c\LRp{u} = 0$, where $c:\Xs \to \Ys$. Assume that both $f$ and $c$ are continuously
Fr\'echet differentiable in an open set containing $u_0$, and $\D
c\LRp{u_o}$ maps $\Xs$ {\bf onto} $\Ys$. Then, it is necessary that
\[
\D f\LRp{\u_0}\v = 0, \quad \forall \v \in \Xs \text{ such that } \D c\LRp{\u_0}\v = \theta,
\]
or equivalently
\[
\LRp{\Grad f\LRp{\u_0},\v}_\Xs = 0, \quad \forall \v \in \Xs \text{ such that } \D c\LRp{\u_0}\v = \theta.
\]
\label{lem:firstOptFunc}
\end{lemma}
\begin{proof}
We follow closely the proof by contradiction in \cite[Lemma 1 of Chapter 9]{Luenberger69}. Without lost of generality, asssume $u_0$ is a minimizer.  Let us consider the transformation $g\LRp{u} =
\LRp{f\LRp{u}, c\LRp{u}}: \Xs \to \R\times\Ys$. Assume that there
exists $h$ such that $\D c\LRp{u_0; h} = 0$ but $\D f\LRp{u_0;h} \ne
0$. Then the function $\LRp{\D f\LRp{u_0}, \D c\LRp{u_0}}$ maps $\Xs$ onto $\R
\times \Ys$ since $\D c\LRp{u_0}$ maps $\Xs$ onto $\Ys$. By the
  inverse function \cref{theo:inverse} there exists $\vareps$ and $u$ with
  $\nor{u-u_0}_\Xs < \vareps$ such that $g\LRp{u} =
  \LRp{f\LRp{u_0}-\delta, 0}$ for some small $\delta > 0$. Thus, $f\LRp{u} = f\LRp{u_0} - \delta < f\LRp{u_0}$: a contradiction. 
\end{proof}

\begin{example}
Let $f : \R^n \to \R$ be continuously differentiable and $\Am \in \R^{m\times n}$, where $\R^n$ and $\R^m$ are endowed with the standard Euclidean inner products. We consider the following problem
\[
\min_{\xb\in \R^n} f\LRp{\xb}, \quad \SubjectTo \Am\xb = \bb.
\]
From
\cref{lem:firstOptFunc}, \cref{ex:gradFinite}, and \cref{ex:gradFiniteVector} we can write the the first order optimality condition as
\[
\Grad f^T\LRp{\xb}\vb = 0, \quad \forall \vb \in \R^n \text{ such that } \Am\vb = \mb{0},
\]
or equivalently
\[
\Grad f^T\LRp{\xb}\vb = 0, \quad \forall \vb \in \N\LRp{\Am},
\]
i.e., \emph{due to the constraint, the gradient of $f$ at an optimum $\xb$ does not vanish but is orthogonal to the nullspace of the gradient $\Am$ of the constraint.}  In other words, for constrained optimization problems, at an optimum the projection of the gradient of the objective function in the nullspace of the gradient of the constraints vanishes. If we define $\Zm$ with columns comprising a basis of the nullspace of $\Am$,
then $\vb = \Zm\mb{r}$ for some vector $\mb{r}$ whose dimension is
the dimension of the nullspace. As a result, the constraint is completely eliminated and the first order optimality
condition  now reads
\[
\Grad f^T\LRp{\xb}{\Zm} = 0.
\]
Note that $\mb{g}_r\LRp{\xb} \equaldef \Grad f^T\LRp{\xb}{\Zm}$\textemdash the coordinates of the gradient in the nullspace of the constraint gradient\textemdash is known as the
reduced gradient \cite{Frank1956,Abadie69}. The reduced gradient is nothing more than the total gradient of the objective function with respect to the reduced variable $\mb{r}$ as we show below.  The proof for an important class of constrained optimization is presented in  \cref{rema:adjoint}. One of the reasons for its name is that its
dimension is smaller than the dimension $n$ of the original gradient vector $\Grad\f\LRp{\xb}$. \emph{Another important point one can draw from the
  optimality condition for the reduced gradient is that in the reduced
  optimization variables $\mb{r}$, the optimization problem
  becomes implicitly unconstrained} (see the explicit transformation to the reduced space at the end of the example). Now from  \cref{coro:optRN} we know that
the derivative of the reduced gradient $\mb{g}_r$, namely the reduced
Hessian $\Hm_r$, is necessary to be semi-positive definite in any
direction in the reduced space. By the chain rule we have
\[
\mb{r}^T\Hm_r\LRp{\xb}\mb{r} = \eval{\dd{\mb{g}_r\LRp{\xb +
      t{Z}\mb{r}}}{t}}_{t=0} = \mb{r}^T{\Zm}^T\Grad f^2\LRp{\xb}{\Zm}\mb{r}, \quad \forall \mb{r},
\]
from which it follows that
\[
{\Hm}_r\LRp{\xb} = {\Zm}^T\Grad f^2\LRp{\xb}{\Zm}.
\]
Note that if QR factorization of $\Am$ is feasible, then ${\Zm}$ can be easily found. ${\Zm}$ can also be explicitly identified for the case when  $m \le n $ and $\Am$
has linearly independent rows. Indeed, up to a permutation of columns, we can rewrite $\Am$ as
\[
\Am = \LRs{{\Um}, {\Vm}},
\]
where ${\Um} \in \R^{m\times m}$ is a nonsingular matrix. Then ${\Zm}$ can be written as
\[
{\Zm} = \LRs{
\begin{array}{c}
-{\Um}^{-1}{\Vm} \\
\Im
\end{array}
},
\]
where $\Im$ is the $\LRp{n-m}\times\LRp{n-m}$ identity matrix. In this case, we can write $\xb = \LRs{\xb_{\Um},\xb_{\Vm}}^T$ and we can eliminate $\xb_{\Um}$ from the constraint as $\xb_{\Um} = \Um^{-1}\bb -\Um^{-1}\Vm\xb_{\Vm}$. The reduced optimization variable is thus $\xb_\Vm$ and the original constraint optimization problem is now unconstrained with respect to $\xb_\Vm$.
\label{exa:linearConstraint}
\end{example}

In order to provide further insights and  make the optimality condition practical for large-scale
computation we need a Lagrangian formalism, and {\em this is where the adjoint plays the key role}. Thanks to the closed range \cref{theo:CRT}, the Lagrangian multiplier theorem \cite{Luenberger69,Bertsekas99, NocedalWright06, Bertsekas82b} is a straightforward equivalence to \cref{lem:firstOptFunc}.

\begin{theorem}[Lagrangian multiplier theorem]
Assume that $f: \Xs \to \R$ is continuously Fr\'echet differentiable
and it attains the extremum at $u_0$ subject to the constraint
$c\LRp{u} = 0$, where $c: \Xs \to \Ys$ is continuously Fr\'echet
differentiable. Suppose that $\D c\LRp{u_0}$ maps $\Xs$ {\bf onto}
$\Ys$. Then, there exists an element $y \in \Ys$ such that the
following Lagrangian functional
\[
L\LRp{u} \equaldef f\LRp{u} + \LRp{y,c\LRp{u}}_\Ys
\]
is stationary at $u_0$, i.e.,
\begin{equation}
\D L\LRp{u_0}\h = \D f\LRp{u_0}\h + \LRp{y,\D c\LRp{u_0}\h}_\Ys = 0, \quad \forall h \in \Xs,
\label{eq:firstOptInf}
\end{equation}
or equivalently
\begin{equation}
\Grad L\LRp{u_0} = \Grad f\LRp{u_0} + \LRs{\D c\LRp{u_0}}^*y = 0.
\label{eq:firstOptInfGrad}
\end{equation}
\label{theo:LagrangeMult}
\end{theorem}
\begin{proof}
From \cref{lem:firstOptFunc} we see that $\Grad f\LRp{u_0}$
is orthogonal to the nullspace of $\D c\LRp{u_0}$. Since $\D
c\LRp{u_0}$ maps $\Xs$ onto $\Ys$,   the range of $\D
c\LRp{u_0}$
closed and  the closed range 
\cref{theo:CRT}  gives
\[
\Grad f\LRp{u_0} \in \Range\LRp{\LRs{\D c\LRp{u_0}}^*},
\]
which implies that there exists $y \in \Ys$ such that
\[
\Grad f\LRp{u_0} = -\LRs{\D c\LRp{u_0}}^*y,
\]
and this ends the proof.
\end{proof}

\begin{remark}
    The appealing feature of the Lagrangian approach in \cref{theo:LagrangeMult} is that the first order optimality condition \cref{eq:firstOptInfGrad} is the standard optimality condition for unconstrained problem in \cref{lem:firstOrderOpt}, but for the Lagrangian instead of the original objective function $f$. The key implication in the Langrangian formalism is thus the optimization problem is unconstrained in the original optimization variable $\u$ plus the Lagrange multiplier $\y$, as far as the first order optimality condition is concerned. \emph{The Lagrangian approach is also known as the adjoint approach as it involves the adjoint of the gradient of the constraint in the first order optimality condition \cref{eq:firstOptInfGrad}}. If further structures of the constraints and/or optimization variables are given, the Lagrangian approach can lead to an efficient reduced space approach as we will discuss in the below examples, including the derivation and insights into backpropagation of neural network in \cref{sect:DNN}.
\end{remark}

\begin{example}
    Back to \cref{exa:linearConstraint}. Applying the Lagrangian multiplier \cref{theo:LagrangeMult} together with the Riesz representation \cref{theo:Riesz}, the equivalent first order optimality condition \cref{eq:firstOptInfGrad} reduces to
    \[
    \Grad f\LRp{\xb} + \Am^T\yb = 0,
    \]
    which says that {\em at an optimum of a constrained optimization problem, the gradient of the objective function $f\LRp{\xb}$ does not vanish but is a linear combination of the gradient of the constraints}. Again, this is the same as saying that the gradient of the objective function $f\LRp{\xb}$ at an optimum is orthogonal to the nullspace of the gradient of the constraints.
\end{example}

\begin{example}
Let $f : \R^n \to \R$ and $\mb{c}\LRp{\xb}:\R^n \to \R^m$, where $\R^n$ and $\R^m$ are endowed with the standard Euclidean inner products. We consider the following optimization problem
\[
\min_{\xb\in \R^n} f\LRp{\xb}, \quad \SubjectTo \mb{c}\LRp{\xb} = \bb.
\]
From
\cref{lem:firstOptFunc}, \cref{ex:gradFinite}, \cref{ex:gradFiniteVector}, and the Lagrangian multiplier \cref{theo:LagrangeMult} we can write the first order optimality condition either as
\begin{equation}
\label{eq:NonFirst}
\Grad f\LRp{\xb} +  \Grad \mb{c}^T\LRp{\xb}\yb = \bs{0},
\end{equation}
or as
\[
\mb{g}\LRp{\xb} \equaldef \Grad f^T\LRp{\xb}\vb\LRp{\xb} = 0, \quad \forall \vb\LRp{\xb} \in \R^n \text{ such that } \Grad \mb{c}\LRp{\xb}\vb\LRp{\xb} = \bs{0}.
\]
Note that unlike the linear problem, $\vb$ is a function of $\xb$
since the nullspace of $\Grad \mb{c}\LRp{\xb}$ depends on $\xb$. As a
result, \emph{the reduced Hessian is different from
  that of the linear constraint counterpart in \cref{exa:linearConstraint} as we now show}. To that end, we compute the Fr\'echet derivative of $\mb{g}$ in any
direction $\mb{p}$ in the reduced space $\N\LRp{\Grad \mb{c}\LRp{\xb}}$:
\begin{equation}
\eqnlab{NonSecond}
\mb{p}^T\Hm\mb{h} = \mb{p}^T\Grad f^2\vb + \Grad f^T\D\vb\LRp{\xb, \mb{p}}.
\end{equation}
To compute $\D\vb\LRp{\xb, \mb{p}}$ we take the Fr\'echet derivative both sides of
$\Grad \mb{c}\LRp{\xb}\vb\LRp{\xb} = \mb{0}$, row by row, in direction $\mb{p}$ to arrive at
\begin{equation}
\eqnlab{NonCon}
\mb{p}^T\Grad^2c_i\vb + \Grad c_i^T\D\vb\LRp{\xb, \mb{p}} = 0.
\end{equation}
Combining \cref{eq:NonFirst}--\cref{eq:NonCon} gives
\[
\mb{p}^T\Hm\mb{h} = \mb{p}^T\LRp{\Grad f^2 -\sum_{i=1}^m \Grad^2c_i y_i}\vb.
\]
Since both $\mb{p}$ and $\vb$ belong to the nullspace
$\Zm\LRp{\xb}$ of $\Grad \mb{c}\LRp{\xb}$, the reduced Hessian $\Hm_r$ is then given by
\[
\Hm_r = \Zm^T\LRp{\Grad f^2 -\sum_{i=1}^m \Grad^2c_i y_i}\Zm.
\]
As can be seen, the Hessians of the
constraints (which is zero for linear constraint case in \cref{exa:linearConstraint}) contribute to the reduced Hessian. Unlike \cref{exa:linearConstraint} in which the reduced space, and hence the reduced optimization variables, is fixed if an iterative gradient-based algorithm is employed, the reduced space for this example changes at each optimization step due the nonlinear nature of the constraint $\mb{c}\LRp{\xb} = \bb$. 
\end{example}

To the end of this section, we consider an important class of constrained optimization problems in which the constraints are equalities and the optimization variables are separable in the sense that from the constraint one can solve for one sub-variable as a function of the other.   As we shall see in \cref{sect:DNN}, training deep neural network with backpropagation is a special case of this class. PDE-constrained optimization problem is another special case as shown in \cref{sect:PDEconstrained}.

\begin{corollary}[Optimization with special equality constraint]
 Consider optimization problems that can be expressed in the following form
\[
\min_{\u\in \Xs, \z \in \Zs} f\LRp{\u,\z}, \quad \SubjectTo c\LRp{\u,\z} =
0, \text{ where } c\LRp{\cdot,\cdot}: \Xs\times\Zs \to \Ys,
\]
where \emph{the Fr\'echet derivative of the constraint with respect to
  $\u$, i.e. $\D_\u c\LRp{\u,\z}: \Xs \to \Ys$, is invertible} at an optimum $\LRs{\u_0,\z_0}^T$. The first
order optimality condition \cref{eq:firstOptInfGrad}, together with the constraint, can be written as
\begin{subequations}
\begin{align}
\eqnlab{forward}
c\LRp{\u_0,\z_0} &= 0, & \text{ Forward equation}, \\
\eqnlab{adjointEq}
\Grad_\u f\LRp{\u_0,\z_0} + \LRs{\D_\u c\LRp{\u_0,\z_0}}^*y &= 0, & \text{ Adjoint equation},\\
\eqnlab{control}
\Grad_\z f\LRp{\u_0,\z_0} + \LRs{\D_\z c\LRp{\u_0,\z_0}}^*y &= 0, & \text{ Control equation},
\end{align}
\label{eq:KKTinf}
\end{subequations}
\label{coro:adjointReduced}
\end{corollary}
where $\Grad_\u$ and $\Grad_\z$ denote the Fr\'echet derivative with respect to $\u$ and $\z$, respectively.
\begin{proof}
    The proof is a straightforward application of \cref{eq:firstOptInfGrad} to the group optimization variable $\LRs{\u,\z}^T$, and thus $\Grad f = \LRs{\Grad_\u f, \Grad_\z f}^T$ and $\D c = \LRs{\D_\u c, \D_\z c}^T$.
\end{proof}
Note that the left-hand side (LHS) of the forward problem \cref{eq:forward} is simply the derivative of the Lagrangian with respect to the adjoint variable $\y$. The LHS of the adjoint equation \cref{eq:adjointEq} is nothing more than the derivative of the Lagrangian with respect to $\u$, and the LHS of the control equation \cref{eq:control} is the derivative of the Lagrangian with respect to $\rr$. For this class of optimization problems, we can eliminate both the ``state" variable $\u$ and adjoint variable $\y$ so that the optimization problem is genuinely unconstrained in only the control variable $\z$ around a neighborhood of the optimum $\LRs{\u_0,\z_0}^T$. Indeed, from the implicit function \cref{theo:implicit} there exists $\Ball{\delta}{\z_0}$ and a continuously differentiable function $\g: \Ball{\delta}{\z_0} \to \Xs$ such that $\u = \g\LRp{\z}$ solve the constraint $c\LRp{\u,\z} = 0$ for any $\z \in \Ball{\delta}{\z_0}$. The objective function becomes $f\LRp{\g\LRp{z},\z}$, and thus a function of only $\z$, for all $\z \in \Ball{\delta}{\z_0}$. The optimization variable is reduced to only $\z$ and this is known as the {\em reduced space approach} \cite{Frank1956,Abadie69}. Clearly, we do not know $\u = \g\LRp{\z}$ explicitly, and the question is how to compute the reduced gradient $\Grad\f\LRp{\g\LRp{\z},\z}$, which is needed for any gradient-based approach? Note that by $\Grad\f\LRp{\g\LRp{\z},\z}$ we mean the total derivative with respect to $\z$.

\begin{lemma}[Reduced gradient in constrained optimization with equality constraints]
    With the same setting as in \cref{coro:adjointReduced}, there exists a neighborhood $ \Ball{\delta}{\z_0}$ such that the reduced gradient at any $\z \in \Ball{\delta}{\z_0}$ is given by
    \begin{equation}
    \Grad\f\LRp{\g\LRp{\z},\z} = \Grad_\z f\LRp{\u,\z} + \LRs{\D_\z c\LRp{\u,\z}}^*y,
    \label{eq:abstractReducedGradient}
    \end{equation}
    where $\u = \g\LRp{\z}$ and $\y$ satisfy the following forward and adjoint equations:
\begin{subequations}
    \begin{align}
    \label{eq:reducedForward}
c\LRp{\u,\z} &= 0, & \text{ Forward equation}, \\
\label{eq:reducedAdjoint}
\Grad_\u f\LRp{\u,\z} + \LRs{\D_\u c\LRp{\u,\z}}^*y &= 0, & \text{ Adjoint equation}.
\end{align}
\label{eq:reduced}
\end{subequations}
\label{lem:reducedGradient}
\end{lemma} 
\begin{proof}
Note that the invertibility of a linear operator $\A$ implies the invertibility\footnote{
    Indeed, suppose $\A: \Xs \to \Ys$ is invertible, then $\LRa{\u,\A^*\LRp{\A^{-1}}^*\w}_\Xs = \LRa{\A\u,\LRp{\A^{-1}}^*\w}_\Ys = \LRa{\A^{-1}\A\u,\w}_\Xs = \LRa{\u,\w}_\Xs$ for any $\u, \w \in \Xs$. Thus, $\LRp{\A^*}^{-1} = \LRp{\A^{-1}}^*$.
    } of its adjoint $\A^*$ with $\LRp{\A^*}^{-1} = \LRp{\A^{-1}}^*$.
For simplicity, we use $\LRp{\A}^{-*}$ to denote $\LRp{\A^{-1}}^*$.
We have
\[
\begin{array}{cr}
\LRp{\Grad\f\LRp{\g\LRp{\z},\z},\h}_\Zs &  \\
\verteq{} & \text{Gradient \cref{def:gradient}}\\
\LRa{\D\f\LRp{\g\LRp{\z},\z},\h}_\Zs &  \\
\verteq{}  & \text{Chain rule}
\\
\LRa{\D_\z f\LRp{\u,\z},\h}_\Zs + \LRa{\D_\u f\LRp{\u,\z},\D_\z\g\LRp{\z}\h}_\Xs & \\
\verteq{} & \text{Gradient \cref{def:gradient}}\\
\LRp{\Grad_\z f\LRp{\u,\z},\h}_\Zs + \LRp{\Grad_\u f\LRp{\u,\z},\D_\z\g\LRp{\z}\h}_\Xs & \\
\verteq{}       &  \text{Derivative of the constraint} \\
\LRp{\Grad_\z f\LRp{\u,\z},\h}_\Zs - \LRp{\Grad_\u f\LRp{\u,\z},\LRs{ \D_\u c\LRp{\u,\z}}^{-1}\D_\z c\LRp{\u,\z}\h}_\Xs  & \\
      \verteq{}       & \text{Adjoint \cref{defi:adjoint}} \\
\LRp{\Grad_\z f\LRp{\u,\z},\h}_\Zs - \LRp{\LRs{\D_\z c\LRp{\u,\z}}^*\LRs{ \D_\u c\LRp{\u,\z}}^{-*}\Grad_\u f\LRp{\u,\z},\h}_\Zs  & \\
      \verteq{}       & \LRs{ \D_\u c\LRp{\u,\z}}^{*} \y := -\Grad_\u f\LRp{\u,\z}\\
\LRp{\Grad_\z f\LRp{\u,\z},\h}_\Zs + \LRp{\LRs{\D_\z c\LRp{\u,\z}}^*y,\h}_\Zs  &
\end{array} 
\]
where, as in the first equality, the derivative of the constraint in the third equality is given by the chain rule\footnote{The chain rule for Fr\'echet deriviation can be derived from \cref{eq:TaylorFrechet}. Let $\f: \Xs \ni \u  \mapsto \f\LRp{\u} \in \Ys$ and $\g: \Zs \ni \z \mapsto \g\LRp{\z} \in \Xs$. We have
$\f\LRp{\g\LRp{\z + \epsilon\h}} = \f\LRp{\g\LRp{\z} + \epsilon\D\g\LRp{\z}\h + o\LRp{\epsilon}} = \f\LRp{\g\LRp{\z}} + \epsilon\D_\u\f\LRp{\g\LRp{\z}}\LRp{\D\g\LRp{\z}\h + \epsilon^{-1}o\LRp{\epsilon}} + o\LRp{\epsilon}$. Thus, $\D_z\f\LRp{\g\LRp{\z}}\h = \lim_{\epsilon \to 0} \frac{\f\LRp{\g\LRp{\z + \epsilon\h}}  - \f\LRp{\g\LRp{\z}}}{\epsilon} = \D_\u\f\LRp{\g\LRp{\z}}\D\g\LRp{\z}\h = \LRa{\D_\u\f\LRp{\g\LRp{\z}}, \D\g\LRp{\z}\h}_\Xs$}: 
\begin{align*}
    \D_\z c\LRp{\u,\z} + \D_\u c\LRp{\u,\z}{\D_z\g\LRp{\z}}&= \theta,
\end{align*}
which, due to the invertibility of $\D_\u c\LRp{\u,\rr}$, allows us to solve for ${\D_z\g\LRp{\z}}$.
\end{proof}

\begin{remark}
Note that in practice, approximating the minimum $\u_0$ is typically done using a gradient descent algorithm and \cref{lem:reducedGradient} shows that the reduced gradient  
can be computed in each iteration for a given $\rr$ via three steps: 1) solve the
\emph{forward equation} \cref{eq:reducedForward} for $u\LRp{\rr}$, 2) solve
the \emph{adjoint equation} \cref{eq:reducedAdjoint} for $y\LRp{u(\rr),\rr}$, and 3) 
substitute $\u(\z)$ and $\y(\u\LRp{\z},\z)$ into \cref{eq:abstractReducedGradient} 
to obtain the
reduced gradient. Moreover, the adjoint equation is always linear in the adjoint variable $\y$ regardless the linear or nonlinear nature of the forward equation. Note that full space iteration  based on the first order optimality condition \cref{eq:KKTinf} is also possible and can be consulted from \cite{NocedalWright06} and the references therein.

    
    \label{rema:adjoint}
\end{remark}
\begin{example}
    To appreciate the adjoint approach, let us apply \cref{lem:reducedGradient} to identify the forward equation, the adjoint equation, and the reduced gradient of the following finite dimensional optimization problem
\[
\min_{\ub\in \R^n, \rrb \in \R^p} f\LRp{\ub, \rrb}, \quad \SubjectTo \mb{c}\LRp{\ub, \rrb} = \bs{0},
\]
where $f : \R^n\times \R^p \to \R$ and $\mb{c}\LRp{\ub,
  \rrb}:\R^n\times\R^p \to \R^n$, and all spaces are equipped with the standard Euclidean inner products. We assume that
$\det\LRp{\Grad_{\ub}\mb{c}} \ne 0, \forall \ub, \rrb$ so that the
implicit function theorem allows us to compute $\ub$ as a function of $\rrb$ from the constraint. Applying \cref{lem:reducedGradient} the reduced gradient reads
\begin{equation}
    \Grad\f = \Grad_{\rrb} f\LRp{\ub_0,\rrb_0} + \LRs{\Grad_{\rrb} \mb{c}\LRp{\ub_0,\rrb_0}}^T\yb,
    \label{eq:controlFinite}
\end{equation}
where $\ub$ and $\yb$ are computed from
\begin{subequations}
\eqnlab{KKTfinite}
\begin{align}
  \eqnlab{forwardFinite}
  \mb{c}\LRp{\ub_0,\rrb_0} &= \mb{0}, & \text{ Forward equation}, \\
  \eqnlab{adjointFinite}
  \Grad_{\ub} f\LRp{\ub_0,\rrb_0} + \LRs{\Grad_{\ub} \mb{c}\LRp{\ub_0,\rrb_0}}^T\yb &= \mb{0}, & \text{ Adjoint equation}.
\end{align}
\end{subequations}
\label{exa:adjointMethodFinite}
\end{example}

\subsection{Application of adjoint to backpropagation in deep learning}
\label{sect:DNN}

In this section, we consider standard fully connected deep neural
network and use the adjoint method in \cref{exa:adjointMethodFinite} to
derive the backpropagation method for computing the gradient of the
loss function with respect to the weights and biases of a general fully-connected deep neural work (DNN). Excellent review papers on deep learning can be found in \cite{LeCun2015,Schmidhuber2015}, and the history of back-propagation can be traced back to \cite{linnainmaa1970representation,Linnainmaa1976}. The gradient is needed for gradient-based methods (see, e.g., \cite{NocedalWright06} and the references therein) such as stochastic gradient descent \cite{Kiefer1952,Robbins1951, ShapiroDentchevaRuszczynski09}. The extension of the adjoint method for other type of neural networks such as ResNet \cite{He2016} and CNN \cite{Fukushima1980,Hubel1968} are straightforward. We are going to show that {\em the backpropagation is nothing more than a reduced space approach to compute the gradient using adjoint method.}

\begin{definition}[$L$-layer Neural network]
  Given $n_\ell,s_0, s_1, \hdots, s_{n_\ell} \in \mathbb{N}$, an $n_\ell$-layer  neural network is defined as the following series of composition
  \begin{equation}
\begin{aligned}
  &\text{Input layer}:  \ab^0 - \xb = \bs{0}, \\
  & \text{The $i$th layer}:   \ab^i - \sigma \LRp{\Wm^i\ab^{i-1}+\bb^{i}} = \bs{0}, \quad i = 1,\hdots,n_\ell,
\end{aligned}
\eqnlab{DNN}
\end{equation}
 where $\xb \in \real^{s_0}$; $\Wm^i \in \real^{s_i} \times \real^{s_{i-1}}$ and  $\bb^i \in \real^{s_i}$, $i = 1, \hdots,n_\ell$,  are weight matrix and bias vector of the $i$th layer;  $\ab^i \in \real^{s_{i}}$ is the output of the $i$th layer; and the activation function, $\sigma$, {\bf acts component-wise} when its argument is a vector.
\end{definition}

Let us define $\ub := \LRs{\ab^0, \hdots, \ab^{n_\ell}}^T$, $\rrb := \LRs{\Wm^1, \bb^1,
  \hdots,\Wm^{n_\ell}, \bb^{n_\ell}}^T$, and $\bs{c}\LRp{\ub,\rrb} = \bs{0}$ as the concatenation of all the sub-equations in \cref{eq:DNN}.
For concreteness, let us consider the loss (objective) function to be:
\[
f\LRp{\ub, \rrb}=\frac{1}{2} \nor{{\ab}^{obs}-\ab^{n_\ell}}^2,
\]
where $\ab^{obs}$ is a given data (label). The neural network training
problem is exactly the constrained optimization
problem in \cref{exa:adjointMethodFinite}. Thus
\begin{itemize}
\item The forward equation \cref{eq:forwardFinite}, by definition, are nothing more than the neural network description in \cref{eq:DNN}. For example, the $i$th block of forward sub-equations (the $i$th layer equation) are 
\[
\bs{c}^i\LRp{\ub,\rrb} = \ab^i - \sigma \LRp{\Wm^i\ab^{i-1}+\bb^{i}} = \bs{0},
\]
the corresponding $i$th block of forward solution is $\ub^i = \ab^i$, and the $i$th block of parameter is $\rrb^i = \LRs{\rrb_\Wm^i, \rrb_{\bb}^i}^T := \LRs{\Wm^i, {\bb}^i}^T$.
Clearly, the Jacobian $\Grad_{\ub}\bs{c}$ is a lower block bi-diagonal matrix with identity blocks on the diagonal, and is thus invertible for all $\ub$ and $\rrb$. Consequently, all results in \cref{exa:adjointMethodFinite} hold. 
\item To unfold the adjoint equation \cref{eq:adjointFinite}, we note that the whole adjoint vector $\yb$ is the concatenation of the adjoint sub-vector $\yb^i$ corresponding to the $i$th layer equation in \cref{eq:DNN}, for $i = 0, \hdots, n_\ell$. Thus, the $i$th adjoint equation corresponds to the derivative with respect to $\ub^i = \ab^i$ in \cref{eq:adjointFinite}, and it reads
\begin{equation}
\begin{aligned}
    \yb^{n_\ell} &= \ab^{obs} - \ab^{n_\ell}, \\
    \yb^{i} &=\LRp{\Wm^{i+1}}^T \LRs{\sigma'\LRp{\Wm^{i+1} \ab^{i} + \bb^{i+1} } \circ \yb^{i+1}}, \quad i = n_\ell - 1,\hdots,0
\end{aligned}
\label{eq:adjointDNN}
\end{equation}
where $\sigma'$ is the derivative of $\sigma$, and $\circ$ denotes the component-wise multiplication of two vectors. Note that since $\sigma$ acts componentwise (when its input is a vector), so is its derivative $\sigma'$: in particular $\sigma'\LRp{\Wm^{i+1} \ab^{i} + \bb^{i+1} }$ is a vector in \cref{eq:adjointDNN}.
Thus, \cref{eq:adjointDNN} provides explicit expressions for the adjoint equations for a general fully connected DNN. Again, note that the adjoint equations are linear in terms of adjoint variables $\yb^i$, $i = 0,\hdots,n_\ell$.

\item To unfold the control equation \cref{eq:controlFinite} to explicitly see the derivative of the objective function with respect to the weights and biases, we take a block of control equations corresponding to sub-blocks of $\rrb^i = \LRs{\rrb_\Wm^i, \rrb_{\bb}^i}^T$ in $\rrb$. For DNN, these derivatives are given as: for $i = 1,\hdots,n_\ell$,
\begin{equation}
\begin{aligned}
    \pp{f}{\Wm^i} &=- \LRs{\vec{\yb^i}\circ \sigma'\LRp{\Wm^{i} \ab^{i-1} + \bb^{i} } }\LRp{{\ab}^{i-1}}^T, \\
    \pp{\f}{{\bb}^i} &=- {\yb^i} \circ \sigma'\LRp{\Wm^{i} \ab^{i-1} + \bb^{i} }.
\end{aligned}
\label{eq:gradientDNN}
\end{equation}

The {\bf backpropagation} nature of the network gradient is now clearly seen from the gradient expressions in \cref{eq:gradientDNN} and the adjoint equations \cref{eq:adjointDNN}. Indeed, in \cref{eq:gradientDNN} we need the $i$th adjoint state $\yb^i$ in order to compute the gradients with respect to the weights and biases in the $i$th layer. The  $i$th adjoint state $\yb^i$ in turn depends on the $(i+1)$th adjoint state $\yb^{i+1}$, which depends on the $(i+2)$th adjoint state $\yb^{i+2}$, etc, all the way to the last adjoint state $\yb^{n_\ell}$ corresponding to the network output layer. In other words, using the adjoint equations \cref{eq:adjointDNN} we backpropagate to compute the adjoint solution from the output layer to the $i$th layer, and then compute the gradients using \cref{eq:gradientDNN}. From the backpropagation point of view, $\yb^i$, $i = 1,\hdots,n_\ell$ are simply the temporary variables to help compute/write the chain rule in a succinct manner. {\em The adjoint approach, however, reveals their precise role as the adjoint solutions\textemdash also known as the Lagrangian multipliers\textemdash of the adjoint equations steming from the first order optimality condition using the reduced space approach in \cref{rema:adjoint}.}
Another important view point that we have exploited here  is that the DNN training problem, from the adjoint point of view, is a constrained optimization problem with the forward pass as the forward equations. {\em The backpropagation is thus nothing more than a reduced space approach to compute the gradient using adjoint method.}

\end{itemize}

\subsection{Application of adjoint to the stability of ordinary differential equations}
\label{sect:ODEs}
In this section, we provide a brief view on the role of adjoint in the study of stability of the equilibria of ordinary differential equations (ODEs). Most of  our mathematical exposition follows \cite{Logemann2014}, and we limit ourselves to autonomous systems of the form
\begin{equation}
\dot{\xb} := \dd{\xb}{t} = \fb\LRp{\xb},
\label{eq:ODEsystem}
\end{equation}
where $\xb \in \real^n$ and $\fb: \Gs \subset \real^n \to \real^n$ is assumed to be continuous and locally Lipschitz. The domain $\Gs$ of $\fb$ is assumed to be a nonempty open subset of $\real^n$. Due to the translational invariance of autonomous system, without loss of generality, we can assume that $\bs{0} \in \Gs$ is an equilibrium point, i.e. $\fb\LRp{\bs{0}} = \bs{0}$. We use $\xb_0$ to denote the initial condition and  $I$ to denote the maximal interval of existence for a solution of \cref{eq:ODEsystem}. All the norms $\nor{\cdot}$ and inner products $\LRp{\cdot,\cdot}$ in this section are the standard Euclidean ones, and in this case \cref{exa:matrixAdjoint} shows that the adjoint of a real matrix is simply its transpose. For matrices, the norm is the induced operator norm. 

\begin{definition}[Lyapunov stability]
The equilibrium point $\bs{0}$ is stable (in the sense of Lyapunov) if for any $\epsilon > 0$, $\exists \delta > 0$ such that for every (maximal) solution $\xb: I \to \Gs$ such that $\xb\LRp{0} \le \delta $, we have $\xb\LRp{t} \le \epsilon$ for all $t \in I \cap \LRp{0,\infty}$. 
\end{definition}

\begin{theorem}[Lyapunov direct method]
    If there exists an open neighborhood $\Uset$ of $\bs{0}$ and a continuous differentiable function $V$ such that
    \begin{enumerate}
        \item $V\LRp{\bs{0}} = 0$ and $V\LRp{\bs{z}} > 0$ for all $\zb \in \Uset \setminus \LRc{\bs{0}}$, and
        \item $V_{\fb}\LRp{\zb} := \LRp{\Grad V\LRp{\zb}, \fb\LRp{\zb}} := \sum_{i=1}^n \pp{V}{\zb_i}\fb_i\LRp{\zb} \le 0$ for all $\zb \in \Uset$.
    \end{enumerate}
    Then $\bs{0}$ is a stable equilibrium  point of \cref{eq:ODEsystem}.
\end{theorem}
\begin{proof}
    See \cite[Theorem 5.2]{Logemann2014}.
\end{proof}

\begin{definition}[Asymptotic stability]
The equilibrium $\bs{0}$ is attractive if there exists $\delta > 0$ such that for every $\x_0 \in \Gs$ such that $\nor{\xb_0} \le \delta$, then the solution $\xb\LRp{t} \to \bs{0}$  as $t \to \infty$. We say $\bs{0}$ asymptotically stable (in the sense of Lyapunov) if it is both stable and attractive.
\end{definition}

\begin{theorem}[A sufficient condition for asymptotic stability]
Assume that there exists a neighborhood $\Uset$ of $\xb_{0}$ and a continuously differentiable function $V$ such that 
\begin{itemize}
    \item[i] $V\LRp{\bs{0}} = 0$ and $V\LRp{\bs{z}} > 0$ for all $\zb \in \Uset \setminus \LRc{\bs{0}}$, and $V_{\fb}\LRp{\zb}  \le 0$ for all $\zb \in \Uset$, and
    \item[ii] $\bs{0}$ is the inverse image of $V_{\fb}\LRp{\zb} = 0$, i.e., $V^{-1}_{\fb}\LRp{0} = \bs{0}$.
\end{itemize}
    Then $\bs{0}$ is asymptotically stable.
    \label{theo:AsymptoticStability}
\end{theorem}
\begin{proof}
    See \cite[Theorem 5.15]{Logemann2014}.
\end{proof}

{\em We next study the stability of systems of linear ODEs, and this is where the {\em adjoint} comes into the picture}. We then
infer the stability of nonlinear systems using the stability of their linearizations. To that end, we consider linear systems with $\Gs = \real^n$, $\Am \in \real^{n\times n}$, and
\begin{equation}
    \dot{\xb} = \fb\LRp{\xb} = \Am\xb.
    \label{eq:ODElinear}
\end{equation}
Clearly, the solution of the \cref{eq:ODElinear} can be written as matrix exponential \cite{Hall2015,Moler2003,HornJohnson91}
\[
\xb\LRp{t} = \exp\LRp{\Am t}\xb_0.
\]

\begin{definition}[Exponential stability]
    The equilibrium $\bs{0}$ is called exponentially stable if there exist $M \ge 1$ and $\alpha > 0$ such as
    \[
    \nor{\exp\LRp{\Am t}\xb_0} \le M \exp\LRp{-\alpha t} \nor{\xb_0}, \quad \forall t \ge 0 \text{ and } \forall \xb_0 \in \real^n.
    \]
\end{definition}

\begin{definition}[Hurwitz matrices]
Let $\sigma\LRp{\Am}$ denote the spectrum (the collection of all eigenvalues) of $\Am$. $\Am$ is Hurwitz if $\sigma\LRp{\Am} \subset \LRc{\lambda \in \mbb{C}: \Re\LRp{\lambda} < 0}$.
\end{definition}

\begin{proposition}
Let $\Am \in \real^{n \times n}$. The following statements are equivalent:
\begin{itemize}
    \item[i)] $\Am$ is Hurwitz.
    \item[ii)] $\bs{0}$ is an exponentially stable equilibrium of \cref{eq:ODElinear}.
    \item[iii)] $\bs{0}$ is an asymptotically stable equilibrium of \cref{eq:ODElinear}.
\end{itemize}
\label{propo:HurwitzExponentialStability}
\end{proposition}
\begin{proof}
    See \cite[Proposition 5.25]{Logemann2014}.
\end{proof}

We are in the position to discuss one of the main results of this section.
\begin{theorem}[Necessary and sufficient conditions for exponential stability]
    $\Am \in \real^{n \times n}$ is Hurwitz iff for each symmetric positive definite (SPD) matrix $\mc{Q} \in \real^{n \times n}$, the matrix equation
    \[
    \mc{P} \Am + \Am^*\mc{P} + \mc{Q} = 0
    \]
    has a SPD solution $\mc{P} \in \real^{n\times n}$.
    \label{theo:Hurwitz}
\end{theorem}
\begin{proof}
    For the necessity, suppose $\Am$ is Hurwitz. It follows from \cref{propo:HurwitzExponentialStability} that $\bs{0}$ is an exponentially stable equilibrium, that is, there exists $M > 0$ and $\alpha>0$ such that
    \[
    \nor{\exp\LRp{\Am t}} = \sup_{\xb_0 \in \real^n} \frac{\nor{\exp\LRp{\Am t}\xb_0}}{\nor{\xb_0}} \le M \exp\LRp{-\alpha t}, \quad t \ge 0.
    \]
    Now for any SPD matrix $\mc{Q}$, let us define
    \[
    \mc{P} := \int_0^\infty \exp\LRp{\Am^*}\mc{Q} \exp\LRp{\Am t} \,dt,
    \]
    which is a well-defined matrix since
    \[
    \nor{\mc{P}} \le \int_0^\infty \nor{\exp\LRp{\Am^* t}} \nor{\mc{Q}} \nor{\exp\LRp{\Am t}}\,dt \le M \nor{\mc{Q}}\int_0^\infty \exp\LRp{-2\alpha t} \,dt < \infty,
    \]
    where in the second inequality we have used the fact from \cref{propo:adjointExistence} that the norm of a linear continuous operator is equal to the norm of its adjoint. $\mc{P}$ is SPD as
    \[
    \LRp{\xb, \mc{P}\xb} = \int_0^\infty \LRp{\xb,\exp\LRp{\Am^*t}\mc{Q} \exp\LRp{\Am t}\xb} \,dt = \int_0^\infty \LRp{\exp\LRp{\Am t}\xb,\mc{Q}\exp\LRp{\Am t}\xb} \,dt,
    \]
    and the fact that $\mc{Q}$ is SPD. Furthermore,
    \[
    \mc{P}\Am + \Am^*\mc{P} = \int_0^\infty \dd{\LRp{\exp\LRp{\Am^*t}\mc{Q} \exp\LRp{\Am t}}}{t}\,dt = - \mc{Q},
    \]
    where we have used the fact that $\exp\LRp{\Am t}$, and hence $\exp\LRp{\Am^* t}$, decays exponentially to $0$.

For the sufficiency, \cref{propo:HurwitzExponentialStability} says that we only need to show that $\bs{0}$ is an asymptotically stable equilibrium. To that end, let us construct the following function
\[
V\LRp{\zb} := \LRp{\zb, \mc{P}\zb}. 
\]
Clearly $V\LRp{\zb} \ge 0$ for all $\zb \in \real^n$, and $V\LRp{\zb} = 0$ iff $\zb = \bs{0}$ as $\mc{P}$ is SPD. Furthermore
\[
V_{\fb}\LRp{\zb} = 2\LRp{\mc{P}\zb,\Am\zb} = \LRp{\zb, \LRp{\mc{P}\Am + \Am^*\mc{P}}\zb}= -\LRp{\zb,\mc{Q}\zb} \le 0.
\]
Thus, by \cref{theo:AsymptoticStability}, $\bs{0}$ is asymptotically stable, and this ends the proof. 
\end{proof}

Let us now use \cref{theo:Hurwitz} to study the stability of the equilibrium $\bs{0}$ of the general nonlinear system \cref{eq:ODEsystem}. 
\begin{hypothesis}[Nonlinearly purturbed linear ODE systems]
We assume that $\fb\LRp{\xb} = \Am\xb + \hb\LRp{\xb}$ where $\Am \in \real^{n\times n}$ and $\hb: \Gs \to \real^n$ is continuous with
\begin{equation}
\lim_{\zb \to \bs{0}} \frac{\nor{\hb\LRp{\zb}}}{\nor{\zb}} = 0,
\label{eq:linearizedAssumption}
\end{equation}
that is, $\hb\LRp{\zb} = o\LRp{\zb}$. In other words, $\hb$ approaches $\bs{0}$ faster than $\zb$. 
\label{hypo:NonlinearPerburbation}
\end{hypothesis}
\begin{theorem}[Linear stability implies nonlinear stability]
    Assume  \cref{hypo:NonlinearPerburbation} holds. If $\bs{0}$ is an asymptotically stable equilibrium of \cref{eq:ODElinear}, it is also an asymptotically stable equilibrium of \cref{eq:ODEsystem}.
    \label{theo:linearStabilityImpliesNonlinearStability}
\end{theorem}
\begin{proof}
Suppose $\bs{0}$ is an asymptotically stable equilibrium of \cref{eq:ODElinear}. Using \cref{theo:Hurwitz} we can pick $\mc{Q} = \Im$, and form $V\LRp{\zb} := \LRp{\zb,\mc{P}\zb}$. It follows that $V\LRp{\bs{0}} = 0$ and $\V\LRp{\zb} > 0$ for all $\zb \ne \bs{0}$ owing to the SPD property of $\mc{P}$. In order to use \cref{theo:AsymptoticStability} to conclude the proof, we just need to show that $V_{\fb}\LRp{\zb} < 0$. We have
\[
V_{\fb}\LRp{\zb} = 2\LRp{\mc{P}\zb,\Am\zb + \hb\LRp{\zb}}= 
-\nor{\zb}^2 + 2\LRp{\mc{P}\zb,\hb\LRp{\zb}} \le -\nor{\zb}^2 + 2\nor{\mc{P}}\nor{\zb}\nor{\hb\LRp{\zb}},
\]
where we have used $\mc{P}\Am + \Am^*\mc{P} = -\Im$ in the second equality. Next, using \cref{eq:linearizedAssumption} we can pick a sufficient small neighborhood $\Uset = \Ball{\epsilon}{\bs{0}}$ such that $\nor{\hb\LRp{\zb}} \le \frac{\nor{\zb}}{4\nor{\mc{P}}}$. Thus $V_{\fb}\LRp{\zb} < 0$ for all $\zb \in \Uset$, and this ends the proof.
\end{proof}

\cref{theo:linearStabilityImpliesNonlinearStability} reduces the stability analysis of an equilibrium of a nonlinear ODE system to an appropriate linear ODE counterpart, which is much simpler as linear algebra is then all we need for studying the stability. {\em We have also seen in \cref{theo:Hurwitz} and \cref{theo:linearStabilityImpliesNonlinearStability} that   adjoint plays the key role in establishing a necessary and sufficient condition for the stability of an equilibrium of an ODE system.}

\begin{example}[Linearization of nonlinear ODE systems]
    We assume $\fb: \Gs \to \real^n$ is differentiable at $\bs{0}$, and $\bs{0} \in \Gs$ is an equilibrium. Let us set
    \[
    \Am := \Grad \fb\LRp{\bs{0}},
    \]
    where the gradient is defined in \cref{def:gradient} and \cref{ex:gradFiniteVector}. The definition of Fr\'echet derivative in \cref{eq:TaylorFrechet} implies \cref{eq:linearizedAssumption} with $\hb\LRp{\xb} := \fb\LRp{\xb} - \Am\xb$. Then by \cref{theo:linearStabilityImpliesNonlinearStability}, the asymptotic stability of $\bs{0}$ for the linearized system implies the asymptotic stability of $\bs{0}$ for the original nonlinear system.
\end{example}

\begin{example}[The asymptotic stability meaning of the basic reproduction in epidemic modeling]
With the ever-increasing human population on every part of the earth, the shrinkage in the natural habitat for plants and animals, and the shortage of natural resources such as water and food, the emergence of new and re-emergence of old infectious diseases are inevitable. Epidemic modeling plays a key role in forecasting how an infectious disease (such as SARS and COVID-19) spreads. This in turn facilitates informative decision-making to prevent a disease outbreak. The reproduction number $\mathfrak{R}_0$ provides epidemiologically meaningful criteria to predict an outbreak \cite{Perasso2018,Heffernan2005,Brauer2019,Martcheva2015,Dobson2020-um,Ledder2023,Daley1984,BrauerEtal2008,Hethcote2000,vandenDriessche2002}. Its popular definition is ``the number of secondary cases one infected individual produces in a population consisting of only susceptibles". As a result, if $\mathfrak{R}_0 < 1$ the disease dies out but persists as an endemic (or goes on extinction) if $\mathfrak{R}_0 > 1$. 

One            
    of the most popular approaches to model disease dynamics is to use ODEs \cite{Brauer2019,Martcheva2015,Dobson2020-um,Ledder2023,Daley1984,BrauerEtal2008,Hethcote2000} in which, for example in an Susceptible-Exposed-Infectious-Recovered-Susceptible (SEIRS) model, the components of $\xb$ in \cref{eq:ODEsystem} are typically the fraction of susceptible, exposed, infected, and recovered within the  population under consideration. The dying out of a disease corresponds to an outbreak returning to a disease-free state, while persistence to an endemic corresponds to a disease that remains in the population. Clearly, the disease-free state is an equilibrium of the ODE system if it is a meaningful representation of the disease dynamics. The dying out of a disease, therefore, corresponds to a perturbation from and then a return to disease-free equilibrium (DFE). This in turn should correspond to the asymptotic stability of the DFE. This is exactly a mathematical justification of the reproduction number. \cref{fig:reproductionNumberMathematics} is our effort\footnote{Unpublished notes.} in sketching the association, via the next generation matrix approach \cite{Hurford2009,vandenDriessche2002,Diekmann2009},  of the reproduction number being less than unity and the asymptotic stability  of $\bs{0}$ as a DFE  of an abstract epidemic ODE model $\dot{\xb} = \fb\LRp{\xb}$. Note that the top arrow is the implication (due to \cref{theo:linearStabilityImpliesNonlinearStability}) and the rest are equivalences. (See our work \cite{akuno2022multipatch}, and the references therein, for a detailed exposition of this correspondent for a new SEIRS epidemic model.) As can be seen, when the number of secondary cases one infected individual produces is less than one, the DFE is asymptotically stable and the solution of the epidemic ODE model approaches the DFE as time goes on. Epidemically speaking, the disease dies out.

\tikzstyle{box} = [rectangle, minimum width=3cm, minimum height=0.6cm,text centered, draw=black]
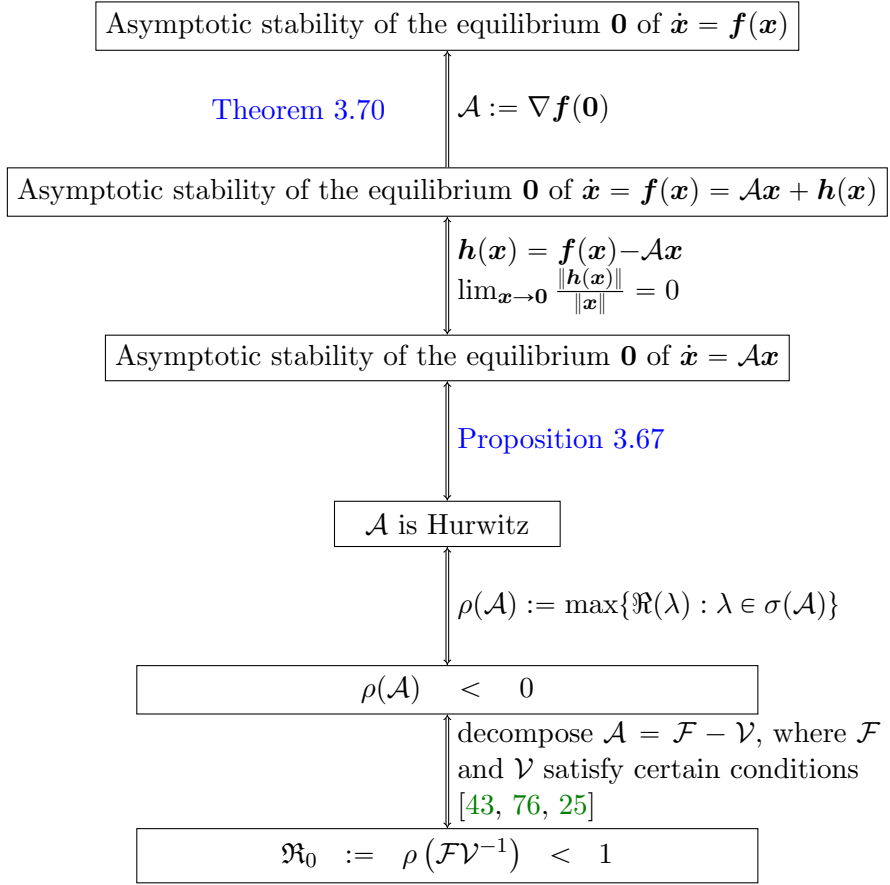
\begin{figure}[h!t!b!]
    \centering
    \begin{tikzpicture}[node distance=2.2cm]

  \node (I) [box] {Asymptotic stability of the equilibrium $\bs{0}$ of  $\dot{\xb}=\fb(\xb)$};
  \node (II) [box, below of=I] {Asymptotic stability of the equilibrium $\bs{0}$ of $\dot{\xb} = \fb(\xb) = \Am\xb + \hb(\xb)$};
  \node (III) [box, below of=II] { Asymptotic stability of the equilibrium $\bs{0}$ of $\dot{\xb} = \Am\xb$};
  \node (IIIa) [box, below of=III] { $\Am \text{ is Hurwitz}$};
  \node (IV) [box, text width=8cm, below of=IIIa] {
    $\rho(\Am) < 0 $
  };
  \node (V) [box, text width=8cm, below of=IV] {
    $ \mathfrak{R}_0 := \rho\LRp{\Fm\Vm^{-1}} < 1$
  };

  \draw [{implies}-{}, double] (I) -- 
  node[anchor=east, text width= 3cm] {\cref{theo:linearStabilityImpliesNonlinearStability}}
  node [anchor=west, text width=3cm] {$\Am :=\Grad\fb(\bs{0})$}(II);
  
  \draw [{implies}-{implies}, double] (II) -- node[anchor=west, text width=3cm] {$\hb(\xb)=\fb(\xb)-\Am\xb$\\$\lim_{\xb\to \bs{0}} \frac{\nor{\hb(\xb)}}{\nor{\xb}} = 0$} (III);

  \draw [{implies}-{implies}, double] (III) -- node[anchor=west] {\cref{propo:HurwitzExponentialStability}} (IIIa);
  
  \draw [{implies}-{implies}, double] (IIIa) -- node[anchor=west] {$\rho(\Am) := \max\{\Re(\lambda): \lambda \in \sigma(\Am)\}$} (IV);
  
  \draw [{implies}-{implies}, double] (IV) -- node[anchor=west, text width=6cm] {decompose $\Am=\Fm-\Vm$, where $\Fm$ and $\Vm$ satisfy certain conditions \cite{Hurford2009,vandenDriessche2002,Diekmann2009}} (V);
\end{tikzpicture}
    \caption{A sketch of the association, via the next generation matrix approach,  of the reproduction number being less than unity and the asymptotic stability  of the disease-free equilibrium $\bs{0}$ of an abstract epidemic ODE model $\dot{\xb} = \fb\LRp{\xb}$. Note that the top arrow is the implication (due to \cref{theo:linearStabilityImpliesNonlinearStability}) and the rest are equivalences.}
    \label{fig:reproductionNumberMathematics}
\end{figure}
    
\end{example}


\section{Part II: Adjoint operators in Infinite dimensional Hilbert spaces}
\label{sect:PartII}
 In this second part,
 we shall  consider  linear mappings between infinite dimensionsal Hilbert spaces. Of our particular interest are densely defined differential operators. 
 A sufficient general definition of adjoint operator that could embrace a variety of problems is the following (see, e.g., \cite{Showalter98}).
 


\begin{definition}[Adjoint operators]
$\A^*$ is called the adjoint of $\A \in \L\LRp{\Us,\Vs}$ iff
\begin{equation}
    \label{eq:adjointInfiniteDef}
    \LRp{\A\u,\v}_\Vs = \LRp{\u,\A^*\v}_\Us \quad \forall \u \in \Do\LRp{\A} \text{ and } \forall \v \in \Do\LRp{\A^*},
\end{equation}
where 
\[
\Do\LRp{\A^*}:= \LRc{\v: \text{ the map } \u \mapsto \LRp{\A\u,\v}_\Vs \text{ is continuous on } \Us}.
\]
\label{defi:adjointInfinite}
\end{definition}

Clearly, \cref{defi:adjointInfinite} reduces to \cref{defi:adjoint} when $\A$ is continuous and is defined on the whole space $\Us$.

\begin{example}
We consider an elliptic differential operator in $n$ dimensions over an open and bounded domain $\Omega \subset \real^n$. In this case, we define
$\Ls^2-$inner product of two functions $\u(\xb),\v(\xb)$ in $\Ls^2(\Omega)$ over $\R$ as
\begin{equation}
    (\u,\v)_{L^2(\Omega)} = \int_\Omega \u(\xb) \v(\xb) d\Omega.
    \label{eq:L2inner}
\end{equation}
Consider the following parametrized linear operator $\A: \Do\LRp{\A} \subset \Ls^2\LRp{\Omega} \to \Ls^2\LRp{\Omega}$ as
\[ 
\A \u = \begin{cases}
-\Div\LRp{e^{\z}\Grad \u} & \text{ in } \Omega,\\ 
\u = 0 & \text{ on } \pOmega,
\end{cases}
\]
where 
$\rr \in \Cs^1\LRp{{\Omega}} \subset \Ls^2\LRp{\Omega}$
is a distributed parameter over $\Omega$ with sufficiently smooth boundary $\pOmega$. 
We define the graph space $\Hs_\A$ as
\[
\Hs_\A := \LRc{\u \in \Hs^1\LRp{\Omega}: \A\u \in \Ls^2\LRp{\Omega}} =  \LRc{\u \in \Hs^1\LRp{\Omega}: -\Div\LRp{e^{\z}\Grad \u} \in \Ls^2\LRp{\Omega}}.
\]
The domain of $\A$ is chosen to be a subset of its graph space:
\[
\Do\LRp{\A} := \LRc{\w \in \Hs_\A: \u = 0 \text{ on } \pOmega}
\]
equipped with the graph norm $\nor{\u}_{\mc{G}} := \sqrt{\nor{\u}^2_{\Hs^1} + \nor{\A\u}^2_{\Ls^2}} = \sqrt{\nor{\u}^2_{\Hs^1} + \nor{\Div\LRp{e^{\z}\Grad \u}}^2_{\Ls^2}}$. Clearly, $\Do\LRp{\A}$ is dense in $\Ls^2\LRp{\Omega}$ and  $\A$ is continuous on $\Do\LRp{\A}$.
The proof of the self-adjointness of $\A$ is straightforward using basic facts from weak/distributional  derivative and a standard distributional arugment and it is provided in \cref{proof:closedEllipticOp} of \cref{sect:appendix}. 
\label{exa:ellipticOp}
\end{example}

\begin{example}
    We now consider a weak formulation for the elliptic differential operator in \cref{exa:ellipticOp}. This  will be important for studying the well-posedness of the associated partial differential equations in \cref{exa:ellipticPDEs}. Multiplying $-\Div\LRp{e^{\z}\Grad \u}$ by a test function and integrating by parts allow us to define the following bilinear form $a\LRp{\cdot,\cdot}: \Hs^1_0\LRp{\Omega}\times\Hs^1_0\LRp{\Omega} \to \real$
\[
a\LRp{\u,\v} := \LRp{e^\z\Grad\u,\Grad\v}_{\Ls^2(\Omega)}.
\]
By the Cauchy-Schwarz inequality we have 
\[
\snor{a\LRp{\u,\v}} \le \nor{e^\z}_{\Ls^\infty\LRp{\Omega}} \nor{\Grad\u}_{\Ls^2(\Omega)} \nor{\Grad\v}_{\Ls^2(\Omega)} \le \nor{e^\z}_{\Ls^\infty\LRp{\Omega}} \nor{\u}_{\Hs^1(\Omega)} \nor{\v}_{\Hs^1(\Omega)}, 
\]
and thus $a\LRp{\cdot,\cdot}$ is continuous on $\Hs^1_0\LRp{\Omega}\times\Hs^1_0\LRp{\Omega}$. Here, $\Ls^\infty\LRp{\Omega}$ is the space of essentially bounded functions on $\Omega$.
As a result, it implicitly defines a unique linear and continuous operator\footnote{Note that for every continuous sequilinear form $a: \Xs \times \Ys \to \Fs$, there exists a unique linear and continuous operator $\A: \Xs \to \Ys$ such that
$\LRp{\A\u,\v}_\Ys := a\LRp{\u,\v}$ for all $\u \in \Xs$ and $\v \in \Ys$. Indeed, by the linearity and continuity of $a\LRp{\cdot,\cdot}$ with respect to its first argument, the Riesz representation \cref{theo:Riesz} ensures that there exists a unique $\A{\u} \in \Ys$ such that
$\LRp{\A\u,\v}_\Ys = a\LRp{\u,\v}$. The continuity of $\A$ is from the continuity of $a\LRp{\u,\v}$: $\nor{\A\u} = \sup_{\v \in \Ys} \frac{\LRp{\A\u,\v}}{\nor{\v}_\Ys} = \sup_{\v \in \Ys} \frac{a\LRp{\u,\v}}{\nor{\v}_\Ys} \le \beta \nor{\u}_\Xs$. The uniqueness of $\A$ is straightforward due to its definition, as if there were another linear and continuous operator $\B$, then we would have
\[
\nor{\B - \A} = \sup_{\u \in \Xs} \frac{\nor{\B\u-\A\u}_\Ys}{\nor{\u}_\Xs} = \sup_{\u \in \Xs}\sup_{\v \in \Ys} \frac{\LRp{\B\u-\A\u,\v}_\Ys}{\nor{\u}_\Xs\nor{\v}_\Ys}  = \sup_{\u \in \Xs}\sup_{\v \in \Ys} \frac{a\LRp{u,\v} - a\LRp{u,\v}}{\nor{\u}_\Xs\nor{\v}_\Ys} =  0,
\]
and thus $\B = \A$.
\label{fn:bilinear}
} 
$\A: \Hs^1_0(\Omega) \to \Hs^1_0(\Omega)$ as
\[
\LRp{\A\u,\v}_{\Hs^1_0(\Omega)} := a\LRp{\u,\v}.
\]
Due to the symmetry of the bilinear form, we have $\A^* = \A$.
    \label{exa:ellipticOpWeak}
\end{example}

\begin{example}
    Consider an open and bounded domain $\Omega \subset \real^n$ and the $\Ls^2-$inner product of two functions $\u(\xb),\v(\xb)$ is given in \cref{eq:L2inner}. Consider the following parametrized linear operator 
\[ 
    \A \u = \begin{cases}
    \bs{\beta} \cdot \Grad \u + \lambda \u & \text{ in } \Omega, \\
\bs{\beta}\cdot\nb\u = 0 & \text{ in } \pOmega_{\text{in}},
    \end{cases}
\]
where $\bs{\beta} \in \LRs{\Cs^1\LRp{{\Omega,\R}}}^n$ and $\Div\bs{\beta} = 0$, $\lambda > 0$, $\nb$ is the unit outward normal vector of the boundary $\pOmega$, and $\pOmega_{\text{in}} := \LRc{\xb \in \pOmega: \bs{\beta}\cdot\nb <0}$ is the inflow boundary. We consider  the graph space $  \Hs_\A := \Hs^1_{\bs{\beta}}\LRp{\Omega}\equaldef\LRc{\u: \u \in
  \mbb{L}^2\LRp{\Omega} \text { and } \bs{\beta}\cdot \Grad \u \in
  \mbb{L}^2\LRp{\Omega}}$ which is dense in $\Ls^2\LRp{\Omega}$ under sufficient regularity\footnote{Assume $\Omega$ has segment property \cite{AntonicBurazin09}.} of the domain $\Omega$ \cite{ErnGuermond06}.
  The domain of $\A$ is defined as a subset of the graph space, namely, $\Do\LRp{\A} := \LRc{\u \in \Hs^1_{\bs{\beta}}\LRp{\Omega}: \bs{\beta}\cdot\nb\u = 0  \text{ in } \pOmega_{\text{in}}}$.
  It is clear that $\A: \Do\LRp{\A} \to \Ls^2\LRp{\Omega}$ is linear and continuous owing to the definition of $\Hs^1_{\bs{\beta}}\LRp{\Omega}$ and the intrinsic graph norm $\nor{\u}_{\Hs^1_{\bs{\beta}}\LRp{\Omega}} := \sqrt{\nor{\u}_{\Ls^2}^2 + \nor{\bs{\beta}\cdot \Grad \u}^2}$. Using the definition of weak/distributional derivative and integration by parts once, we can show (see \cref{exa:Friedrichs} for a more general differential operator) that
  \[ 
    \A^* \v = \begin{cases}
    -\bs{\beta} \cdot \Grad \v + \lambda \v & \text{ in } \Omega, \\
\bs{\beta}\cdot\nb\v = 0 & \text{ in } \pOmega_{\text{out}},
    \end{cases}
\]
where $\pOmega_{\text{out}} := \LRc{\xb \in \pOmega: \bs{\beta}\cdot\nb > 0}$ is the outflow boundary and 
\[ 
\Do\LRp{\A^*} = \LRc{\v \in \Hs^1_{\bs{\beta}}\LRp{\Omega}: \bs{\beta}\cdot\nb\u = 0  \text{ in } \pOmega_{\text{out}}}.
\]
\label{exa:advectionOp}
\end{example}

\begin{example}[Friedrichs' systems]
 We are interested in
Friedrichs' system that  embraces a large class of elliptic,
parabolic, hyperbolic, and mixed-type PDEs operators \cite{Friedrichs58}: 
\begin{equation}
\eqnlab{LPDE}
\A\u := \sum_{k=1}^{\n}\Am_k \partial_k \u + \Cm \u 
 \text{ in } \in \Omega,
\end{equation}
where $\d$ is the spatial dimension, $\u$
the unknown solution with values in $\real^m$, $\f$ the forcing term, and $\Omega $ is an open and bounded subset of $ \R^n$ with sufficient regular boundary $\pOmega$. The matrices $\Am_k$
and $\Cm$ are assumed to be constant and continuous across $\Omega$. Here,
$\partial_k$ is understood as the $k$th partial derivative. 
We start with the standard assumptions
(see, e.g., \cite{Friedrichs58, ErnGuermond06, Jensen04, ErnGuermondCaplain07}):
\begin{subequations}
\label{eq:Friedrichs}
\begin{align}
&\Cm \in \LRs{\Ls^\infty\LRp{\Omega}}^{\m,\m}, \label{eq:FriedrichsCcond}\\ 
&\Am_k \in
  \LRs{\Ls^\infty\LRp{\Omega}}^{\m,\m}, \quad k = 1,\hdots,n, \quad \text{and }
  \sum_{k=1}^{\n}\partial_k\Am_k \in \LRs{L^\infty\LRp{\Omega}}^{\m,\m},
\label{eq:FriedrichsAcond}  \\
& \Am_k = \LRp{\Am_k}^T \text{ in } \Omega, \quad k= 1,\hdots,n, \eqnlab{Asym}, \\
\label{eq:FriedrichsCoercivity}
& \Cm+ \Cm^T + \sum_{k=1}^{\n}\partial_k\Am_k \ge 2\alpha_0 \Im \text{ in } \Omega,
\end{align}
where $\alpha_0 > 0$ is some coercivity constant.
In this paper, we consider the following abstract boundary condition for $\A$ \cite{ErnGuermond06}: let $\Dm := \sum_{k=1}^\n{\Am_k\nb_k}$, where $\nb_k$ is the $k$th component of the unit outward normal vector $\nb$ on $\pOmega$, and assume there exists $\Mm \in \LRs{\Ls^\infty\LRp{\pOmega}}^{\m,\m}$ such that
\end{subequations}
\begin{subequations}
\label{eq:M}
\begin{align}
\label{eq:FriedrichsBCM}
\Mm + \Mm^T &\ge 0  \text{  on } \pOmega, \\
\label{eq:FriedrichsBCu}
\LRp{\Dm - \Mm}\u &= 0  \text{  on } \pOmega, \\
\label{eq:FriedrichsBCMD}
\N\LRp{\Dm - \Mm} + \N\LRp{\Dm + \Mm} &= \R^\m, \text{  on } \pOmega.
\end{align}
\end{subequations}
Following \cite{ErnGuermond04}, we define\footnote{Note that we pick the $\Ls^2$-setting here for concreteness, but all the results for Friedrichs' system/operator in this paper hold for the general Hilbert space setting in \cite{ErnGuermond04, ErnGuermond06}.} 
the graph space of $\A$ using its differential part $\B := \A - \Cm$:
\[
\Hs_\A := \LRc{\u \in \LRs{\Ls^2\LRp{\Omega}}^\m: \B\u \in \LRs{\Ls^2\LRp{\Omega}}^\m},
\]
which is dense \cite{ErnGuermondCaplain07} in $\LRs{\Ls^2\LRp{\Omega}}^\m$ when $\Omega$ is sufficiently regular (see \cref{exa:advectionOp} for an example). Furthermore, 
from the definition, it is easy to see that $\A$ is linear and continuous on $\Do\LRp{\A} \equaldef 
\LRc{\x \in \Hs_\A: \LRp{\Dm - \Mm}\x = 0 \text{ on } \pOmega}$ equipped with intrinsic graph norm.
Using the definition of weak/distributional derivative and integration by parts  we can show  that its adjoint is found to be (see \cref{proof:adjointFriedrichsOp} in \cref{sect:appendix}):
\[
\A^*\v = -\sum_{k=1}^{\n}\partial_k \LRp{\Am_k\v} + \Cm^T \v \in \Ls^2\LRp{\Omega}, 
\]
for any $\v \in \Do\LRp{\A^*}$, where
\[
\Do\LRp{\A^*} = \LRc{\v \in \Hs_\A: \Dm\v \in \LRs{\N\LRp{\Dm - \Mm}}^\perp}.
\]
\label{exa:Friedrichs}
\end{example}

\begin{remark}
    Note that we have considered Friedrichs systems with full coercivity in \cref{exa:Friedrichs}. Applying such a general setting to various concrete PDEs \cite{ErnGuermond04, ErnGuermond06} will reveal concrete adjoints and their domains, but we omit the details here. The elliptic operator in \cref{exa:ellipticOp}, when written in the first order form, is a particular two-field Friedrichs system \cite{ErnGuermond04, ErnGuermond06}, whose adjoint can be derived similarly.
\end{remark}

\subsection{Application of adjoint to ill-posed problems}
\label{sect:illposedApplication}
In this section, we consider continuous linear operator defined on the whole space and we shall extend the spectral decomposition in \cref{coro:spectralFinite} and SVD decomposition in \cref{theo:SVDfinite} to compact (linear) operators in infinite dimensions. This allows us to show that inverting a compact operator is an ill-posed problem. We then explain rigorously how  the standard Tikhonov regularization could overcome the ill-posedness. We begin by recalling the definition of compact linear operators and some of its consequences.

\begin{definition}[Compact operator]
Let $\A: \Xs \to \Ys$ be linear. We say that $\A$ is a compact operator if for every bounded sequence $\LRc{\u_i}_{i=1}^\infty \subset \Xs$, the sequence $\LRc{\A\u_i}_{i=1}^\infty \subset \Ys$ has a convergent subsequence.
\label{defi:compact}
\end{definition}
A direct consequence of \cref{defi:compact} is that any compact operator is a linear and continuous map.
\begin{corollary}
If $\A: \Xs \to \Ys$ is compact, then $\A \in \Bf\LRp{\Xs,\Ys}$.
\label{coro:compactContinuous}
\end{corollary}
Self-adjoint compact operators in Hilbert spaces possesses many desirable properties among which are countable real spectrum, finite dimensional eigenspaces for non-zero eigenvalues, and the convergence to zero of eigenvalues when the number of them are infinite (see, e.g., \cite{ArbogastBona08,brezis2010functional,OdenDemkowicz10, Royden2017-yo}). One of the important consequences is the Hilbert-Schmidt theorem  (see, e.g., \cite{ArbogastBona08,brezis2010functional,OdenDemkowicz10, Royden2017-yo,Showalter77}), which is a generalization of \cref{coro:spectralFinite}.
\begin{theorem}[Hilbert-Schmidt theorem for self-adjoint compact operators]
Let $\B:\Xs \to \Xs$ be a self-adjoint compact operator. Then there exists an orthonormal set of
eigen-functions $\varphi_i$ corresponding to non-zero eigenvalues
$\lambda_i$ of $\B$ such that for any $x \in \Xs$ we have a unique
expansion of the form
\begin{equation}
\x = \sum_{i} \LRp{\varphi_i, \x}_{\Xs}\varphi_i + \P\x,
\label{eq:spectralExpansion}
\end{equation}
where  $\P$ is an orthogonal projection from $\Xs$ to the
nullspace $\N\LRp{\B}$.

Furthermore, we have
\[
\B x = \sum_{i}\lambda_i \LRp{\varphi_i, x}_{\Xs} \varphi_i,
\]
that is, the set of all eigenfunctions $\LRc{\varphi_i}$ forms a basis for $\Range\LRp{\B}$.
\label{theo:HilbertSchmidt}
\end{theorem}

\begin{remark}
The fact that $\B$ is compact when $\Xs$ is a finite dimensional space implies that \cref{theo:HilbertSchmidt} is a generalization of \cref{coro:spectralFinite}. Indeed, let $\dim\LRp{\Xs} = n$. In this case, the eigenspace corresponding to the zero eigenvalue is spanned by finite number of (say $d$) orthonormal eigenfunctions $\LRc{\varphi_j^0}_{j=1}^d$ and thus $\P\x = \sum_{j=1}^d \LRp{\varphi_j^0, x}_{\Xs}\varphi_j^0$, which can be absorbed into the first sum on the right side of \cref{eq:spectralExpansion} so that we can write
\[
\x = \sum_{i=1}^n \LRp{\varphi_i, \x}_{\Xs}\varphi_i,
\]
after renaming the eigenfunctions corresponding to zero eigenvalues. This is exactly \cref{coro:spectralFinite}.
\end{remark}

We now follow the exposition in \cref{sect:svdFinite} to construct the SVD for compact operators. Let $\A: \Xs \to \Ys$ be a compact operator, then it can be shown that
 that $\B := \A^* \A: \Xs \to
  \Xs$ is a self-adjoint compact operator \cite{OdenDemkowicz10, ArbogastBona08}.  The Hilbert-Schmidt 
\cref{theo:HilbertSchmidt} says that there exists orthonormal
eigenfunctions $\varphi_i$ corresponding to nonzero eigenvalues $\lambda_i$ of $\B$ such that
\[
\B\varphi_i = \A^* \A \varphi_i = \lambda_i \varphi_i,
\]
which implies
\[
\LRp{\A^* \A \varphi_i, \varphi_i}_{\Xs} = \lambda_i \LRp{\varphi_i,\varphi_i}_{\Xs},
\]
which, by definition of $\A^*$, in turn can be written as
\[
\nor{\A \varphi_i}^2_{\Xs} = \lambda_i \nor{\varphi_i}^2_{\Xs},
\]
 which shows that $\lambda_i
\ge 0$. Let us define the \emph{singular value} $\sigma_i$ of $\A$ as
\begin{equation}
\sigma_i := \sqrt{\lambda_i}.
\label{eq:singularValues}
\end{equation}
We are now in the position to study the singular value decomposition for compact operators (see, e.g.,\cite{ColtonKress98}) that is a direct extension of \cref{theo:SVDfinite}.
\begin{theorem}[Singular value decomposition for compact operators]
Let $\LRc{\sigma_i}$ be the sequence of non-zero singular values (defined in \cref{eq:singularValues}) of a
compact operator $\A: \Xs \to \Ys$ and be ordered as
\[
\sigma_1 \ge \sigma_2 \ge\hdots,
\]
then there exist two orthonormal sequences $\LRc{\varphi_i}$ and $\LRc{\phi_i}$ such that
\begin{enumerate}
\item $\A\varphi_i = \sigma_i\phi_i$ and $\A^*\phi_i = \sigma_i\varphi_i$.
\item $\forall \x \in \Xs$, we have
\[
\x = \sum_i\LRp{\x, \varphi_i}_{\Xs}\varphi_i + \P\varphi,
\]
$\P: \Xs \to \N\LRp{\A}$ is an orthogonal projection.
\item There holds
\[
\A\x = \sum_i\sigma_i\LRp{\x, \varphi_i}_{\Xs}\phi_i.
\]
We call $\LRc{\sigma_i, \varphi_i,\phi_i}$, $i = 1, 2,\hdots$, {\bf the singular system} of $\A$.
\end{enumerate}
\label{theo:SVD}
\end{theorem}
\begin{proof}
The proof of this theorem is similar to its finite dimensional counterpart \cref{theo:SVDfinite}. The key that we exploit is the Hilbert-Schmidt theorem \cref{theo:HilbertSchmidt}.
\begin{enumerate}
\item By \cref{theo:HilbertSchmidt} we have
\[
\A^*\A\varphi_i = \sigma_i^2\varphi_i,
\]
where $\LRc{\varphi_i}$ is an orthonormal set in $\Xs$. Let us define
\[
\sigma_i\phi_i := \A\varphi_i,
\]
then
\[
\LRp{\phi_i,\phi_j}_{\Ys} = \frac{1}{\sigma_i\sigma_j}\LRp{\A\varphi_i,\A\varphi_j}_{\Ys} =
\frac{1}{\sigma_i\sigma_j}\LRp{\A^*\A\varphi_i,\varphi_j}_{\Xs} =   \frac{\sigma_i}{\sigma_j}\LRp{\varphi_i,\varphi_j}_{\Xs} = \delta_{ij}.
\]
That is, $\LRc{\phi_i}$ is an orthonormal set in $\Ys$. By definition we have
\[
\A^*\phi_i = \frac{1}{\sigma_i} \A^*\A\varphi_i = \sigma_i \varphi_i.
\]
\item Again, by Hilbert-Schmidt \cref{theo:HilbertSchmidt} we have
\[
\forall \x \in \Xs: \quad \x = \sum\LRp{\x, \varphi_i}_{\Xs}\varphi_i + \P\varphi,
\]
where $\P$ is an orthonormal projection from $\Xs$ to
$\N\LRp{\A^*\A}$. The second assertion is now clear owing to the fact that $\N\LRp{\A}
= \N\LRp{\A^*\A}$.
\item We start with the partial sum
\[
s_N := \sum_{i = 1}^N\LRp{\x,\varphi_i}_{\Xs}\varphi_i,
\]
and thus
\[
\A s_N = \sum_{i = 1}^N\sigma_i\LRp{\x,\varphi_i}_{\Xs}\phi_i.
\]
Now passing to the limit we obtain
\[
\lim_{N \to \infty}\A s_N = \A\LRp{\x - \P\x} = \A\x.
\]
Consequently,
\[
\A\x = \sum_i\mu_i\LRp{\x,\varphi_i}_{\Xs}\phi_i.
\]
\end{enumerate}
\end{proof}

The next result \cite{ColtonKress98}, due to Picard, tells us the conditions  under which inverting a compact operator is well-defined. {\em As we will see, the adjoint plays a key role.}
\begin{theorem}[Picard]
Suppose $\A: \Xs \to \Ys$ is a compact operator. The equation
\[
\A\x = \y
\]
is solvable iff
\begin{itemize}
\item[i)] $\y \in \N\LRp{\A^*}^\perp$, and
\item[ii)] $\sum\frac{1}{\sigma^2_i}\snor{\LRp{\y,\phi_i}_{\Ys}}^2 < \infty$,
\end{itemize}
where $\LRc{\sigma_i,\varphi_i,\phi_i}$ is the singular system of $\A$. In this case the solution is given by
\[
\x = \sum_i \frac{1}{\sigma_i}\LRp{\y,\phi_i}_{\Ys}\varphi_i.
\]
\label{theo:Picard}
\end{theorem}
\begin{proof}
    The SVD  \cref{theo:SVD} provides a simple proof for this theorem.
\begin{itemize}
\item[$\Rightarrow$ ] Solvability implies that $\y$ belongs to the
  range of $\A$, i.e. $\y \in \Range\LRp{\A}$. From the closed range \cref{theo:CRT} \footnote{Note that $\Range\LRp{\A}$ cannot be closed since $\A$
    is compact. Assume, on the contrary, it is, then by the bounded
    inverse theorem \cite{Rudin73} we know that $\A^{-1}$ is
    continuous and hence $\I = \A^{-1}\A$ is also a compact
    operator. But, this is a contradiction since identity operator in
    infinite dimensional space cannot be a compact operator
    \cite{OdenDemkowicz10, ArbogastBona08}. If $R\LRp{\A}$ were closed,
    then \emph{i)} would be both necessary and sufficient. Since this is not
    true for compact operators, we have to replace the closedness by
    the smooth property \emph{ii)} of the right hand side.}
   we know that $\Range\LRp{\A}
  \subset \overline{\Range\LRp{\A}} = \N\LRp{\A^*}^\perp$, and hence \emph{i)} holds. On the other
  hand, by \cref{theo:HilbertSchmidt} we can express a solution $\x$ as
\[
\x = \sum_i\LRp{\x, \varphi_i}_{\Xs}\varphi_i + \P\x,
\]
which, together with the Parseval identity, implies
\[
\nor{\x}_\Xs^2 = \sum_i\snor{\LRp{\x, \varphi_i}_{\Xs}}^2 + \nor{\P\x}_\Xs^2,
\]
which in turns implies
\[
\sum_i\snor{\LRp{\x, \varphi_i}_{\Xs}}^2 \le \nor{\x}_\Xs^2 < \infty.
\]
Since 
\[
\LRp{\x, \varphi_i}_{\Xs} = \frac{1}{\sigma_i}\LRp{\x,
  \A^*\phi_i}_{\Xs} = \frac{1}{\sigma_i}\LRp{\A\x, \phi_i}_\Ys =\frac{1}{\sigma_i}\LRp{\y, \phi_i}_\Ys,
\]
the assertion \emph{ii)} holds. 
\item[$\Leftarrow$ ]
Since $\y \in \N\LRp{\A^*}^\perp$,  Hilbert-Schmidt \cref{theo:HilbertSchmidt} gives
\[
\y = \sum_i\LRp{\y,\phi_i}_\Ys\phi_i.
\]
Now, from \emph{ii)} the following definition
\[
\x := \sum_i\frac{1}{\sigma_i}\LRp{\y, \phi_i}_\Ys\varphi_i
\]
is meaningful.
Together with \cref{coro:compactContinuous}, we have
\[
\A\x = \sum_i\frac{1}{\sigma_i}\LRp{\y, \phi_i}_\Ys\A\varphi_i = \sum_i\LRp{\y, \phi_i}_\Ys\phi_i = \y,
\]
where we have used in the last equality the Hilbert-Schmidt \cref{theo:HilbertSchmidt} for $\B = \A\A^*$, the fact that $\phi_i$ are eigenfunctions of $\B$, and $\y \in \N\LRp{\A^*}^\perp$.
This concludes the proof.
\end{itemize}
\end{proof}

{We now discuss the important consequence of the Picard
\cref{theo:Picard}}, that is, {\em inverting a compact operator is an ill-posed problem}. We observe that
\[
\y = \A\x = \sum_i\sigma_i\LRp{\x,\phi_i}_\Xs\phi_i,
\]
where we have used the second assertion of  \cref{theo:SVD}.
Since $\A$ is compact, and hence $\sigma_i\to 0$ as $i\to \infty$, $\A$ smoothes out the
contribution from the ``high frequency" mode: i.e. $\varphi_i$ for large $i$. In other words, the
output $\y$ is insensitive to high frequency modes $\varphi_i$ when $i \to \infty$.

Conversely,  let us perturb the right hand side $\y$ as
\[
\tilde{\y} = \y + \delta \phi_N,
\]
where $\delta \in \real$ and $N \in \mathbb{N}$,
then the corresponding solution reads
\[
\tilde{\x} = \sum_i\frac{1}{\mu_i}\LRp{\tilde{\y},
  \phi_i}_\Ys\varphi_i  = \x + \frac{\delta}{\mu_N}\varphi_N.
\]
Thus,
\[
\frac{\nor{\tilde{\x}-\x}_\Xs}{\nor{\tilde{\y}-\y}_\Ys} = \frac{1}{\mu_N} \to \infty, \text{ as } N \to \infty,
\]
\emph{which is exactly the subtle instability problem of inverting a
  compact operator, namely, small changes in the input can lead to
  very large change in the solution}. Most of linear inverse problems (such as deconvolution) fall into this category, and practical nonlinear inverse problems do too (see, e.g., \cite{Bui-ThanhGhattas12f, Bui-ThanhGhattas12, Bui-ThanhGhattas12a} and the references therein).

\begin{definition}[Well-posedness]
In Hadamard's sense \cite{hadamard02}, the problem $\A\x = g$ is well-posed if 
\begin{enumerate}
\item $\A$ is surjective (\emph{there exists a solution}: {\bf existence}), 
\item $\A$ is injective (\emph{there is at most one solution}: {\bf uniqueness}), and
\item $\A^{-1}$ is continuous (\emph{the solution depends continuously on the data}: {\bf stability}).
\end{enumerate}
\label{defi:wellposedness}
\end{definition} 

\begin{example}[Inverse of the fundamental theorem of calculus]
Let $\A: \Xs \to \Ys$ and consider the fundamental theorem of calculus in the following form
\[
\y\LRp{t} = \A\x := \int_0^t \x\LRp{s}\,ds, \quad 0\le t \le 1,
\]
and the inverse problem is to find $\x$ given its anti-derivative $\y$.
We are going to show that, depending on $\Xs, \Ys$, this inverse problem can be ill-posed or well-posed.
\begin{itemize}
\item First let us consider $\Xs = \Cs\LRp{\LRs{0,1}}, \Ys = \Cs\LRp{\LRs{0,1}}$. 
Let $\A\x = \y$ and consider
\[
\tilde{\y} := \y - \frac{\alpha}{N} + \frac{\delta}{N}\cos\LRp{N t}.
\]
Then, by the fundamental theorem of calculus, the corresponding solution is given by
\[
\tilde{\x} = \x + \delta \sin\LRp{N t}.
\]
Clearly
\[
\nor{\tilde{\y} - \y}_{\Cs\LRp{ \LRs{0,1}}} 
 := \sup_{t \in \LRs{0,1}}\snor{\tilde{\y}(t) - \y(t)}\to 0, \quad \text{ as } N \to \infty,
\]
but
\[
\nor{\tilde{\x} - \x}_{\Cs\LRp{ \LRs{0,1}}} = \alpha \quad \forall N.
\]
We conclude that $\A$ does not distinguish $\x$ and $\tilde{\x}$,
and as the result the inverse problem does not have a unique
solution. In fact, $\A$ is a compact operator\footnote{By the Ascoli-Arzela
  theorem \cite{OdenDemkowicz10, ArbogastBona08}.} and, as we have discussed above, it ``smoonthes'' out the difference in $\x$ and
$\tilde{\x}$ so that the observation $\y$ is the same. Intuitively, a compact operator
``squeezes'' its domain into ``smaller" range: for the above example
$\A$, as an integral operator, maps $\Cs\LRp{\LRs{0,1}}$ into $\Cs^1\LRp{\LRs{0,1}}
\subset \Cs\LRp{\LRs{0,1}}$. Since the inverse of a compact operator is unbounded, inverting the fundamental theorem of calculus is unstable
by the Picard \cref{theo:Picard}. The setting $\Xs = \Cs\LRp{\LRs{0,1}}, \Ys = \Cs\LRp{\LRs{0,1}}$ thus leads to an ill-posed problem.

\item Now let us consider $\Xs = \Cs \LRp{\LRs{0,1}}, \Ys = \Cs ^1\LRp{\LRs{0,1}}$. In this case we have
\[
\nor{\tilde{\y} - \y}_{\Cs^1\LRp{\LRs{0,1}}} := \nor{\tilde{\y} - \y}_{\Cs \LRp{\LRs{0,1}}} + \nor{\tilde{\y}' - \y'}_{\Cs\LRp{ \LRs{0,1}}} = \alpha, \quad \forall N,
\]
and since
\[
\nor{\tilde{\x} - \x}_{\Cs \LRp{\LRs{0,1}}} = \alpha \quad \forall N.
\]
we conclude that a small change in $\y$ leads to a small (in fact the
same) change in $\x$. The inverse problem is thus stable. The
uniqueness is also trivial due to the fact  that $\dd{\y}{t} = \x$.
The surjectivity is also clear. Consequently, the inverse problem is
well-posed in the Hadamard's sense.\footnote{Note that this is an
  instance of the Tikhonov theorem \cite{RoyCouchm01} since $\Cs^1\LRp{\LRs{0,1}}$ is
  compactly embedded in $\Cs\LRp{\LRs{0,1}}$.}

\end{itemize}
\end{example}

\begin{remark}
In practice, we do not solve $\A\x = \y$ directly on the infinite dimensional setting but via some discretization approach to obtain a finite dimensional problem to solve (on computer). This does not go around the ill-posedness issue. Indeed,
in this case, the compactess  of $\A$ is manifest in the ill-conditioning of its discrete counterpart whose smallest singular value could be very small. Inverting the discrete system is thus an ill-conditioned problem\textemdash a discrete way of saying ill-posedness.    
\end{remark}

We have seen that there could be multiple (or there is no) solutions to the linear
problem of interest $\A\x = \y$.
The main reason is that the nullspace of $\A$ is non-trivial or $\y$
is not in the range of $\A$. The question is if we can find a ``useful solution" in this case? One way to address this question is to look for the
solution that minimizes the residual, such as the least squares problem in \cref{coro:linearLS}:
\begin{equation}
\min_{\x}\half\nor{\A\x - \y}_\Ys^2. 
\label{eq:minNorm}
\end{equation}
However, when $\A: \Xs \to \Ys$ is compact, the ill-posedness nature of our inverse problem does not
go away as the normal equation \cref{eq:LSsolution} is still ill-posed due to the fact that $\A^*\A$ is compact. In other words, we still have problem with the
uniqueness  if $\N\LRp{\A}$ is not trivial, and the (bigger) problem with
instability due to inverting the compact operator $\A^*\A$. 
However, the optimization idea paves the way
for using optimization technique to overcome the ill-posedness problem, as we now discuss. Note that
the objective function in \cref{eq:minNorm} is quadratic in $\x$, and hence a
``parabola". Clearly, if it  is a
well-behaved parabola, then the minimizer is unique. This immediately
suggests that we should add a quadratic term to the objective function to improve its behavior, and hence removing the uniqueness issue: as will be shown, this also addresses the stability. This is
essentially the idea behind the {\em Tikhonov regularization} \cite{Tikhonov1977-me,Tikhonov1995-bm}, which
proposes to solve the following nearby problem
\begin{equation}
\min_{\x \in \Xs}\half\nor{\A\x - \y}_\Ys^2 + \frac{\kappa}{2} \nor{\x - \x_0}_\Xs^2,
\label{eq:minNormTik}
\end{equation}
where $\x_0$ is some ``prior'' reference function and $\kappa$ is
known as the \emph{regularization parameter}.
To show that the
\emph{regularized optimization} problem \cref{eq:minNormTik} is
well-posed, we need the  projection theorem
\cref{theo:projection} and the following key result from the Riesz-Fredholm theory
\cite{ColtonKress83}.
\begin{lemma}
Let $\A$ be a compact operator from $\Xs$ to $\Xs$. If $\LRp{I + \A}$
is injective, then $\LRp{I + \A}$ is continuously invertible.
\label{lem:RF}
\end{lemma}
\begin{theorem}
For any $\kappa > 0$, the regularized optimization problem \cref{eq:minNormTik} is well-posed.
\end{theorem}
\begin{proof}
Without loss of generality, assume $\kappa = 1$. We begin by
rewrite the optimization \eqnref{minNormTik} into the following equivalent form
\begin{equation}
\eqnlab{minNormTikLS}
\min_{\rr}\half\nor{\B \rr - w}_{\Ys \times \Xs}^2,
\end{equation}
where we have defined $\B: \Xs \ni \rr \mapsto \LRs{\A, \I}\rr:=\LRs{\A\rr, \rr} \in
\Ys \times \Xs$, and $w := \LRs{w_1, w_2}:=\LRs{\y, \rr_0}$. 
The inner product of $z,w \in \Ys \times \Xs$ is defined as $\LRp{z,w}_{\Ys \times \Xs} := \LRp{z_1,w_1}_\Ys + \LRp{z_2,w_2}_\Xs$, and the induced norm for any $z = \LRs{z_1,z_2} \in \Ys \times \Xs$ is given by $\nor{z}_{\Ys
  \times \Xs}^2 := \nor{z_1}_\Ys^2 + \nor{z_2}_\Xs^2$. 
From the definition of the inner product in $\Ys
  \times \Xs$, the definition of adjoint, and the fact that the identity operator $\I$ is self-adjoint,  we have $\B^*z = {\A^*z_1 + z_2}$. Next, from
\cref{coro:linearLS} we know that the minimizer satisfies
\[
\B^*\B \rr = \B^*w,
\]
which is equivalent to
\[
\LRp{\A^*\A + I}\rr = \A^*\y + \rr_0.
\]
Since $\LRp{\A^*\A + I}$ is injective\footnote{From $\LRp{\A^*\A + I}\rr = 0$ we have $ 0 =\LRp{\x,\LRp{\A^*\A + I}\rr}_\Xs = \nor{\A\x}_\Ys^2 + \nor{\x}^2_\Xs$ and thus $\x = \theta$.},
\cref{lem:RF} shows that it is continuously invertible,
i.e. $\nor{\LRp{\A^*\A + I}^{-1}} < \infty$. Hence,
\[
\nor{\rr}_\Xs = \nor{\LRp{\A^*\A + I}^{-1}\LRp{\A^*\y + \rr_0}}\le \nor{\LRp{\A^*\A + I}^{-1}}\LRp{\nor{\A^*} \nor{\y}_\Xs + \beta \nor{\rr_0}_\Xs},
\]
that is, the solution $\x$ of the Tikhonov regularization \cref{eq:minNormTik} is not only unique but also depends continuously on the data $\y$, and this concludes the proof.
\end{proof}

\subsection{Application of adjoint in the wellposedness of linear operator equation}
\label{sect:PDEapplication}

In this section we are interested in the well-posedness (in the sense of Hadamard in  \cref{defi:wellposedness}) of operator equation $\A\x = \y$, where $\A: \Xs \to \Ys$ is linear and continuous. Our exposition follows \cite{ErnGuermond04} closely. We begin with a key result \cite{ArbogastBona08,Abramovich2002-jt, ErnGuermond04}.

\begin{lemma}
    Let $\A: \Xs \to \Ys$ be linear and continuous. Then 
    \[
    \A \text{ is bounded below } \Leftrightarrow \exists \alpha > 0: \nor{\A\u}_\Ys \ge \alpha \nor{\u}_\Xs \Leftrightarrow
    \begin{cases}
    \A \text{ is injective}, \\
    \Range\LRp{\A} \text{ is closed}.
    \end{cases}
    \]
    \label{lem:boundedBelow}
\end{lemma}
\begin{proof}
    We provide a proof in \cref{proof:boundednessBelow}.
\end{proof}

The following result highlights  the role of the adjoint operator $\A^*$ on the injectivity and the closedness of $\Range\LRp{\A}$, and hence the boundedness below of $\A$.
\begin{theorem}
    Let $\A: \Xs \to \Ys$ be linear and continuous. The following are equivalent:
    \begin{enumerate}[label = \arabic*)]
        \item $\A^*: \Ys \to \Xs$ is surjective.
        \item $\A$ is injective and $\Range\LRp{\A}$ is closed.
        \item There exists $\alpha > 0$ such that
        \[
        \nor{\A\u}_\Ys \ge \alpha \nor{\u}_\Xs, \quad \forall \u \in \Xs.
        \]
        \item There exists $\alpha > 0$ such that
        \[
        \inf_{\u\in \Xs}\sup_{\v\in \Ys}\frac{\LRp{\A\u,\v}_\Ys}{\nor{\u}_\Xs\nor{\v}_\Ys} \ge \alpha.
        \]
    \end{enumerate}
    \label{theo:surjectiveAdjoint}
\end{theorem}
\begin{proof}
    We only need to show $1) \Leftrightarrow 2)$ as the equivalence between $2)$ and $3)$ is due to \cref{lem:boundedBelow}, and $4)$ is simply a restatement of $3)$. We have
\[
    \begin{array}{cr}
    \A^* \text{ is surjective} & \\
    \Updownarrow & \\
    \Range\LRp{\A^*} = \Xs \text{ and thus } \Range\LRp{\A^*} \text{ closed} & \\
    \Updownarrow & \text{The closed range \cref{theo:CRT}}\\
    \N\LRp{\A} = \Range\LRp{\A^*}^\perp = \LRc{\theta} \text{ and } \Range\LRp{\A} \text{ closed } & \\
    \Updownarrow & \\
    \A \text{ is injective and } \Range\LRp{\A} \text{ closed.} & 
    \end{array}
\]
\end{proof}
The following twin counterpart of \cref{theo:surjectiveAdjoint} characterizes the surjectivity of $\A$ via the adjoint $\A^*$.
\begin{theorem}
    Let $\A: \Xs \to \Ys$ be linear and continuous. The following are equivalent:
    \begin{enumerate}[label = \arabic*)]
        \item $\A: \Ys \to \Xs$ is surjective.
        \item $\A^*$ is injective and $\Range\LRp{\A^*}$ is closed.
        \item There exists $\alpha > 0$ such that
        \[
        \nor{\A^*\v}_\Xs \ge \alpha \nor{\v}_\Ys, \quad \forall \v \in \Ys.
        \]
        \item There exists $\alpha > 0$ such that
        \[
        \inf_{\v\in \Ys}\sup_{\u\in \Xs}\frac{\LRp{\u,\A^*\v}_\Xs}{\nor{\u}_\Xs\nor{\v}_\Ys} \ge \alpha.
        \]
    \end{enumerate}
    \label{theo:surjectiveForward}
\end{theorem}

Combining \cref{theo:surjectiveAdjoint} and \cref{theo:surjectiveForward} we see that $\A$ is bijective iff $\A^*$ is bijective. The more popular statement that leads to the 
Banach-Ne$\check{\text{c}}$as-Babu$\check{\text{s}}$ka
 theorem for the variational equation is the following
\begin{lemma}
    Let $\A: \Xs \to \Ys$ be linear and continuous. The following are equivalent:
    \begin{enumerate}[label = \arabic*)]
        \item $\A$ is bijective
        \item 
        \begin{itemize}
            \item $\exists \alpha > 0$ such that  $
        \inf_{\u\in \Xs}\sup_{\v\in \Ys}\frac{\LRp{\A\u,\v}_\Ys}{\nor{\u}_\Xs\nor{\v}_\Ys} \ge \alpha,
        $ and
            \item If $\LRp{\A\u,\v}_\Ys = 0, \forall \u \in \Xs,$ then $ \v = \theta$.
        \end{itemize}
    \end{enumerate}
    \label{lem:bijectivityForward}
\end{lemma}
\begin{proof}
    The proof is straightforward. Indeed, the first statement of $2)$ is the injectivity of $\A$ plus the closedness of $\Range\LRp{\A}$ due to \cref{theo:surjectiveAdjoint}, and the second statement of $2)$ can be written equivalently in terms of $\A^*$ as: ``if $\LRp{\u,\A^*\v}_\Ys = 0, \forall \u \in \Xs,$ then $ \v = \theta$", which is equivalent to ``if $\A^*\v = \theta$ then $ \v = \theta$", which in turn simply means $\N\LRp{\A^*} = \LRc{\theta}$, which then means $\A$ is surjective owing to \cref{theo:surjectiveForward}.
\end{proof}

\begin{example}
    Consider $\Am:\real^n\mapsto \real^m$. We are interested in applying \cref{lem:bijectivityForward}  to find conditions for the linear system of equations $\Am \ub=\yb$ to have a unique solution. To that end, we suppose $\Am$ is bijective. Thus, both $\Am$ and $\Am^*$ are injective. By the rank-nullity theorem \cref{eq:RankNullity} we have
    \[n=dim(\N(\Am))+dim(\Range(\Am)). \]
    Since $A$ is injective, we have $dim(\N(\Am))=0$, and hence
    \[n=dim(\Range(\Am))\leq m.\]
Following similar arguments for the injectivity of $\Am^*$, we have:
\[ m=dim(\R(\Am^*))\leq n.\]
Therefore, it is necessary that $n=m$ for the bijectivity of $\Am$. Further,  the inf-sup condition in \cref{lem:bijectivityForward} says:
\[ 0<\alpha \leq \inf_{\ub\in \real^n}\sup_{\vb\in \real^n} \frac{\LRp{\Am\ub,\vb}}{\nor{\ub}_{\real^n}\nor{\vb}_{\real^n}}\leq \inf_{\ub\in \real^n} \frac{\nor{\Am\ub}_{\real^n}}{\nor{\ub}_{\real^n}}=\sigma_{min}(\Am).\]
\[ \implies 0<\alpha\leq \sigma_{min}(\Am),\]
where $\sigma_{min}(\Am)$ denotes the smallest singular value of $\Am$. We conclude the necessary and sufficient for a linear system of equations $\Am \ub=\yb$ to have a unique solution is that the matrix $\Am$ is square and invertible. This is consistent with what we know from linear algebra.
\end{example}

\begin{theorem}[Banach-Ne$\check{\text{c}}$as-Babu$\check{\text{s}}$ka]
Let $\A: \Xs \to \Ys$ be the unique linear and continuous operator (see \cref{fn:bilinear}) associated 
with a continuous sequilinear form $a: \Xs\times\Ys \to \Fs$ such that
\[
\LRp{\A\u,\v}_\Ys := a\LRp{\u,\v}, \quad \forall \u \in \Xs \text{ and } \v \in \Ys,
\]
where $\snor{a\LRp{\u,\v}} \le \beta \nor{\u}_\Xs\nor{\v}_\Ys$ and $0 < \beta <\infty$.
The following are equivalent:
\begin{enumerate}[label = \arabic*)]
    \item For all $\y \in \Ys$, there exists a unique solution $\u \in \Xs$ such that 
    \[
    a\LRp{\u,\v} = \LRp{\y,\v}_\Ys, \quad \forall \v \in \Ys.
    \]
    \item There exists $\alpha > 0$ such that
    \begin{enumerate}[label = \bf{C\arabic*})]
            \item \label{enum:infsup} $\exists \alpha > 0$ such that  $
        \inf_{\u\in \Xs}\sup_{\v\in \Ys}\frac{a\LRp{\u,\v}}{\nor{\u}_\Xs\nor{\v}_\Ys} \ge \alpha,
        $ and
            \item \label{enum:AdjointInjectivity} If $a\LRp{\u,\v} = 0, \forall \u \in \Xs,$ then $ \v = \theta$.
        \end{enumerate}
\end{enumerate}
Furthermore, when either of the statements holds then the unique solution is stable in the following sense:
\[
\nor{\u}_\Ys \le \frac{1}{\alpha} \nor{\y}_\Ys.
\]
\label{theo:BNB}
\end{theorem}
\begin{proof}
    The proof is obvious due to the definition $\LRp{\A\u,\v}_\Ys := a\LRp{\u,\v}$, and thus the equivalent of the variational equation $a\LRp{\u,\v} = \LRp{\y,\v}_\Ys$ and $\A\u = \y$.
    Specifically, owing to \cref{lem:bijectivityForward}, statement $2)$ is equivalent to the bijectivity of $\A$. The stability of the solution $\u$ is the direction consequence of the boundedness from below of $\A$:
    \[
    \alpha \nor{\u}_\Xs \le \nor{\A\u}_\Ys = \sup_{\v \in \Ys} \frac{\LRp{\A\u,\v}_\Ys}{\nor{\v}_\Ys}  = \sup_{\v \in \Ys} \frac{a\LRp{\u,\v}_\Ys}{\nor{\v}_\Ys} =   \sup_{\v \in \Ys} \frac{\LRp{\y,\v}_\Ys}{\nor{\v}_\Ys} = \nor{\y}_\Ys.
    \]
\end{proof}

\begin{remark}
The condition $
\inf_{\u\in \Xs}\sup_{\v\in \Ys}\frac{a\LRp{\u,\v}}{\nor{\u}_\Xs\nor{\v}_\Ys} \ge \alpha,
$  is known as the inf-sup condition, and, as we have shown, it is nothing more than the restatement of the boundedness from below of the associated linear operator $\A$ or equivalently the injectivity of $\A$ plus its closed range.
\end{remark}

\begin{example}
    Consider solving the diferential equation $\u' := \dd{\u}{\x}=\f$ in $(0,1)$ with $\u(0)=0$ and $f\in \Ls^2(0,1)$. The corresponding weak form of the problem is formulated as: seek $\u\in \Hs_0^1(0,1) := \{ \u\in \Ls^2(0,1),\ \ \u'\in \Ls^2(0,1), \ \u(0)=0\}$ such that:
    \[\LRp{\u',\v}_{\Ls^2}=\LRp{\f,\v}_{\Ls^2},\quad \forall \v\in \Ls^2(0,1).\]
    We choose
    $\Xs=\Hs_0^1(0,1)$, $\Ys=\Ls^2(0,1)$ and $\Fs = \real$, and are going to use \cref{theo:BNB} to show that the differential equation is well-posed.
    
    The continuity of the bilinear form $a(\u,\v)$ is clear as
    \[\snor{a(\u,\v)}=\snor{(\u',\v)}\leq \nor{\u'}_{\Ls^2}\nor{\v}_{\Ls^2}\leq \nor{\u}_{\Hs^1}\nor{\v}_{\Ls^2}.\]

We next verify the inf-sup condition \cref{enum:infsup}. By a simple integration and the Cauchy-Schwarz inequality (see also \cref{exa:SLI}) we obtain the following the Poincar\'e-Friedrichs inequality
\[
\nor{\u'}_{\Ls^2}\geq\ \nor{\u}_{\Ls^2},
\]
and thus the inf-sup condition holds since
\[
   \inf_{\u\in \Hs^1}\sup_{\v\in \Ls^2}\frac{a(\u,\v)}{\nor{\u}_{\Hs^1}\nor{\v}_{\Ls^2}} = \inf_{\u\in \Hs^1} \frac{\nor{\u'}_{\Ls^2}}{\nor{\u}_{\Hs^1}} \ge \frac{1}{\sqrt{2}}.
\]

Now we verify the injectivity of the adjoint, i.e \cref{enum:AdjointInjectivity}.
We start from
    \begin{equation}
    \LRp{\u',\v}_{\Ls^2}=0,\ \ \forall \u\in \Hs_0^1\LRp{0,1}.
    \label{eq:adjointInj}
    \end{equation}
Since $\Cs^\infty_0(0,1)\subset \Hs_0^1(0,1)$, we have:
    \[ \LRp{\psi',\v}_{\Ls^2}=0,\quad \forall \psi \in \Cs^\infty_0(0,1).\]
    By definition of the distributional derivative we arrive at
    \[\LRa{\psi,\v'}_{\Ls^2}=0,\quad \forall \psi\in \Cs^\infty_0(0,1),\]
    which implies that $\v$ is a constant function. Now in \eqref{eq:adjointInj} taking $\u = \x$ we have
    \[
    \int_{0}^1 \v dx = 0,
    \]
    which means $\v = 0$. Thus, the differential equation is well-posed with the setting $\Xs=\Hs_0^1(0,1)$, $\Ys=\Ls^2(0,1)$.
\end{example}

\begin{example}[Friedrichs' system]
    We consider the abstract problem $\A\u = \y$ where $\A$ is the  Friedrichs' operator defined in \cref{exa:Friedrichs}. We choose $\Ys = \Ls^2\LRp{\Omega}$ and 
    \[
    \Xs = \LRc{\x \in \Hs_\A: \LRp{\Dm - \Mm}\x = 0 \text{ on } \pOmega},
    \]
    with the inner product $\LRp{\u,\w}_\Xs := \LRp{\u,\w}_\Ys + \LRp{\B\u,\B\w}_\Ys$, and hence the induced graph norm $\nor{\u}_\Xs = \sqrt{\nor{\u}_\Ys^2 + \nor{\B\u}_\Ys^2}$. We can show that $\A: \Xs \to \Ys$ is bijective (see \cite[Theorem 5.7]{ErnGuermond04}), and thus applying \cref{theo:BNB} shows that $\A\u = \y$ is well-posed for any $\y \in \Ys$. The beauty here is that this single proof is applicable for a large class of PDEs \cite{Friedrichs58, ErnGuermond06, Jensen04, ErnGuermondCaplain07}.
\label{exa:FriedrichsWellposedness}
\end{example}

When $\Ys = \Xs$ and  the sequilinear form is symmetric, the inf-sup condition is both necessary and sufficient for bijectivity.
\begin{lemma}
    Consider the variational equation $a\LRp{\u,\v} = \LRp{\y,\v}_\Xs,$ $ \forall \v \in \Xs$, with a continuous sequilinear form $a:\Xs\times \Xs \to \Fs$. Suppose $a\LRp{\cdot,\cdot}$ is symmetric, i.e., $a\LRp{\w,\v} = \overline{a\LRp{\v,\w}}$. Then, there exists a unique solution $\u$  iff
    \[
    \exists \alpha > 0 \text{ such that }  
        \inf_{\u\in \Xs}\sup_{\v\in \Ys}\frac{a\LRp{\u,\v}}{\nor{\u}_\Xs\nor{\v}_\Ys} \ge \alpha.
    \]
\end{lemma}
\begin{proof}
    We need to prove only the second condition in the second statement of \cref{theo:BNB}, namely the injectivity of the adjoint in \cref{enum:AdjointInjectivity}. But this is obvious due to symmetry:
    \[
    0 = \overline{a\LRp{\w,\v}} = a\LRp{\v,\w}, \quad \forall \w \in \Xs \implies 0 = 
    \sup_{\w \in \Xs} \frac{a\LRp{\v,\w}}{\nor{\w}_\Xs} \ge \alpha \nor{\v}_\Xs \implies \v = \theta.
    \]
\end{proof}
\begin{remark}
    Note that the symmetry of $a\LRp{\w,\v}$ is equivalent to the self-adjointness of its associate linear and continuous operator $\A$ defined in \cref{theo:BNB}. Indeed, since $a:\Xs\times \Xs \to \Fs$, we have
    \[
    \LRp{\A\w,\v}_\Xs = a\LRp{\w,\v} = \overline{a\LRp{\v,\w}} = \overline{\LRp{\A\v,\w}}_\Xs = \LRp{\w, \A\v}_\Xs,
    \]
    which means $\A^* = \A$.
\end{remark}

On the other hand, when $\Ys = \Xs$ and the sequilinear form is coercive, the condition for bijectivity is simpler.
\begin{lemma}[The Lax-Milgram lemma]
Consider the variational equation $a\LRp{\u,\v} = \LRp{\y,\v}_\Xs,$ $ \forall \v \in \Xs$, with a continuous sequilinear form $a:\Xs\times \Xs \to \Fs$.  If
\begin{equation}
    \tag{Coercivity} \a\LRp{\v,\v} \ge \alpha \nor{\v}_\Xs^2, \quad \forall \v \in \Xs,
    \label{eq:coercivity}
\end{equation}
then there exists a unique solution $\u$ and $\nor{\u} \le \frac{1}{\alpha}\nor{\y}_\Xs$.
\label{lem:LaxMilgram}
\end{lemma}
\begin{proof}
    We need to verify the two conditions in the second statement of \cref{theo:BNB}. The inf-sup condition in \cref{enum:infsup} is clear from \cref{eq:coercivity} as
    \[
    \alpha \nor{v}_\Xs \le \frac{\a\LRp{\v,\v}}{\nor{v}_\Xs} \le \sup_{\w \in \Xs} \frac{\a\LRp{\v,\w}}{\nor{w}_\Xs}.
    \]
    For the injectivity of the adjoint in \cref{enum:AdjointInjectivity}, we note that
    \[
    a\LRp{\w,\v} = 0, \quad \forall \w \in \Xs \implies 0 = \sup_{\w \in \Xs} a\LRp{\w,\v} \ge a\LRp{\v,\v} \ge \alpha \nor{\v}^2_\Xs \implies \v = \theta.
    \]
\end{proof}

When the sequilinear form is symmetric and positive, it turns out that \cref{eq:coercivity} is both sufficient and necessary.
\begin{corollary}
    Consider the variational equation $a\LRp{\u,\v} = \LRp{\y,\v}_\Xs,$ $ \forall \v \in \Xs$, with a continuous sequilinear form $a:\Xs\times \Xs \to \Fs$, with continuity constant $\beta$. Suppose $a\LRp{\cdot,\cdot}$ is symmetric, i.e., $a\LRp{\w,\v} = \overline{a\LRp{\v,\w}}$ and positive, i.e., $a\LRp{\v,\v} > 0, \forall \v \ne \theta$. Then, there exists a unique solution $\u$ iff the coercivity condition \cref{eq:coercivity} holds.
\end{corollary}
\begin{proof}
    \cref{lem:LaxMilgram} proves the sufficiency, and we need to show the necessity. Since the sequilinear form $a\LRp{\cdot,\cdot}$ is symmetric and positive, it defines an inner product in $\Xs$ and the induced norm is 
    \[
    \nor{\v}_a := \sqrt{a\LRp{\v,\v}}.
    \]
    Thus, by the inf-sup condition, the Cauchy-Schwarz inequality, and the continuity of $a\LRp{\cdot,\cdot}$, we have
    \[
    \alpha \nor{v}_\Xs\le \sup_{\w\in \Xs} \frac{a\LRp{\w,\v}}{\nor{\w}_\Xs} \le \sup_{\w\in \Xs} \frac{\sqrt{a\LRp{\v,\v}}\sqrt{a\LRp{\w,\w}}}{\nor{\w}_\Xs} \le \sqrt{\beta}  \sqrt{a\LRp{\v,\v}},
   \]
    and this ends the proof.
\end{proof}
\begin{example}
    We now consider the variational equation $a\LRp{\u,\v} = \LRp{\y,\v}_\Ys$ where the bilinear form $a\LRp{\u,\v}$, the spaces $\Xs$ and $\Ys$, and other specifications are described in \cref{exa:ellipticOpWeak}.
We have shown that $a\LRp{\u,\v}$ is symmetric and continuous on $\Hs^1_0\LRp{\Omega}$. It is clearly positive if we assume that $z$ is bounded away from $-\infty$. What remains is to show that $a\LRp{\u,\v}$ is coercive. Recall the Poincar\'e-Friedrichs inequality (see, e.g., \cite{ErnGuermond04,ArbogastBona08,brezis2010functional}) for $\Hs^1_0\LRp{\Omega}$: there exists a constant $c$ depending on only $\Omega$ such that
\[
c\nor{\u}_{\Ls^2} \le \nor{\Grad\u}_{\Ls^2}, \quad \forall \u \in \Hs^1_0\LRp{\Omega}.
\]
It follows that for any $\v \in \Hs^1_0\LRp{\Omega}$ we have
\[
a\LRp{\v,\v} = \LRp{e^\z\Grad\v,\Grad\v}_{\Ls^2}^2 \ge e^{\inf\z}\nor{\Grad\v}_{\Ls^2}^2 \ge \min\LRc{1,c}\half e^{\inf\z} \nor{\v}_{\Hs^1}^2,
\]
and this ends the proof.
    \label{exa:ellipticPDEs}
\end{example}

\subsection{Application of adjoint to the Sturm-Liouville problem and generalized Fourier series}
\label{sect:SLapplication}

In this section, we are interested in self-adjoint operator $\A$ in infinite dimensions.
Our exposition mostly follows \cite{Showalter98}. To the end of this section, for any $\A \in \L\LRp{\Xs,\Ys}$ we assume that its domain $\Do\LRp{\A}\subseteq \Us$ is dense in $\Xs$.

\begin{definition}[Closed linear operators]
A linear operator $\A \in \L\LRp{\Xs,\Ys}$ is called closed iff its graph 
\[
\mc{G}_\A :=\LRc{\LRs{\u,\A\u}: \u \in \D\LRp{\A}}
\]
is closed.
\end{definition}
Clearly, any continuous linear operator is necessarily closed.
The following are standard results for closed operators \cite{Showalter98}.
\begin{lemma}
The following hold:
\begin{itemize}

    \item If $\A$ is closed, so is $\A^*$.
    \item $\A$ is closed and $\Do\LRp{\A} = \Us$ iff $\A \in \B\LRp{\Us,\Vs}$. 
    \item If $\Do\LRp{\A} = \Us$, then $\A^*$ is continuous, and hence $\Do\LRp{\A^*}$ is closed.
    \item If $\A$ is closed, then $\Do\LRp{\A^*}$ is dense in $\Vs$.
\end{itemize}
\label{lem:closedOptMore}
\end{lemma}

We consider operators with $\Xs = \Ys$ in this section.  Suppose $\Vs$ is dense in $\Xs$ and the injection $\Vs \to \Xs$ is compact, and for simplicity in writing we denote $\Vs \xhookrightarrow[\text{dense}]{\text{compact}} \Xs$. Let $\A$ be a linear operator $\A :\Do\LRp{\A} \subset \Vs \xhookrightarrow[\text{dense}]{\text{compact}} \Xs \to \Xs$ with the domain defined as
\[
\Do\LRp{\A} := \LRc{\x \in \Vs: \A\x \in \Xs}.
\]
We also assume that the associate sequilinear form $a\LRp{\u,\v} := \LRp{\A\u,\v}_\Xs$ is defined for $\u \in \Do\LRp{\A}$ and $\v \in \Vs$, and is continuous on $\Vs$. Furthermore, we assume that the $a\LRp{\u,\v}$ is $\Vs$-elliptic in the following sense: there exists $c > 0 $ such that
\[
a\LRp{\v,\v} + \overline{a\LRp{\v,\v}} \ge 2 c \nor{\v}_{\Vs}^2. 
\]
Finally, we assume that $a\LRp{\cdot,\cdot}$ is symmetric, i.e.,
\[
a\LRp{\u,\v} = \overline{a\LRp{\v,\u}}, \quad \forall \u,\v \in \Vs,
\]
which is, again, equivalent to the self-adjointness of $\A$ on $\Vs$. 
Recall that \cref{lem:realEigen} holds in this case: in particular, eigenvalues of $\A$ are real and its eigenfunctions corresponding to distinct eigenvalues are orthogonal to each other. The following theorem provides further characteristics of  eigenpairs of $\A$.
\begin{theorem}
    Suppose all the aforementioned assumptions hold for $\A$. Then, there is a countable sequence of eigenpairs $\LRc{\lambda_n,\v_n}_{n = 1}^\infty$ such that 
    \begin{itemize}
        \item $\A\v_n = \lambda_n \v_n$, where $\nor{\v}_\Xs = 1$,
        \item $\LRp{\v_n,\v_m}_\Xs = 0$ for all $n \ne m$,
        \item $0 < \lambda_1 \le \lambda_2 \le \hdots \le \lambda_n \to \infty$ when $n \to \infty$, and 
        \item $\LRc{\v_n}_{n = 1}^\infty$ is a basis of $\Xs$.
    \end{itemize}
    \label{theo:generalizedFourier}
\end{theorem}
\begin{proof}
    We sketch the proof and more details can be found in \cite{Showalter77}. The keys are: i) the $\Vs$-elliptic condition implies coercivity of $\A$ on $\Vs$. By Lax-Milgram \cref{lem:LaxMilgram}, $\A: \Vs \to \Vs' \equiv \Vs$ is a continuous bijection. The restriction of $\A$ on $\Do\LRp{\A}$ is thus surjective on $\Xs$, and $\A^{-1}: \Xs \to \Do\LRp{\A} \subset \Vs$ exists and is continuous by the coercivity; ii) the compact injection of $\Vs$ in $\Xs$ then implies that $\A^{-1}: \Xs \to \Xs$ is compact.  Hilbert-Schmidt \cref{theo:HilbertSchmidt} then ensures  the existence of the eigenpairs $\mu_n,\v_n$ of $\A^{-1}$, where $\LRc{\v_n}_n$ is a basis of $\Do\LRp{\A}$, and hence the eigenpairs $\lambda_i,\v_n$ of $\A$ where $\lambda_n = 1/\mu_n$; and iii) the $\Vs$-elliptic condition also implies  $\A$ is closed operator and this leads to the conclusion that $\LRc{\v_n}_{n = 1}^\infty$ is a also basis of $\Xs$.
\end{proof}
\cref{theo:generalizedFourier} provides a generalized Fourier series in $\Xs$. In particular, for any function $\f \in \Xs$, we have
\[
\f = \sum_{n = 1}^\infty \LRp{\f,\v_n}_\Xs\v_n,
\]
where equality means convergence in the topology generated by the $\Xs$-norm. 
Below are a few special cases that lead to the standard Fourier series and a general Sturm-Liouville problem.
\begin{example}
    Consider the following Sturm-Liouville problem: seek $\lambda$ and $\v$ such that
    \begin{equation}
    \begin{cases}
    - \frac{d^2 \v}{d \x^2} =  \lambda \v & \text{ in } \Omega = \LRp{0,1}, \\
\v = 0 & \text{ on } \pOmega = \LRc{0,1},
    \end{cases}
    \label{eq:SLI}
    \end{equation}
where the derivative is understood in the classical sense sense. We solve this problem using \cref{theo:generalizedFourier}. To that end, we define the symmetric sequilinear form as 
\[
a\LRp{\u,\v} := \LRp{\dd{\u}{\x},{\dd{\v}{\x}}}_{\Ls^2\LRp{\Omega}} = \int_0^1 \dd{\u}{\x}\overline{\dd{\v}{\x}}\,d\x,
\]
for all $\u, \v \in \Vs$ with $\Vs := \Hs^1_0\LRp{0,1}$ and thus the derivatives in $a\LRp{\cdot,\cdot}$ are understood in the weak sense. Note that $\Vs \xhookrightarrow[\text{dense}]{\text{compact}} \Xs := \Ls^2\LRp{0,1}$ \cite{Showalter98}, and the continuity on $\Vs$ of $a\LRp{\cdot,\cdot}$ is straightforward.  By the fundamental theorem of calculus and Cauchy-Schwarz inquality we can easily arrive at a Poincar\'e-Friedrichs inequality: $\forall \v \in \Cs^1_0\LRs{0,1}$
\[
\nor{\dd{\v}{\x}}_{\Ls^2\LRp{0,1}} \ge \nor{\v}_{\Ls^2\LRp{0,1}},
\]
which also holds for any $\v \in \Hs^1_0\LRp{0,1}$ due to the density of $\Cs^1_0\LRp{0,1}$ in $\Vs$. This leads to the $\Vs$-ellipticity of $a\LRp{\u,\v}$ as
\[
a\LRp{\v,\v} + \overline{a\LRp{\v,\v}} \ge \nor{\v}_{\Hs^1_0\LRp{0,1}}^2.
\]

Thus, \cref{theo:generalizedFourier} ensures that there is a complete orthonormal eigenfunctions of $\A: \Do\LRp{\A}\to \Ls^2\LRp{0,1}$ in $\Ls^2\LRp{0,1}$, where $\Do\LRp{\A} := \LRc{\w \in \Hs^1_0\LRp{0,1}: \A\w \in \Ls^2\LRp{0,1}}$ and $\LRp{\A\u,\v}_{\Ls^2\LRp{0,1}} := a\LRp{\u,\v}$ for all $\u \in \Do\LRp{\A}$ and $\v \in \Hs^1_0\LRp{0,1}$. 
Let us determine what $\Do\LRp{\A}$ is. We have, by definition of $a\LRp{\cdot,\cdot}$,
\[
\begin{array}{clcr}
\LRp{\A\w,\v}_{\Ls^2\LRp{0,1}}     &=& \int_0^1\dd{\w}{\x}\overline{\dd{\v}{\x}}\,d\x  & \forall \v \in \Hs^1_0\LRp{0,1}\\
& \Updownarrow & &  \text{ density of } \Cs_0^\infty\LRp{0,1} \text{ in } \Hs^1_0\LRp{0,1}\\
 \LRp{\A\w,\varphi}_{\Ls^2\LRp{0,1}}    &= &   \int_0^1\dd{\w}{\x}\overline{\dd{\varphi}{\x}}\,d\x & \forall \varphi \in \Cs_0^\infty\LRp{0,1}\\
   &\Updownarrow &  & \text{definition of distributional derivative}\\
\LRp{\A\w,\varphi}_{\Ls^2\LRp{0,1}}    &= &  \LRa{- \frac{d^2 \w}{d \x^2},\varphi}   & \forall \varphi \in \Cs_0^\infty\LRp{0,1}\\
&\Updownarrow &  & \\
\A\w &=& - \frac{d^2 \w}{d \x^2} & \text{ in } \Ls^2\LRp{0,1}
\end{array}
\]
As a result, $\Do\LRp{\A} =  \Hs^2_0\LRp{0,1}$.

Each eigenpair $\lambda_n$ and $\v_n \in \Do\LRp{\A}$ satisfies
\[
\A\v_n = \lambda_n \v_n \text{ in } \Ls^2\LRp{0,1},
\]
i.e,
\begin{equation}
- \frac{d^2 \v_n}{d \x^2}  =  \lambda_n \v_n, \text{ in } \Ls^2\LRp{0,1},
\label{eq:eigenSLI}
\end{equation}
which is equivalent to
\[
\v_n = - \int\int\lambda_n\v_n\,d\x\,d\y,
\]

Owing $\v_n \in \Do\LRp{\A}$, and the embedding 
of $\Hs^1_0\LRp{0,1}$ in $
\Cs_0\LRp{0,1}$ (see, e.g., \cite[Theorem 8.2]{brezis2010functional})
 the eigenvalue problem \cref{eq:eigenSLI} holds in the classical sense, which is exactly the Sturm-Liouville problem \cref{eq:SLI}. Thus, eigenfunctions of the Sturm-Liouville problem \cref{eq:SLI} forms a complete basis for $\Ls^2\LRp{0,1}$. By a simple integrations, we have $\v_n = \sqrt{2}\sin\LRp{n\pi \x}$, $n = 1, 2, \hdots$, which is exactly a Fourier basis ({\em sine series}) for $\Ls^2\LRp{0,1}$.
\label{exa:SLI}
\end{example}

\begin{corollary}
    Let $\Vs$ and $\Xs$ be given in \cref{theo:generalizedFourier}. Suppose the sequilinear form $a\LRp{\cdot,\cdot}$ is symmetric  and continuous on $\Vs$. Assume that there exists $c > 0 $ such that for some $\lambda \in \real$
\[
a\LRp{\v,\v} + \overline{a\LRp{\v,\v}} + 2\lambda\nor{\v}_\Xs^2 \ge 2 c \nor{\v}_{\Vs}^2, \quad \forall \v \in \Vs.
\]
Then, there exists an orthonormal sequence of eigenfunctions of $\A$, which is a basis for $\Xs$ and the corresponding eigenvalues satisfies $-\lambda < \lambda_1 \le \lambda_2 \le \hdots \le \lambda_n \to \infty$ when $n \to \infty$.
\label{coro:generalizedFourierCoro}
\end{corollary}
\begin{proof}
    Let us define the sequilinear  $b\LRp{\u,\v} := a\LRp{\u,\v} + \lambda \LRp{\u,\v}_\Xs$. The continuity and symmetry of $b\LRp{\cdot, \cdot}$ are clear. The linear operator $\B$ associated with $b\LRp{\u,\v} $ 
    is given by $\B\u = \A\u + \lambda\u$ for any $\u \in \Do\LRp{\B} \equiv \Do\LRp{\A}$, and $b\LRp{\u,\v} = \LRp{\B\u,\v}_{\Xs}$ for all $\u\in \Do\LRp{\B}$ and $\v \in \Vs$. Furthermore, we have
    \[
    b\LRp{\v,\v} + \overline{b\LRp{\v,\v}} \ge 2 c \nor{\v}_{\Vs}^2, \quad \forall \v \in \Vs,
    \]
    i.e., $b\LRp{\u,\v}$ is $\Vs$-elliptic. \cref{theo:generalizedFourier} thus applies to $b\LRp{\cdot,\cdot}$. In particular, there exist eigenpairs $\LRc{\gamma_n,\v_n}_{n=1}^\infty$ of $\B$ such that $\A\v_n + \lambda\v_n = \B\v_n = \gamma_n\v_n$ and $0 < \gamma_1\le \gamma_2\le \hdots \le \gamma_n \to \infty$ when $n \to \infty$. As a result, $\LRc{\lambda_n,\v_n}_{n=1}^\infty$ with $\lambda_n := \gamma_n -\lambda$, and thus $-\lambda < \lambda_1 \le \lambda_2 \le \hdots \le \lambda_n \to \infty$ when $n \to \infty$, and $\LRc{\v_n}_{n=1}^\infty$  being a basis of $\Xs$, and this concludes the proof.
\end{proof}

\begin{remark}
    We can use \cref{coro:generalizedFourierCoro} for \cref{exa:SLI} as well. In particular, we can choose $\Vs = \Hs^1\LRp{0,1}$, and $b\LRp{\u,\v} := a\LRp{\u,\v} + \lambda\LRp{\u,\v}_{\Ls^2\LRp{0,1}}$ is coercive on $\Vs$ for any $\lambda > 0$ with $c = \min\LRc{1,\lambda}$. The rest of the arguments are similar and the same results are obtained.
\end{remark}

\begin{example}
    Consider the following Sturm-Liouville problem: seek $\kappa$ and $\v$ such that
    \begin{equation}
    \begin{cases}
    - \frac{d^2 \v}{d \x^2} =  \kappa \v & \text{ in } \Omega = \LRp{0,1}, \\
\dd{\v}{\x}= 0 & \text{ on } \pOmega = \LRc{0,1},
    \end{cases}
    \label{eq:SLII}
    \end{equation}
where the derivative is understood in the classical sense. The sequilinear $a\LRp{\cdot,\cdot}$ is defined the same as in \cref{exa:SLI} with $\Vs := \Hs^1\LRp{0,1}$ and $\Xs := \Ls^2\LRp{0,1}$. We still have $\Vs \xhookrightarrow[\text{dense}]{\text{compact}} \Xs := \Ls^2\LRp{0,1}$ \cite{Showalter98}. The linear operator $\A$ associated with $a\LRp{\cdot,\cdot}$ defined via the identity $\LRp{\A\u,\v}_\Xs = a\LRp{\u,\v}$ for all $\u \in \Do\LRp{\A}:=\LRc{\w \in \Vs: \A\w \in \Xs \text{ and } \dd{\w}{\x}= 0  \text{ on } \pOmega = \LRc{0,1}}$, and  $\v \in \Vs$. Note that unlike \cref{exa:SLI} in which the boundary conditions are naturally incorporated in $\Vs$, we have to build  the boundary conditions in the definition of the domain of $\A$ in order to associate $a\LRp{\cdot,\cdot}$ with the eigenvalue problem \cref{eq:SLII}. 
The continuity and symmetry of $a\LRp{\cdot,\cdot}$ on $\Vs$ are clear.
  Furthermore, $a\LRp{\cdot,\cdot}$ satisfies \cref{coro:generalizedFourierCoro} for any $\lambda > 0$ and $c = \min\LRc{1,\lambda}$. Using a similar argument as in \cref{exa:SLI}, the eigenpairs $\lambda_n,\v_n$ of $\A$ are exactly the solutions of \cref{eq:SLII}, and in particular $\v_n = \sqrt{2}\cos\LRp{n\pi\x}$ for $n \ge 1$ and $\v_0 = 1$. By \cref{coro:generalizedFourierCoro}, $\LRc{\v_n}_{n=0}^\infty$ is another Fourier basis ({\em cosine series}) for $\Ls^2\LRp{0,1}$.
  \label{exa:SLII}
\end{example}

\begin{example}[A more general Sturm-Liouville problem]
We now generalize \cref{exa:SLI} and \cref{exa:SLII}: seek $\kappa$ and $\v$ such that
    \begin{equation}
    \begin{cases}
    \frac{1}{\rho}\LRs{\dd{}{\x}\LRp{p\dd{\v}{\x}} +q \v}=  \kappa \v & \text{ in } \Omega = \LRp{0,1}, \\
\alpha\v\LRp{0}+ \beta\dd{\v}{\x}\LRp{0}= 0, & \gamma\v\LRp{1}+ \delta\dd{\v}{\x}\LRp{1}= 0,
    \end{cases}
    \label{eq:generalSL}
    \end{equation}
where $\rho \in \Cs\LRs{0,1}$ and $\rho > 0$, $p \in \Cs^1\LRs{0,1}$ and $p < 0$, and $q \in \Cs\LRs{0,1}$. Here, the constants $\alpha, \beta, \gamma,$  and $\delta$ satisfy $\alpha^2 + \beta^2 \ne 0$, and $\gamma^2 + \delta^2 \ne 0$. 
We take $\Xs := \Ls^2_\rho\LRp{0,1}$ is the $\Ls^2\LRp{0,1}$ space with the weighted inner product: $\LRp{\u,\v}_{\Ls^2\LRp{0,1},\rho} := \LRp{\rho\u,\v}_{\Ls^2\LRp{0,1}}$.
We choose the sequilinear form $a\LRp{\cdot,\cdot}$ as
\[
a\LRp{\u,\v} := \LRp{p\dd{\u}{\x},\dd{\v}{\x}}_{\Ls^2\LRp{0,1}} + \LRp{q\u,\v}_{\Ls^2\LRp{0,1}},
\]
with the derivative understood in the weak sense. By taking $\Vs = \Hs^1_\rho\LRp{0,1}$, where $\Hs^1_\rho\LRp{0,1}$ is $\Hs^1\LRp{0,1}$  based on $\Ls^2_\rho\LRp{0,1}$, it is clear  that $a\LRp{\cdot,\cdot}$ is continuous and  symmetric on $\Vs$. The linear operator $\A$ associated with $a\LRp{\cdot,\cdot}$ is defined via the identity $\LRp{\A\u,\v}_\Xs = a\LRp{\u,\v}$ for all $\u \in \Do\LRp{\A}:=\LRc{\w \in \Hs^1_\rho\LRp{0,1}: \A\w \in \Xs, \alpha\w\LRp{0}+ \beta\dd{\w}{\x}\LRp{0}= 0 \text{ and }  \gamma\w\LRp{1}+ \delta\dd{\w}{\x}\LRp{1}= 0}$, and  $\v \in \Vs$. By a similar distributional argument as in \cref{exa:SLI}, one can show that $\A\w = \frac{1}{\rho}\LRs{\dd{}{\x}\LRp{p\dd{\w}{\x}} +q \w}$ in $\Ls^2\LRp{0,1}$. Furthermore, it is clear that $a\LRp{\cdot,\cdot}$ satisfies \cref{coro:generalizedFourierCoro} for any $\lambda > - q$ and $c = \min\LRc{\nor{\frac{\rho}{p}}_\infty^{-1},\nor{\frac{\rho}{\lambda+q}}_\infty^{-1}}$. We thus conclude that the eigenpairs $\lambda_n$ and $\v_n \in \Do\LRp{\A}$ satisfy
\begin{equation}
\frac{1}{\rho}\LRs{\dd{}{\x}\LRp{p\dd{\v_n}{\x}} +q \v_n} = \lambda_n\v_n,
\label{eq:eigenSLIIgeneral}
\end{equation}
which is equivalent to
\[
\v_n = \int\frac{1}{p}\int\rho\LRp{\lambda_n - q}\v_n\,d\x\,d\y,
\]
which, together with the embedding result in
\cite[Theorem 8.2]{brezis2010functional}, shows that
 the eigenvalue problem \cref{eq:eigenSLIIgeneral} holds in the classical sense, which is exactly the Sturm-Liouville problem \cref{eq:generalSL}. Thus, eigenfunctions of the Sturm-Liouville problem \cref{eq:generalSL} forms a complete basis for $\Ls^2\LRp{0,1}$. However, it is not so clear how to calculate the eigenpairs analytically as we have done for the previous examples. 
    \label{exa:generalSL}
\end{example}

\begin{example}
Next, we consider the  following Sturm-Liouville problem that is not covered by \cref{exa:generalSL}: seek $\kappa$ and $\v$ such that
    \begin{equation}
    \begin{cases}
    - \frac{d^2 \v}{d \x^2} =  \kappa \v & \text{ in } \Omega = \LRp{0,1}, \\
\v(0) = \v(1) &\text{ and } \dd{\v}{\x}\LRp{0}= \dd{\v}{\x}\LRp{1},
    \end{cases}
    \label{eq:SLIII}
    \end{equation}
where the derivative is understood in the classical sense. The sequilinear $a\LRp{\cdot,\cdot}$ is defined the same as in \cref{exa:SLI} with $\Vs := \Hs^1\LRp{0,1}$ and $\Xs := \Ls^2\LRp{0,1}$. The linear operator $\A$ associated with $a\LRp{\cdot,\cdot}$ defined via the identity $\LRp{\A\u,\v}_\Xs = a\LRp{\u,\v}$ for all $\u \in \Do\LRp{\A}:=\LRc{\w \in \Vs: \A\w \in \Xs, \v(0) = \v(1) \text{ and } \dd{\v}{\x}\LRp{0}= \dd{\v}{\x}\LRp{1}}$, and  $\v \in \Vs$. Similar to \cref{exa:SLII}, $a\LRp{\cdot,\cdot}$ is obviously symmetric and continuous on $\Vs$. Furthermore, $a\LRp{\cdot,\cdot}$ satisfies \cref{coro:generalizedFourierCoro} for any $\lambda > 0$ and $c = \min\LRc{1,\lambda}$. Using a similar distributional argument as in \cref{exa:SLI}, the eigenpairs $\lambda_n,\v_n$ of $\A$ are exactly the solutions of \cref{eq:SLIII}, and in particular $\v_0 ={1}$ and $\v_n \in \LRc{\sqrt{2}\cos\LRp{2\n\pi\x}, \sqrt{2}\sin\LRp{2\n\pi\x}}$ for $\n\ge 1$, which 
is the usual Fourier basis for $\Ls^2\LRp{0,1}$.
    \label{exa:SLIII}
\end{example}

From the above examples, a few observations are in order. First, 
the view of (generalized) Fourier series from self-adjoint Sturm-Liouville operators immediately provides a rigorous convergence guarantee for the (generalized) Fourier series in the $\Ls^2$-sense. This view also shows that there are other orthogonal bases for $\Ls^2(a,b)$ and we can in principle find them by solving the corresponding Sturm-Liouville eigenvalue problems. Second, the results also show that $\Ls^2(a,b)$ is a separable Hilbert space. Third, the results are not restricted to $
\Ls^2$ spaces over compact subsets in $\real$ but are also valid for any compact subsets in $\real^\n$ using tensor product, dilation, and translation (see, e.g., \cite{ArbogastBona08}). 

\subsection{Application of adjoint  in PDE-constrained optimization}
\label{sect:PDEconstrained}

We have seen the important role of adjoint in constrained optimization in \cref{sect:optimization}, especially constrained optimization with equality constraints that have separable structure (see \cref{coro:adjointReduced} and \cref{lem:reducedGradient}). We have also seen how the adjoint looks like and how it helps compute the gradient of deep neural networks (DNN) efficiently as backpropagation in \cref{sect:DNN}.  In this section, we shall work out the details of adjoint equation and the reduced gradient  for optimization problems constrained by partial differential equations (PDE-constrained optimization). We consider prototype steady state (time-independent) PDEs of elliptic and hyperbolic types. The goal is to show we translate the abstract results in \cref{lem:reducedGradient} to concrete problems. This can serve as the baseline for carrying out the same task for different PDE-constrained optimization problems. Other topics on PDE-constrained optimization can be found in \cite{DelosReyes2015,biegler2012large,biegler2007real,borzi2012computational,PDEConstrained11}.

\begin{example}[Advection-PDE-constrained optimization problem]
Consider the following PDE-constrained optimization problem
\[
\min_{\z,\u} \half\int_\Omega\u^2\,d\Omega
\]
\SubjectTo
\begin{subequations}
\begin{align}
\bs{\beta} \cdot \Grad \u &= 0, \quad \text{ in } \Omega, \\
\bs{\beta}\cdot\nb\u &= \z, \quad \text{ in } \pOmega_{\text{in}},
\end{align}
\label{eq:forwardTransportPDE}
\end{subequations}
where $\u \in \Hs^1_{\bs{\beta}}\LRp{\Omega}\equaldef\LRc{\u: \u \in
  \mbb{L}^2\LRp{\Omega} \text { and } \bs{\beta}\cdot \Grad \u \in
  \mbb{L}^2\LRp{\Omega}}$. See \cref{exa:advectionOp} for the definition of other quantities in the constraint and the associated differential operator together with its adjoint. This optimization problem is a special case of the abstract problem discussed in \cref{lem:reducedGradient}. 
  Note that for $\u \in
\Hs^1_{\bs{\beta}}\LRp{\Omega}$, its trace (in fact weighted trace with
weight $\snor{\bs{\beta}\cdot \nb}$) on $\pOmega_{\text{in}}$ belongs to
$\mbb{L}^2\LRp{\pOmega_{\text{in}}}$. The correct space for $\z$ is
thus $\mbb{L}^2\LRp{\pOmega_{\text{in}}}$ with a weighted norm (see \cref{fn:density}). Since the constraints map 
$\LRs{\u,\z} \in \Xs \times \Zs := \Hs^1_{\bs{\beta}}\LRp{\Omega} \times \Ls^2\LRp{\pOmega_{\text{in}}}$ to $\Ys \equaldef
\Ls^2\LRp{\Omega} \times \Ls^2\LRp{\pOmega_{\text{in}}}$, the
Lagrange multiplier $\y = \LRs{\v,\w}$ has two components $\v \in
\mbb{L}^2\LRp{\Omega}$ and $\w \in \mbb{L}^2\LRp{\pOmega_{\text{in}}}$,
respectively. Our task is to find the explicit form of the first order optimality condition \cref{eq:KKTinf} which, we 
recall, is a special case of the first order optimality condition via the Lagrangian multiplier \cref{theo:LagrangeMult}.
For practical PDE-constrained problem, the adjoint
operators $\LRs{\D_\u c\LRp{\u_0,\z_0}}^*\y$ and $\LRs{\D_\z
  c\LRp{\u_0,\z_0}}^*\y$ are subtlely coupled and we have to go back to the Lagrangian functional in \cref{theo:LagrangeMult}
  to derive the optimality condition. To that end, let us form the Lagrangian functional
\[
L\LRp{\z, \u} = \half\int_\Omega\u^2\,d\Omega + \int_\Omega \LRp{\bs{\beta} \cdot \Grad
\u}\v\,d\Omega + \int_{\pOmega_{\text{in}}}\LRp{\bs{\beta}\cdot\nb\u - \z}\,\w\,ds,
\]
 and note that our optimization variable has two components $\LRs{\u,\z} \in \Hs^1_{\bs{\beta}}\LRp{\Omega} \times \mbb{L}^2\LRp{\pOmega_{\text{in}}}$.
Take an arbitrary direction $\LRs{h,r} \in \Hs^1_{\bs{\beta}}\LRp{\Omega} \times  \mbb{L}^2\LRp{\pOmega_{\text{in}}}$, the first order optimality condition \cref{eq:firstOptInf}, {\em with $\LRs{\u,\z}$ in place of $\u_0$ and $\LRs{h,r}$ in place of $h$}, reads
\[
\LRa{\LRs{\v,\w},\D c\LRp{\LRs{\u,\z}, \LRs{h, r}}}_\Ys = \int_\Omega \LRp{\bs{\beta}\cdot\Grad h}
  \v\,d\Omega + \int_{\pOmega_{\text{in}}}\LRp{\bs{\beta}\cdot\nb h - r}\,\w\,ds,
\]
which after integration by parts becomes
\begin{multline*}
\LRa{\LRs{\v,\w},\D c\LRp{\LRs{\u,\z}, \LRs{h, r}}}_\Ys = -\int_\Omega \LRp{\bs{\beta}\cdot\Grad \v}
  h\,d\Omega + \int_{\pOmega_{\text{in}}}\bs{\beta}\cdot\nb\LRp{\w + \v}h\,ds + \\
\int_{\pOmega_{\text{out}}}\bs{\beta}\cdot\nb \v\,h\,ds - \int_{\pOmega_{\text{in}}}r\,\w\,ds.
\end{multline*}
 Here, we have restricted $\v$ in $\Hs^1_{\bs{\beta}}\LRp{\Omega}$
for the differential and integral operators to make sense. The first order optimality condition
\cref{eq:firstOptInf} in this case reads: $\forall \LRs{h,r} \in \Hs^1_{\bs{\beta}}\LRp{\Omega} \times \mbb{L}^2\LRp{\pOmega_{\text{in}}}$,
\[
\int_\Omega \u h\,d\Omega
-\int_\Omega \LRp{\bs{\beta}\cdot\Grad \v}
  h\,d\Omega + \int_{\pOmega_{\text{in}}}\bs{\beta}\cdot\nb\LRp{\w + \v}h\,ds + \int_{\pOmega_{\text{out}}}\bs{\beta}\cdot\nb \v\,h\,ds - \int_{\pOmega_{\text{in}}}r\,\w\,ds = 0,
\]
which, after taking $r = 0$ and any $h \in
\Hs^1_{\bs{\beta},0}\LRp{\Omega}\equaldef\LRc{\u \in
  \Hs^1_{\bs{\beta}}\LRp{\Omega}: \eval{\u}_{\pOmega} = 0}$, becomes
\[
\int_\Omega \u h\,d\Omega -\int_\Omega \LRp{\bs{\beta}\cdot\Grad \v}
h\,d\Omega = 0, \quad \forall h \in \Hs^1_{\bs{\beta},0}\LRp{\Omega},
\]
which implies\footnote{It is due to the fact that $H^1_{\bs{\beta},0}\LRp{\Omega}$ is dense
in $\mbb{L}^2\LRp{\Omega}$ assuming $\Omega$ has segment property \cite{AntonicBurazin09}.}
\[
{
-\bs{\beta}\cdot\Grad \v + \u = 0.
}
\]
Consequently, the first order optimality condition reduces to: $\forall \LRs{h, r} \in \Hs^1_{\bs{\beta}}\LRp{\Omega} \times \Ls^2\LRp{\pOmega_{\text{in}}},$
\begin{equation}
\int_{\pOmega_{\text{in}}}\bs{\beta}\cdot\n\LRp{\w + \v}h\,ds + \int_{\pOmega_{\text{out}}}\bs{\beta}\cdot\nb \v\,h\,ds - \int_{\pOmega_{\text{in}}}r\,\w\,ds = 0, 
\label{eq:OptimalityAdvectionReduced}
\end{equation}
which, by taking $h = 0$ on $\pOmega_{\text{out}}$ and $r = 0$, becomes
\[
\int_{\pOmega_{\text{in}}}\bs{\beta}\cdot\nb\LRp{\w + \v}h\,ds = 0,
\]
which in turn gives\footnote{Note that this is true due to the fact
  that the trace operator $\gamma: H^1_{\bs{\beta},0}\LRp{\Omega} \to
  \mbb{L}^2_{\bs{\beta}\cdot\n}\LRp{\pOmega}$ is a continuous
  surjection (see, e.g., \cite{Bui-ThanhDemkowiczGhattas11b}), where
  $\mbb{L}^2_{\bs{\beta}\cdot\nb}\LRp{\pOmega}$ is $\mbb{L}^2\LRp{\pOmega}$ with the weighted inner product $\LRp{\u,\v}_{\mbb{L}^2_{\bs{\beta}\cdot\nb}\LRp{\pOmega}} := \int_{\pOmega} \snor{\bs{\beta}\cdot\nb}\u\v\,ds$.
  \label{fn:density}
  }
\[
\w = -\v \text{ on } \pOmega_{\text{in}},
\]
that is, the adjoint variables are not independent. This can be then
substituted into \cref{eq:OptimalityAdvectionReduced} to further reduce the first order optimality condition to
\[
\int_{\pOmega_{\text{out}}}\bs{\beta}\cdot\nb \v\,h\,ds + \int_{\pOmega_{\text{in}}}r\,\v\,ds = 0, \quad \forall \LRs{h, r} \in \Hs^1_{\bs{\beta}}\LRp{\Omega} \times \Ls^2\LRp{\pOmega_{\text{in}}}.
\]
It follows that, by taking $r = 0$ and using the surjectivity in \cref{fn:density}, we conclude
\[
\bs{\beta}\cdot\nb \v = 0 \text{ on } \pOmega_{\text{out}},
\]
and thus
\[ 
\v = 0 \text{ on } \pOmega_{\text{in}}.
\]
In summary, 
 the control equation \cref{eq:control} becomes
\begin{equation}
\v = 0 \text{ on } \pOmega_{\text{in}},
\label{eq:controlTransport}
\end{equation}
and the adjoint equation \cref{eq:adjointEq} reads
\begin{subequations}
\begin{align}
-\bs{\beta}\cdot\Grad \v  &= -\u \text{ in } \Omega,\\
\bs{\beta}\cdot\nb \v &= 0 \text{ on } \pOmega_{\text{out}},
\end{align}
\label{eq:adjointTransportPDE}
\end{subequations}
Note that the differential operator on the left side of the adjoint equation (together with the homogeneous boundary condition) is exactly the adjoint operator we found in \cref{exa:advectionOp}, which is not a surprise.
As can be seen, the adjoint equation describes a reverse flow with $-\bs{\beta}$ velocity (compared to $\beta$ in the forward equation) with (the
derivative of) the objective function, particularly the forward solution $\u$, as the forcing. The control equation says
that at the optimal the forcing of the adjoint equation is such that the
adjoint solution $\v$ on $\pOmega_{\text{in}}$ must vanish. Clearly, one admissible solution is that 
the adjoint is identically zero and the forcing $u$ is
identically zero. It then follows from the forward equation that $\z =
0$. This is not surprising since, by inspection, the quadratic
optimization under consideration has a solution $\u = 0$ and $\z = 0$. The reduced gradient  
can be now computed for a given $\rr$ via three steps: 1) solve the
\emph{forward equation} \cref{eq:forwardTransportPDE} for $u\LRp{\rr}$, 2) solve
the \emph{adjoint equation} \cref{eq:adjointTransportPDE} for $\v\LRp{u(\rr),\rr}$, and 3) 
substitute $\v(\u\LRp{\z},\z)$ into the left hand side of \cref{eq:controlTransport} 
to obtain the
reduced gradient.
\label{exa:advectionControl}
\end{example}

\begin{example}[Elliptic-PDE-constrained optimization problem]
\[
\min_{\z,\u} J\LRp{\u} := \half \int_\Omega\LRp{\u-\u^{obs}}^2\,d\Omega
\]
\SubjectTo
\begin{subequations}
\begin{align}
-\Div\LRp{e^{\z}\Grad \u} &= 0, \quad \text{ in } \Omega, \\
\u &= g, \quad \text{ in } \pOmega,
\end{align}
\label{eq:forwardElliticPDEoptimalControl}
\end{subequations}
where $\u^{obs}\LRp{\xb}$ is some reference/observational data, and the definition of the operator associated with the constraint and its adjoint are given in
\cref{exa:ellipticOp}: in particular, $\rr \in \Cs^1\LRp{{\Omega}} \subset \Ls^2\LRp{\Omega}$. Here, $g$ is the Dirichlet boundary data.
This optimization is a special case of the abstract one in \cref{lem:reducedGradient}.
Thus, the constraint  maps $ \LRs{\u,\z} \in \Xs\times \Zs :=\Hs_\A \times
\Cs^1\LRp{{\Omega}}$ to $\Ls^2\LRp{\Omega} \times
\Ls^2\LRp{\pOmega}$, and the
Lagrange multiplier $\y = \LRs{\v,\w}$ has two components $\v \in
\mbb{L}^2\LRp{\Omega}$ and $\w \in \mbb{L}^2\LRp{\pOmega}$,
respectively. 
Similar to \cref{exa:advectionControl}, we go back to Lagrangian multiplier \cref{theo:LagrangeMult} to derive the explicit form of the first order optimality condition. In this case,
the Lagrangian reads
\[
L\LRp{\rr, \u} = J\LRp{\u} + \int_\Omega \LRs{-\Div\LRp{e^{\rr}\Grad \u}}\v\,d\Omega + \int_{\pOmega}\LRp{\u - g}\,\w\,ds,
\]
and note that our optimization variable has two components $\LRs{\u,\z} \in 
\Hs_\A \times
\Cs^1\LRp{{\Omega}}$.
Take an arbitrary direction $\LRs{h,r} \in 
\Hs_\A \times
\Cs^1\LRp{{\Omega}}$, the first order optimality condition \cref{eq:firstOptInf}, {\em with $\LRs{\u,\z}$ in place of $\u_0$ and $\LRs{h,r}$ in place of $h$}, reads
\begin{multline*}
\int_\Omega\LRp{\u-\u^{obs}}h\,d\Omega + \int_\Omega \LRs{-\Div\LRp{e^{\rr}\Grad h}}\v\,d\Omega + \int_{\pOmega}h\,\w\,ds \\
+\int_\Omega \LRs{-\Div\LRp{e^{\rr}r\Grad u}}\v\,d\Omega = 0.
\end{multline*}
We next restrict $\v \in \Hs_\A$ and integrate the
second term by parts two times  we arrive at: 
\begin{multline}
\int_\Omega\LRp{\u-\u^{obs}}h\,d\Omega + \int_\Omega \LRs{-\Div\LRp{e^{\rr}\Grad
    \v}}h\,d\Omega + \int_{\pOmega}h\,\w\,ds - \int_{\pOmega}e^\rr\Grad h\cdot\nb\,\v\,ds \\+
\int_{\pOmega}e^\rr\Grad \v\cdot\nb\,h\,ds + \int_\Omega
\LRs{-\Div\LRp{e^{\rr}r\Grad u}}\v\,d\Omega = 0, \quad \forall \LRs{h, r} \in \Hs_\A \times
\Cs^1\LRp{{\Omega}}.
\label{eq:firstVariationElliptic}
\end{multline}
Following a similar strategy\footnote{ 
Here we take $\h \in \Cs_0^\infty\LRp{\Omega}$ and $r = 0$ to obtain the equation \cref{eq:adjointEllipticPDE} for $\v$. To get the boundary condition \cref{eq:adjointEllipticPDEbc}, we then take $\h \in \Cs_0^1\LRp{\Omega}$ so
that the normal trace $\Grad h\cdot \nb$ is surjective on $\Ls^2\LRp{\pOmega}$.} as in 
\cref{exa:advectionControl} for \cref{eq:firstVariationElliptic} gives the adjoint equation
\begin{subequations}
\begin{align}
\label{eq:adjointEllipticPDE}
-\Div\LRp{e^{\rr}\Grad \v} &= -\LRp{\u-\u^{obs}}, \quad \text{ in }
\Omega, \\ 
\label{eq:adjointEllipticPDEbc}
\v &= 0, \quad \text{ in } \pOmega.
\end{align}
\label{eq:adjointEllipticPDEoptimalControl}
\end{subequations}
Note that the differential operator on the left side of the adjoint equation (together with the homogeneous boundary condition) is exactly the adjoint operator we found in \cref{exa:ellipticOp}, which is not a surprise.
The first order optimality condition is thus reduced to: $\quad \forall \LRs{h, r} \in \Hs_\A \times
\Cs^1\LRp{{\Omega}},$
\begin{equation}
 \int_{\pOmega}h\,\w\,ds +
\int_{\pOmega}e^\rr\Grad \v\cdot\nb\,h\,ds + \int_\Omega
\LRs{-\Div\LRp{e^{\rr}r\Grad u}}\v\,d\Omega = 0, 
\label{eq:firstVariationEllipticReduced}
\end{equation}
which, by taking $r = 0$,  gives, 
\[
\w = -e^\rr\Grad \v\cdot\nb.
\]
Similar to \cref{exa:advectionControl}, we see that the second component of the adjoint variable $\y$
depends on the first, hence $\v$ is in fact the only adjoint
variable. 
The first order optimality condition is further reduced to
\begin{equation*}
 \int_\Omega
\LRs{-\Div\LRp{e^{\rr}r\Grad u}}\v\,d\Omega = 0, \quad \forall r \in
\Cs^1\LRp{{\Omega}},
\end{equation*}
which\textemdash after integrating by parts, using the fact that $\v = 0$ on $\pOmega$, and using the fact that $\Cs^1\LRp{\Omega}$ is dense in $L^2\LRp{\Omega}$\textemdash gives
the control equation
\begin{equation}
e^{\rr}\Grad u\cdot\Grad \v = 0.
\label{eq:controlEllipticPDE}
\end{equation}

The reduced gradient  
can be now computed for a given $\rr$ via three steps: 1) solve the
\emph{forward equation} \cref{eq:forwardElliticPDEoptimalControl} for $u\LRp{\rr}$, 2) solve
the \emph{adjoint equation} \cref{eq:adjointEllipticPDEoptimalControl} for $\v\LRp{u(\rr),\rr}$, and 3) 
substitute $\u(\z)$ and $\v(\u\LRp{\z},\z)$ into the left hand side of \cref{eq:controlEllipticPDE} 
to obtain the
reduced gradient.
\end{example}

\section{Conclusions}
This paper has put together in one place the roles of adjoint in  various disciplines of mathematics, sciences, and engineering. The objective is to systematically compile these materials on the same mathematical footing starting from the basic definitions. We aim to provide a  unified perspective and understanding of adjoint.
This work could give broader views and better insights into the application of adjoint beyond a single community. We have established general results and then  specified them to each application with sufficient details including the connections among them. This paper is written as an interdisciplinary tutorial on adjoint with results and applications in both finite-dimensional and infinite-dimensional Hilbert spaces.
We have shown how adjoint can solve problems and facilitate progress in various fields through specific  examples including linear algebra (e.g. eigendecomposition and the singular value decomposition), ordinary differential equations (an epidemic model), partial differential equations (of elliptic, hyperbolic, and Friedrichs' types), neural networks (feed-forward deep neural networks), least squares and inverse problems (with Tikhonov regularization), and PDE-constrained optimization.

\section*{Acknowledgments}
An answer to a question from Dr. Jonathan Wittmer in early Fall 2022 on a 6:00am bus from Oak Hill to UT Austin about the usefulness of adjoint  initiated this work. We would like to thank Dr. Wittmer for the fruitful conversation. We also would like to thank Hai Nguyen of the \href{https://phoices.netlify.app/}{Pho-Ices Group} at the Oden Institute, UT Austin, for generating \cref{fig:FTA}.

\appendix

\section{Proofs of some results in the paper}
\label{sect:appendix}

    

\begin{proof}[Proof for \cref{exa:ellipticOp}]
     
     We next find $\A^*$ and $\Do\LRp{\A^*}$. Take any $\v \in \Do\LRp{\A^*}$, by \cref{defi:adjointInfinite} we have
     \begin{align*}
     \LRp{\A\u,\v}_{\Ls^2\LRp{\Omega}} &= \LRp{\u,\A^*\v}_{\Ls^2\LRp{\Omega}}, & \forall \u \in \Do\LRp{\A} &\\ 
      &\Updownarrow &  &\forall \varphi \in \Cs_0^\infty\LRp{\Omega} \text{ dense in } \Do\LRp{\A} \\
      \LRp{\A\varphi,\v}_{\Ls^2\LRp{\Omega}} &= \LRp{\varphi,\A^*\v}_{\Ls^2\LRp{\Omega}}  & & \\
      &\Updownarrow &  &\text{definition of distributional derivative} \\
     \LRa{\varphi, -\Div\LRp{e^{\z}\Grad \v}} &= \LRp{\varphi,\A^*\v}_{\Ls^2\LRp{\Omega}}  & \\
     &\Updownarrow & &  \\
     -\Div\LRp{e^{\z}\Grad \v} &= \A^*\v \in \Ls^2\LRp{\Omega} & 
     &
     \end{align*}
     which shows that
     $\A^* = \A$ on $\Do\LRp{\A^*}$.
     Next, from 
     \[
\LRp{\A\u,\v}_{\Ls^2\LRp{\Omega}} = \LRp{\u,\A^*\v}_{\Ls^2\LRp{\Omega}} =  \LRp{\u,\A\v}_{\Ls^2\LRp{\Omega}}, \quad \forall \u \in \Do\LRp{\A} \text{ and } \v \in \Do\LRp{\A^*},
     \]
and integration by parts we conclude
\[
-\LRp{e^\z\Grad\u \cdot \nb, \v}_{\Ls^2\LRp{\pOmega}} + \LRp{\u,e^\z\Grad\v\cdot\nb}_{\Ls^2\LRp{\pOmega}} = 0, 
\quad \forall \u \in \Do\LRp{\A} \text{ and } \v \in \Do\LRp{\A^*},
\]
with the assumption that the boundary integrals make sense. As a result, in addition to $\Div\LRp{e^{\z}\Grad \v} \in \Ls^2$ we need $\v \in \Hs^1\LRp{\Omega}$ and $\v = 0$ on $\pOmega$. We conclude that $\Do\LRp{\A^*} = \Do\LRp{\A}$, and $\A$ is self-adjoint.
\label{proof:closedEllipticOp}
\end{proof}

\begin{proof}[Proof for \cref{exa:Friedrichs}]
We now derive the adjoint $\A^*$ (and its domain) of the Friedrichs' operator in \cref{exa:Friedrichs}. Proceed as in \cref{proof:closedEllipticOp}, we can easily see that
\[
\A^*\v = -\sum_{k=1}^{\n}\partial_k \LRp{\Am_k\v} + \Cm^T \v \in \Ls^2\LRp{\Omega},
\]
for any $\v \in \Do\LRp{\A^*}$. Thus $\Do\LRp{\A^*} \subseteq \Hs_\A$. Next from 
%
\[
       \LRp{\A\u,\v}_{\Ls^2\LRp{\Omega}} = \LRp{\u,\A^*\v}_{\Ls^2\LRp{\Omega}}, \quad \forall \u \in \Do\LRp{\A} \text{ and } \v \in \Do\LRp{\A^*},
      \]
     and integrating by parts, we obtain
     \[
     \LRp{\u,\Dm\v}_{\Ls^2\LRp{\pOmega}} = 0, \quad \forall \u \in \Do\LRp{\A} \text{ and } \v \in \Do\LRp{\A^*},
     \]
     which in turns yields $\Dm\v \in \N\LRp{\Dm - \Mm}^\perp$ since $\u \in \N\LRp{\Dm-\Mm}$. We thus conclude
     \[
     \Do\LRp{\A^*} = \LRc{\v \in \Hs_\A: \Dm\v \in \LRs{\N\LRp{\Dm - \Mm}}^\perp}.
     \]
    \label{proof:adjointFriedrichsOp}
\end{proof}

\begin{proof}[Proof of \cref{lem:boundedBelow}]
    For the necessary, the injectivity is clear. Now, let $\LRc{\y_i}_{i=1}^\infty \subset \Range\LRp{\A}$ and $\y_i \xrightarrow[]{\Ys} \y$ and we need to show that $\y \in \Range\LRp{\A}$. 
    There exists $\LRc{\x_i}_{i=1}^\infty \subset \Xs: \y_i = \A\x_i$. For any $\epsilon > 0$, there exists an integer $n = n\LRp{\epsilon}$ such that for all $i,j > n$ we have
    \[
     \begin{array}{cr}
     \alpha\nor{\x_i-\x_j}_\Xs \le 
     \nor{\A\x_i - \A\x_j}_\Ys
     \nor{\y_i - \y_j}_\Ys < \epsilon  & \\ 
     \Downarrow &   \\
     \LRc{\x_i}_{i=1}^\infty \text{ is Cauchy } \implies \x_i \xrightarrow[]{\Xs} \x & \\
     \Downarrow &   \text{ continuity of } \A\\
     \y \xleftarrow[]{\Ys} \y_i = \A\x_i \xrightarrow[]{\Ys} \A\x,
     \end{array}
     \]
     and thus $\Range\LRp{\A}$ is closed.

     For the sufficiency, $\Range\LRp{\A}$ is a Banach space due to its closedness. $\A: \Xs \to \Range\LRp{\A}$ is thus bijective, which in turn implies $\A^{-1}$ is linear and  continuous  owing to the Open Mapping Theorem \cite{brezis2010functional,OdenDemkowicz10,ArbogastBona08,Rudin73}. Thus, let $\y = \A\x$, there exists $\beta > 0$ such that
     $\nor{\A^{-1}\y}_\Xs \le \beta \nor{\y}_\Ys \implies \nor{\x}_\Xs \le \beta \nor{\A\x}_\Ys$, and hence $\A$ is bounded below.
     \label{proof:boundednessBelow}
\end{proof}


\bibliographystyle{siamplain}
\bibliography{referencesNew,references}

\end{document}